\documentclass[12pt,reqno]{amsart}
\usepackage{amsmath}
\usepackage{amssymb}
\usepackage{verbatim}
\usepackage{graphicx}
\usepackage[font=footnotesize]{caption}
\usepackage{subcaption}
\usepackage{epstopdf}
\usepackage{multirow}
\usepackage{url}
\usepackage{makecell}
\usepackage{tabularx}
\usepackage{longtable}
\usepackage{graphicx}
\usepackage{tikz}
\usepackage{pgfplots}
\usepackage{xcolor}
\usepackage{pgfplots}
\usepackage{pgfplotstable}
\allowdisplaybreaks

\usepackage{bm}
\usepackage[toc]{appendix}
\usepackage[capitalise]{cleveref}
\usepackage{color}
\allowdisplaybreaks

\usepackage[top=0.5in,bottom=0.5in,left=0.7in,right=0.7in]{geometry}

\newtheorem{lemma}{Lemma}

\newtheorem{theorem}[lemma]{Theorem}

\newtheorem{example}[lemma]{Example}

\newcommand{\N}{\mathbb{N}}    
\newcommand{\NN}{\mathbb{N}_0} 
\newcommand{\R}{\mathbb{R}}    

\newcommand{\nv}{\vec{n}}
\newcommand{\bo}{\mathcal{O}}
\newcommand{\Op}{\Omega^+}
\newcommand{\Om}{\Omega^-}

\newcommand{\ka}{\textsf{k}}
\newcommand{\ia}{\textsf{i}}

\newcommand{\B}{\mathcal{B}}
\newcommand{\rr}{\mathcal{R}} 
\newcommand{\id}{\mathbf{I}_d}
\newcommand{\gd}{g_D}
\newcommand{\gn}{g_N}

\newcommand{\be}{ \begin{equation} }
	\newcommand{\ee}{ \end{equation} }

\newcommand{\ind}{\Lambda}

\numberwithin{equation}{section}
\numberwithin{lemma}{section}

\begin{document}
	
	\title[
Hybrid Finite Difference Scheme
for Elliptic Interface Problems]{
Hybrid Finite Difference Scheme
for Elliptic Interface Problems with Discontinuous and High-Contrast Variable Coefficients}
	
	\author{Qiwei Feng, Bin Han, and Peter Minev}
	
	\thanks{Research supported in part by Natural Sciences and Engineering Research Council (NSERC) of Canada under grants RGPIN-2019-04276 (Bin Han), RGPIN-2017-04152 (Peter Minev),  Westgrid (www.westgrid.ca), and Compute Canada Calcul Canada (www.computecanada.ca)}
	
	\address{Department of Mathematical and Statistical Sciences,
		University of Alberta, Edmonton, Alberta, Canada T6G 2G1.
		\quad {\tt qfeng@ualberta.ca}
		\quad {\tt bhan@ualberta.ca}
		\quad {\tt minev@ualberta.ca}
	}
	
	\makeatletter \@addtoreset{equation}{section} \makeatother
	
	\begin{abstract}
For elliptic interface problems with discontinuous coefficients, the maximum accuracy order for compact 9-point finite difference scheme in irregular points is three \cite{FHM21b}.
The discontinuous coefficients usually have abrupt jumps across the interface curve in the porous medium of realistic problems,
causing the pollution effect of numerical methods.
So, to obtain a reasonable numerical solution of the above problem, the higher order scheme and its effective implementation are necessary.
In this paper, we propose an efficient and flexible way to achieve the implementation of a hybrid (9-point scheme with sixth order accuracy for interior regular points and 13-point scheme with fifth order accuracy for interior irregular points) finite difference scheme in uniform meshes for the elliptic interface problems with discontinuous and high-contrast piecewise smooth coefficients in a rectangle $\Omega$.
We also derive the $6$-point and $4$-point finite difference schemes in uniform meshes with sixth order accuracy for the side points and corner points of various mixed boundary conditions (Dirichlet, Neumann and Robin) of elliptic equations in a rectangle. Our numerical experiments confirm the flexibility and the sixth order accuracy in $l_2$ and $l_{\infty}$ norms of the proposed hybrid scheme.

	\end{abstract}
	
	
	\maketitle
	
	\pagenumbering{arabic}
	
	\section{Introduction and motivations}

Elliptic interface problems with discontinuous coefficients appear in many real-world applications: composite materials, fluid mechanics, nuclear waste disposal, and many others.
One possible avenue to solve such problems, the so-called  immersed interface method (IIM), is proposed by LeVeque and Li.
It has been combined with finite difference, finite volume, and finite element spatial discretizations, with various degree of accuracy.
Some of the most important developments include: the second order IIM \cite{CFL19,LeLi94},
the second order immersed finite volume element methods \cite{EwingLLL99},
the second order immersed finite  element methods \cite{GongLiLi08,HeLL2011},
the second order  fast iterative immersed interface methods of \cite{Li98}, the second order explicit-jump immersed interface methods of \cite{WieBube00},
the third order compact finite difference scheme of \cite{PanHeLi21} and
fourth order IIM of \cite{XiaolinZhong07}.

Another possible approach for the resolution of elliptic interface problems  with discontinuous coefficients is the matched interface and boundary methods (MIB) . The
related papers of MIB for the elliptic interface problems can be summarized as:  second order MIB  \cite{YuZhouWei07},
fourth order MIB  \cite{ZW06},  fourth order MIB with the FFT acceleration  \cite{FendZhao20pp109677}, sixth order  MIB  \cite{YuWei073D,ZZFW06}.
For the anisotropic elliptic interface problems with discontinuous and matrix coefficients, \cite{DFL20} proposed a new finite
element-finite difference (FE-FD) method with a second order of accuracy.

\begin{figure}[h]
	\centering
	\hspace{12mm}	
	\begin{subfigure}[b]{0.4\textwidth}
		\begin{tikzpicture}[scale = 1.5]
			\draw[help lines,step = 1]
			(0,0) grid (4,4);
			\node at (1,1)[circle,fill,inner sep=2.5pt,color=black]{};
			\node at (1,2)[circle,fill,inner sep=2.5pt,color=black]{};
			\node at (1,3)[circle,fill,inner sep=2.5pt,color=black]{};	
			\node at (2,1)[circle,fill,inner sep=2.5pt,color=black]{};
			\node at (2,2)[circle,fill,inner sep=2.5pt,color=black]{};
			\node at (2,3)[circle,fill,inner sep=2.5pt,color=black]{};	
			\node at (3,1)[circle,fill,inner sep=2.5pt,color=black]{};
			\node at (3,2)[circle,fill,inner sep=2.5pt,color=black]{};
			\node at (3,3)[circle,fill,inner sep=2.5pt,color=black]{};	
			\draw[line width=1.5pt, red]  plot [smooth,tension=0.8]
			coordinates {(0,1.1) (1,1.4) (2,2.4) (2.7,4)};
		\end{tikzpicture}
	\end{subfigure}
	\begin{subfigure}[b]{0.4\textwidth}
		\begin{tikzpicture}[scale = 1.5]
			\draw[help lines,step = 1]
			(0,0) grid (4,4);
			\node at (1,1)[circle,fill,inner sep=2.5pt,color=black]{};
			\node at (1,2)[circle,fill,inner sep=2.5pt,color=black]{};
			\node at (1,3)[circle,fill,inner sep=2.5pt,color=black]{};	
			\node at (2,1)[circle,fill,inner sep=2.5pt,color=black]{};
			\node at (2,2)[circle,fill,inner sep=2.5pt,color=black]{};
			\node at (2,3)[circle,fill,inner sep=2.5pt,color=black]{};	
			\node at (3,1)[circle,fill,inner sep=2.5pt,color=black]{};
			\node at (3,2)[circle,fill,inner sep=2.5pt,color=black]{};
			\node at (3,3)[circle,fill,inner sep=2.5pt,color=black]{};	
			\node at (0,2)[circle,fill,inner sep=2.5pt,color=black]{};
			\node at (4,2)[circle,fill,inner sep=2.5pt,color=black]{};
			\node at (2,0)[circle,fill,inner sep=2.5pt,color=black]{};
			\node at (2,4)[circle,fill,inner sep=2.5pt,color=black]{};		
			\draw[line width=1.5pt, red]  plot [smooth,tension=0.8]
			coordinates {(0,1.1) (1,1.4) (2,2.4) (2.7,4)};
		\end{tikzpicture}
	\end{subfigure}
	\caption
	{For irregular points, the 9-point scheme (left) and the 13-point scheme (right). The curve in red color is the interface curve $\Gamma_I$.}
	\label{Extend the compact}
\end{figure}
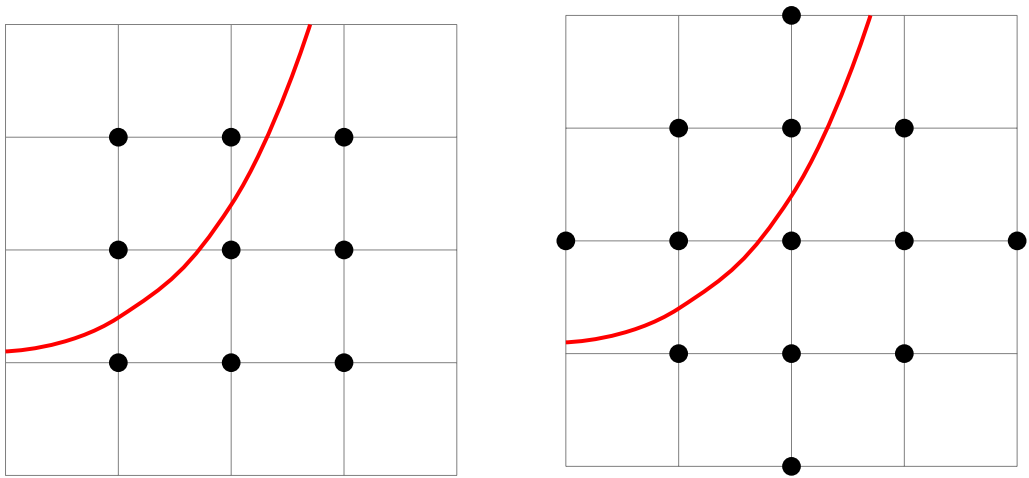

In \cite{FHM21b} we developed a  compact 9-point finite difference scheme for elliptic problems, that is formally fourth order accurate away from the interface of singularity of the solution (regular points), and third order accurate in the vicinity of this interface (irregular points).  The numerical experiments in \cite{FHM21b} demonstrate  that the proposed scheme is fourth order accuracy in the $l_2$ norm. Since the maximum accuracy  for compact 9-point finite difference stencil at regular points is six, and a 13-point  stencil at irregular points can achieve a fifth order of accuracy, in the present paper we derive a hybrid scheme that utilizes a 9-point stencil for  regular points and a 13-point stencil for irregular points, for the case of elliptic problems with discontinuous scalar coefficients.
In  \cite{FHM21b} we demonstrated that if the coefficient of the problem is continuous the stencil of a 9-point scheme in 2D can be partitioned into 72 different configurations by the interface of singularity of the solution.  In the case of discontinuous coefficients, we need to use a 13-point stencil at irregular points and this results in more possibilities for the stencil partitioning (see figure \ref{Extend the compact}). Thus, in the present paper we also derive an efficient way to achieve the implementation of the proposed hybrid scheme.

 A comprehensive literature review of the finite difference approximation of mixed boundary conditions in rectangular domains can be found in \cite{LiPan2021}. In addition, one should also mention the following literature concerned with the discretization of the boundary conditions for elliptic problems:
the sixth order 6-point finite difference scheme for 1-side Neumann and 3-side Dirichlet boundary conditions of Helmholtz equations with constant wave numbers  \cite{Nabavi07},
the sixth order 5-point or 6-point finite difference schemes for 1-side Neumann/Robin and 3-side Dirichlet boundary conditions of Helmholtz equations with variable wave numbers  \cite{TGGT13},
the fourth order MIB for  4-side Robin boundary conditions of  elliptic interface problems \cite{FendZhao20pp109677},
up to 8th order MIB for mixed boundary conditions of Dirichlet, Neumann and Robin with all constant coefficients of  Poisson/Helmholtz equations  \cite{FendZhao20pp109391}.

Compact finite differences have also been successfully applied to elliptic problems with various boundary conditions in non-rectangular domains.
In \cite{RenFengZhao2022} a fourth order MIB for Dirichlet, Neumann, and Robin boundary conditions has been proposed.
  \cite{WieBube00} developed a
second order explicit-jump immersed interface method for problems with Dirichlet and Neumann boundary conditions, and \cite{ItoLiKyei05,LiIto06}
proposed fourth order compact finite difference schemes for various combinations of boundary conditions .

\begin{figure}[htbp]			
	\begin{tikzpicture}
	\draw[domain =0:360,smooth]
	plot({sqrt(0)*cos(\x)}, {sqrt(0)*sin(\x)});
	\draw
	(-pi, -pi) -- (-pi, pi) -- (pi, pi) -- (pi, -pi) --(-pi,-pi);
	\node (A) at (-2.8,0) {$\Gamma_1$};
	\node (A) at (2.8,0) {$\Gamma_2$};
	\node (A) at (0,-2.8) {$\Gamma_3$};
	\node (A) at (0,2.8) {$\Gamma_4$};
	\node (A) at (-5,0) {$\B_1u = \tfrac{\partial u}{\partial \nv}+\alpha u=g_1$};
	\node (A) at (4.4,0) {$\B_2u=u=g_2$};
	\node (A) at (0,-3.5) {$\B_3u=\tfrac{\partial u}{\partial \nv}=g_3$};
	\node (A) at (0,3.5) {$\B_4u = \tfrac{\partial u}{\partial \nv}+\beta u=g_4$};
\end{tikzpicture}	
		\caption
		{An example for the boundary configuration in \eqref{Qeques2}, where $\alpha$ and $\beta$ are two smooth 1D functions in $y$ and $x$ directions respectively.}
		\label{fig:boundary}
	\end{figure}
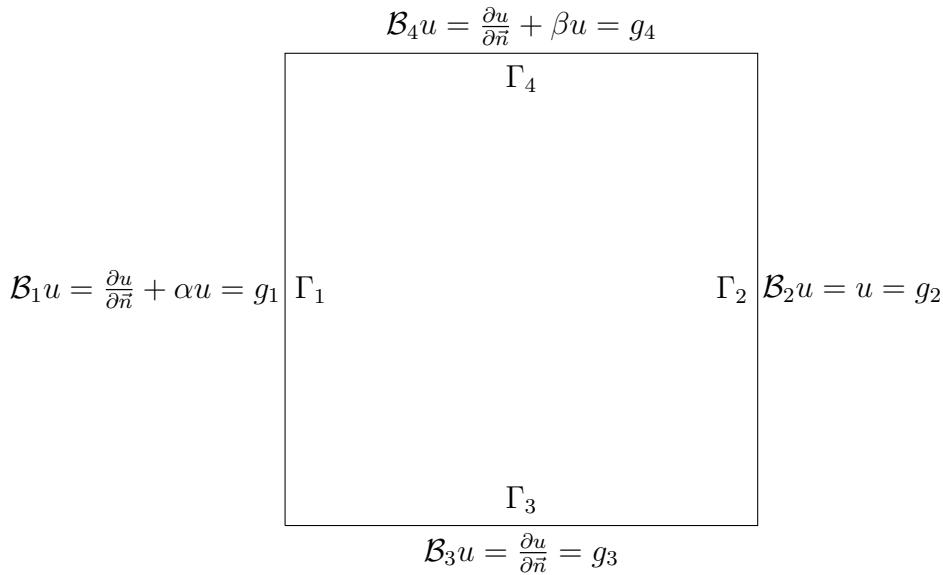

In \cite{FHM21Helmholtz}, we discussed the $6$-point and $4$-point finite difference schemes with sixth order accuracy for the side points and corner points of the Helmholtz equations respectively with a constant wave number $\ka$ in a rectangle. In this paper, we also extend the above results in \cite{FHM21Helmholtz} to the elliptic equations with variable coefficients and mixed combinations  of Dirichlet $u|_{\Gamma_i}=g_i$,  Neumann $\tfrac{\partial u}{\partial \nv}|_{\Gamma_j}=g_j$ and Robin $\tfrac{\partial u}{\partial \nv}+\alpha u|_{\Gamma_k}=g_k$ with  smooth functions $\alpha$, $g_i$, $g_j$ and $g_k$, where  $\Gamma_i/\Gamma_j/\Gamma_k$ for $i,j,k=1,2,3,4$ is one side of the rectangle (see \cref{fig:boundary} for an example of the mixed boundary conditions).

In order to define the subject of the present paper, let $\Omega=(l_1, l_2)\times(l_3, l_4)$ and $\psi$ be a smooth two-dimensional function. Consider a smooth curve $\Gamma_I:=\{(x,y)\in \Omega: \psi(x,y)=0\}$, which partitions $\Omega$ into two subregions:
	$\Op:=\{(x,y)\in \Omega\; :\; \psi(x,y)>0\}$ and $\Om:=\{(x,y)\in \Omega\; : \; \psi(x,y)<0\}$. We also define
	$
	a_{\pm}:=a \chi_{\Omega^{\pm}}$, $f_{\pm}:=f \chi_{\Omega^{\pm}}$ and $u_{\pm}:=u \chi_{\Omega^{\pm}}.
	$ The model problem in this paper is defined as follows:
	\begin{equation} \label{Qeques2}
		\begin{cases}
		-\nabla \cdot( a\nabla u)=f &\text{in $\Omega \setminus \Gamma_I$},\\	
			\left[u\right]=\gd, \quad \left[a\nabla  u \cdot \nv \right]=\gn &\text{on $\Gamma_I$},\\
			\B_1 u =g_1 \text{ on } \Gamma_1:=\{l_{1}\} \times (l_3,l_4), & \B_2 u =g_2 \text{ on } \Gamma_2:=\{l_{2}\} \times (l_3,l_4),\\
			\B_3 u =g_3 \text{ on } \Gamma_3:=(l_1,l_2) \times \{l_{3}\}, &
			\B_4 u =g_4 \text{ on } \Gamma_4:=(l_1,l_2) \times \{l_{4}\},
		\end{cases}
	\end{equation}
	where  $f$ is the source term, and for any point $(x_0,y_0)\in \Gamma_I$,
	\begin{align*}
		[u](x_0,y_0) & :=\lim_{(x,y)\in \Op, (x,y) \to (x_0,y_0) }u(x,y)- \lim_{(x,y)\in \Om, (x,y) \to (x_0,y_0) }u(x,y),\\
		[ a\nabla  u \cdot \nv](x_0,y_0) & :=  \lim_{(x,y)\in \Op, (x,y) \to (x_0,y_0) } a\nabla  u(x,y) \cdot \nv- \lim_{(x,y)\in \Om, (x,y) \to (x_0,y_0) } a \nabla  u(x,y) \cdot \nv,
	\end{align*}
where $\nv$ is the unit normal vector of $\Gamma_I$ pointing towards $\Op$.
	In \eqref{Qeques2}, the boundary operators $\B_1,\ldots,\B_4  \in \{\id,\frac{\partial }{\partial \nv}+ \alpha \id\}$, where
$\id$ represents the Dirichlet boundary condition, when $\alpha=0$, $\frac{\partial }{\partial \nv}$ represents the Neumann boundary condition,  when $\alpha$ is a smooth 1D function, $\frac{\partial }{\partial \nv}+\alpha \id$ represents the Robin boundary condition. An example for the boundary conditions of  \eqref{Qeques2} is shown in \cref{fig:boundary}.

We derive a hybrid   finite difference scheme to solve \eqref{Qeques2} given the following assumptions:
	\begin{itemize}
		
	\item[(A1)] The coefficient	$a$ is positive, piecewise smooth  and has uniformly continuous partial derivatives of (total) orders up to six  in each of the subregions $\Op$ and $\Om$. The coefficient $a$ is discontinuous across the interface $\Gamma_I$.
		
		\item[(A2)] The solution $u$ and the source term $f$ have uniformly continuous partial derivatives of (total) orders up to seven and five respectively in each of the subregions $\Op$ and $\Om$. Both $u$ and $f$ can be discontinuous across the interface $\Gamma_I$.

		\item[(A3)] The interface curve $\Gamma_I$ is smooth in the sense that for each $(x^*,y^*)\in \Gamma_I$, there exists a local parametric equation: $\gamma: (-\epsilon,\epsilon)\rightarrow \Gamma_I$ with $\epsilon>0$ such that $\gamma(0)=(x^*,y^*)$ and $\|\gamma'(0)\|_{2}\ne 0$.
		\item[(A4)] The 1D interface functions $\gd \circ \gamma$ and $\gn \circ \gamma$ have uniformly continuous derivatives of (total) orders up to five and four respectively on the interface $\Gamma_I$, where $\gamma$ is given in (A2).
		\item[(A5)] Each of the 1D boundary functions $g_1,\ldots,g_4$ in \eqref{Qeques2} and $\alpha$ in the Robin boundary conditions
has uniformly continuous derivatives of (total) order up to five on the boundary $\Gamma_j$.
	\end{itemize}

	The organization of this paper is as follows.

	In \cref{subsec:regular},	we derive the compact 9-point finite difference scheme with sixth order accuracy for regular points in \cref{thm:regular:interior}.
	
	In \cref{Boundarypoints},	we propose the $6$-point schemes with sixth order accuracy for  the side points of the boundary conditions  $\tfrac{\partial u}{\partial \nv}+\alpha u|_{\Gamma_1}=g_1$, $\tfrac{\partial u}{\partial \nv}|_{\Gamma_3}=g_3$ and $\tfrac{\partial u}{\partial \nv}+\beta u|_{\Gamma_4}=g_4$  in \cref{hybrid:thm:regular:Robin:1,hybrid:thm:regular:Neu:3,hybrid:thm:regular:Robin:4}
	with two smooth functions $\alpha$ and $\beta$.
	
	In \cref{Cornerpoints},	we construct the $4$-point schemes with sixth order accuracy for  the corner points of the boundary conditions  $\tfrac{\partial u}{\partial \nv}+\alpha u|_{\Gamma_1}=g_1$, $\tfrac{\partial u}{\partial \nv}|_{\Gamma_3}=g_3$ and $\tfrac{\partial u}{\partial \nv}+\beta u|_{\Gamma_4}=g_4$  in \cref{thm:corner:1,thm:corner:2}
with two smooth functions $\alpha$ and $\beta$.

In \cref{hybrid:Irregular:points}, we first propose
a simpler version of the transmission equation for the interface curve $\Gamma_I$ in \cref{thm:interface}. Then the 13-point finite difference scheme with fifth order accuracy  for irregular points is shown in \cref{13point:Inter}.
In order to achieve the implementation effectively for the 13-point scheme, we derive efficient implementation details using \eqref{ChandX} to \eqref{ckln:irregular}.

	 In \cref{sec:numerical}, we present $10$ numerical examples, including $5$ examples with exact known solutions $u$,
for our proposed hybrid finite difference scheme with contrast ratios  $\sup(a_+)/\inf(a_-)=10^{-3},10^{-6},10^{6},10^{7}$.  Our numerical experiments confirm the flexibility and the sixth order accuracy in $l_2$ and $l_{\infty}$ norms of our proposed hybrid scheme.	 For the coefficients $a(x,y)$, two jump functions $g_D,g_N$, interface curves $\Gamma_I$ and  boundary conditions, we test the following cases:

\begin{itemize}	
\item Either $a_+/a_-$ or $a_-/a_+$ is very large on the interface $\Gamma_I$ for high contrast coefficients $a$.
\item The jump functions $g_D$ and $g_N$ are both either constant or non-constant.
\item The interface curve $\Gamma_I$ is either smooth or sharp-edged.
	\item 4-side Dirichlet boundary conditions.
	\item 3-side Dirichlet and 1-side Robin boundary conditions.
	\item 1-side Dirichlet, 1-side Neumann and 2-side Robin boundary conditions.
\end{itemize}
	
In \cref{sec:hybrid:Conclu}, we summarize the main contributions of this paper. Finally, in \cref{hybrid:sec:proofs} we present the proofs for results stated in \cref{sec:sixord}.

	\section{Hybrid finite difference method on uniform Cartesian grids}
	\label{sec:sixord}

	We follow the same setup as in \cite{FHM21a,FHM21b,FHM21Helmholtz}. Let $\Omega=(l_1,l_2)\times (l_3,l_4)$ and we assume $l_4-l_3=N_0 (l_2-l_1)$ for some $N_0 \in \N$. For any positive integer $N_1\in \N$, we define $N_2:=N_0 N_1$ and so the grid size is  $h:=(l_2-l_1)/N_1=(l_4-l_3)/N_2$.
	Let
	\be \label{xiyj}
	x_i=l_1+i h, \quad i=0,\ldots,N_1, \quad \text{and} \quad y_j=l_3+j h, \quad j=0,\ldots,N_2.
	\ee
 Recall that a compact stencil centered at $(x_i,y_j)$ contains nine points $(x_i+kh, y_j+lh)$ for $k,l\in \{-1,0,1\}$. Define
 \be\label{Compact:Set}
	\begin{split}
		& d_{i,j}^+:=\{(k,\ell) \; : \;
		k,\ell\in \{-1,0,1\}, \psi(x_i+kh, y_j+\ell h)\ge 0\}, \quad \mbox{and}\\
		& d_{i,j}^-:=\{(k,\ell) \; : \;
		k,\ell\in \{-1,0,1\}, \psi(x_i+kh, y_j+\ell h)<0\}.
	\end{split}
\ee
	Thus, the interface curve $\Gamma_I:=\{(x,y)\in \Omega \; :\; \psi(x,y)=0\}$ splits the nine points in our compact stencil into two disjoint sets  $\{(x_{i+k}, y_{j+\ell})\; : \; (k,\ell)\in d_{i,j}^+\} \subseteq \Op \cup \Gamma_I$ and
	$\{(x_{i+k}, y_{j+\ell})\; : \; (k,\ell)\in d_{i,j}^-\} \subseteq \Om$. We refer to a grid/center point $(x_i,y_j)$ as
	\emph{a regular point} if  $d_{i,j}^+=\emptyset$ or $d_{i,j}^-=\emptyset$.
	The center point $(x_i,y_j)$ of a stencil is \emph{regular} if all of its nine points are in $\Op \cup \Gamma_I$ (hence $d_{i,j}^-=\emptyset$) or in $\Om$ (i.e., $d_{i,j}^+=\emptyset$).
	Otherwise, if both $d_{i,j}^+$ and $d_{i,j}^-$ are nonempty, the center point $(x_i,y_j)$ of a stencil is referred to as \emph{an irregular point} .
	
	Now, let us pick and fix a base point $(x_i^*,y_j^*)$ inside the open square $(x_i-h,x_i+h)\times (y_j-h,y_j+h)$, which can be written as
	\be \label{base:pt}
	x_i^*=x_i-v_0h  \quad \mbox{and}\quad y_j^*=y_j-w_0h  \quad \mbox{with}\quad
	-1<v_0, w_0<1.
	\ee
Throughout the paper, we shall use the following notations:
	\be\label{ufmn}
		\begin{split}
		&\alpha^{(n)}:=\frac{d^{n} \alpha}{ dy^n }(y_j^*),
		\quad {g_1}^{(n)}:=\frac{d^{n} g_1}{ dy^n }(y_j^*), \\
		& \beta^{(m)}:=\frac{d^{m} \beta}{ dx^m }(x_i^*),
		\quad {g_3}^{(m)}:=\frac{d^{m} g_3}{ dx^m }(x_i^*), 	\quad {g_4}^{(m)}:=\frac{d^{m} g_4}{ dx^m }(x_i^*),\\
		&a^{(m,n)}:=\frac{\partial^{m+n} a}{ \partial^m x \partial^n y}(x_i^*,y_j^*),
		\quad u^{(m,n)}:=\frac{\partial^{m+n} u}{ \partial^m x \partial^n y}(x_i^*,y_j^*),\quad
		f^{(m,n)}:=\frac{\partial^{m+n} f}{ \partial^m x \partial^n y}(x_i^*,y_j^*),
	\end{split}
	\ee
	which are their $(m,n)$th partial derivatives at the base point $(x_i^*,y_j^*)$.
By \cite[(2.13)]{FHM21b}, we have
\be \label{u:approx:key:V}
u(x+x_i^*,y+y_j^*)
=
\sum_{(m,n)\in \ind_{M+1}^{V,1}}
u^{(m,n)} G^V_{M,m,n}(x,y) +\sum_{(m,n)\in \ind_{M-1}}
f^{(m,n)} Q^V_{M,m,n}(x,y)+\bo(h^{M+2}),
\ee
for $x,y\in (-2h,2h)$, where  $u$ is the exact solution for \eqref{Qeques2},
the index sets $\ind_{M-1}$ and $\ind_{M+1}^{V,1}$ are defined in \eqref{Sk} and \eqref{indV12} respectively, and
the functions $G^V_{M,m,n}$ and $Q^V_{M,m,n}$ are defined in \eqref{GVmn} and \eqref{QVmn} respectively.
By \cite[(2.13) and (2.14)]{FHM21Helmholtz}, we also have
\be
\label{u:approx:key:H}
u(x+x_i^*,y+y_j^*)  = \sum_{(m,n)\in \ind_{M+1}^{H, 1}} u^{(m,n)} G^{H}_{M,m,n}(x,y) + \sum_{(m,n)\in \ind_{M-1}} f^{(m,n)} Q^H_{M,m,n}(x,y) + \mathcal{O}(h^{M+2}),
\ee
where the index sets $\ind_{M-1}$ and $\ind_{M+1}^{H,1}$ are defined in \eqref{Sk} and \eqref{indH12} respectively, and the functions $G^H_{M,m,n}$ and $Q^H_{M,m,n}$ are defined in \eqref{GHmn} and \eqref{QHmn} respectively.

For the sake of better readability, all technical proofs of this section are provided in \cref{hybrid:sec:proofs}.


	\subsection{Stencils for regular points (interior)}
	\label{subsec:regular}	
	We now extend the fourth order compact scheme in \cite[Theorem~3.1]{FHM21b} to a sixth order compact scheme. We only need to choose $M=6$ and replace  $G_{m,n}$, $H_{m,n}$ and $\ind_{M+1}^{1}$ in \cite{FHM21b} by $G^V_{M,m,n}$ in \eqref{GVmn}, $Q^V_{M,m,n}$ in \eqref{QVmn}, and  $\ind_{M+1}^{V,1}$ in \eqref{indV12}. We choose $(x_i^*,y_j^*)$ to be the center point of the 9-point compact scheme, i.e., $(x_i^*,y_j^*)=(x_i,y_j)$ and $v_0=w_0=0$ in \eqref{base:pt}.

	 \begin{theorem}\label{thm:regular:interior}
		Let a grid point $(x_i,y_j)$ be a regular point, i.e., either $d_{i,j}^+=\emptyset$ or $d_{i,j}^-=\emptyset$ and $(x_i,y_j) \notin \partial \Omega $. Let $(u_{h})_{i,j}$ denote the numerical approximation of the exact solution $u$ of the elliptic interface problem \eqref{Qeques2} at an interior regular point $(x_i, y_j)$. Then the following difference scheme on a stencil centered at $(x_{i},y_j)$:
		\be\label{stencil:regular:interior:V}
		\mathcal{L}_h u_h :=
		\begin{aligned}	
			&C_{-1,-1}(u_{h})_{i-1,j-1}& 
			+&C_{0,-1}(u_{h})_{i,j-1}&    
			+&C_{1,-1}(u_{h})_{i+1,j-1} \\
			+&C_{-1,0}(u_{h})_{i-1,j}&
			+&C_{0,0}(u_{h})_{i,j}&
			+&C_{1,0}(u_{h})_{i+1,j}\\
			+&C_{-1,1}(u_{h})_{i-1,j+1}&
			+&C_{0,1}(u_{h})_{i,j+1}&
			+&C_{1,1}(u_{h})_{i+1,j+1}\\
		\end{aligned}
		=\sum_{(m,n)\in \ind_{5}} f^{(m,n)}C_{f,m,n},
		\ee
		achieves sixth order of accuracy for $-\nabla \cdot( a\nabla u)=f$ at the point $(x_i,y_j)$, where
\[
C_{f,m,n}:= \sum_{k=-1}^1\sum_{\ell=-1}^1 C_{k,\ell} Q^{V}_{6,m,n}(kh, \ell h), \quad \mbox{for all} \quad (m,n)\in \ind_{5},
\]
\begin{equation}\label{Ckell}
C_{k,\ell}(h):=\sum_{i=0}^{M+1} c_{k,\ell,i}h^i,
\qquad k,\ell \in \{-1,0,1\},
\end{equation}
and $\{c_{k,\ell,i}\}$ is any non-trivial solution to the linear system induced by \cite[(3.5)]{FHM21b} with $M=6$.
		Moreover, the maximum accuracy order of a compact finite difference scheme for  $-\nabla \cdot( a\nabla u)=f$ at the point $(x_i,y_j)$ is six.
	\end{theorem}

To verify \cref{thm:regular:interior} with the numerical experiments in \cref{sec:numerical}, we use the unique solution $\{c_{k,\ell,i}\}$ to \cite[(3.5)]{FHM21b} with $M=6$ and the normalization condition $c_{-1,-1,0}=1$, setting to zero all $c_{-1,0,7}, c_{0,-1,7},
c_{0,0,6}, c_{0,0,7},
c_{-1,1,i_1}, c_{0,1,i_2}, c_{1,-1,i_2}, c_{1,0,i_3}, c_{1,1,i_4}$ for
$i_1=1,6,7$, $i_2=5,6,7$, $i_3=4,5,6,7$ and $i_4=2,3,4,5,6,7$.


\subsection{Stencils for boundary points}\label{All:Boundary}
In this subsection, we extend \cite[Section~2.2]{FHM21Helmholtz} and discuss how to find a compact finite difference scheme with accuracy order six centered at $(x_i,y_j) \in \partial \Omega$.
For clarity of presentation, we consider the following boundary conditions
%
\be\label{model:corner}
\begin{aligned}
	&\B_1u=\tfrac{\partial u}{\partial \nv}+\alpha u=g_1 \;\; \text{on} \;\; \Gamma_1,\qquad
	&& \B_2u=u=g_2 \;\; \text{on} \;\; \Gamma_2,\\
	&\B_3u=\tfrac{\partial u}{\partial \nv}=g_3 \;\; \text{on} \;\; \Gamma_3,
	&&\B_4u=\tfrac{\partial u}{\partial \nv}+\beta u=g_4 \;\; \text{on} \;\; \Gamma_4,
\end{aligned}
\ee
where $\alpha$ and $\beta$ are two smooth 1D functions in $y$ and $x$ directions.
For the 6-point and 4-point schemes in this subsection, we choose $(x_i^*,y_j^*)=(x_i,y_j)$ and $v_0=w_0=0$ in \eqref{base:pt}.
An illustration of \eqref{model:corner} is shown in \cref{fig:boundary}. For the following identities in \eqref{CB1:EQ} and \eqref{CB4:EQ}, we define
\be\label{delta:Fun}
\delta_{a,a}:=1 \quad \mbox{and}\quad
\delta_{a,b}:=0 \quad \mbox{for } a\ne b.
\ee

\subsubsection{Side points on the boundary $\partial \Omega$}\label{Boundarypoints}

\begin{theorem}\label{hybrid:thm:regular:Robin:1}
	Let  $(u_{h})_{i,j}$ denote the numerical approximation of the exact solution $u$ of the elliptic interface problem \eqref{Qeques2} at the point $(x_i, y_j)$. The following discretization on a stencil  centered at $(x_0,y_j)\in \Gamma_1$:
	\be	 \label{stencil:regular:Robin:1}
	\mathcal{L}^{\mathcal{B}_1}_h u_h  :=
	\begin{aligned}	
		&C^{\mathcal{B}_1}_{0,-1}(u_{h})_{0,j-1}&    
		+&C^{\mathcal{B}_1}_{1,-1}(u_{h})_{1,j-1} \\		 +&C^{\mathcal{B}_1}_{0,0}(u_{h})_{0,j}&             +&C^{\mathcal{B}_1}_{1,0}(u_{h})_{1,j}\\
		+&C^{\mathcal{B}_1}_{0,1}(u_{h})_{0,j+1}&           +&C^{\mathcal{B}_1}_{1,1}(u_{h})_{1,j+1}\\
	\end{aligned} =	\sum_{(m,n)\in \ind_{4}} f^{(m,n)}C^{\mathcal{B}_1}_{f,m,n}+\sum_{n=0}^{5}g_{1}^{(n)}C_{g_{1},n}^{\B_1},
	\ee
	achieves sixth order of accuracy for $\B_1u=\frac{\partial u}{\partial \nv}+\alpha u=g_1$ at the point $(x_0,y_j) \in \Gamma_1$,
	where
	 \[
	C^{\mathcal{B}_1}_{f,m,n} := \sum\limits_{k=0}^1 \sum\limits_{\ell=-1}^1 C^{\mathcal{B}_1}_{k,\ell} Q^{V}_{5,m,n}(kh, \ell h),\quad \mbox{for all} \quad (m,n) \in \ind_4,
	\]
	\[
	C_{g_1,n}^{\B_1} :=  -\sum\limits_{k=0}^1 \sum\limits_{\ell=-1}^1 C^{\mathcal{B}_1}_{k,\ell} G^{V}_{5,1,n}(kh, \ell h), \quad \mbox{for all} \quad n=0, \dots, 5,
	\]
\[
	C^{\B_1}_{k,\ell}(h):=\sum_{i=0}^{6} c^{\B_1}_{k,\ell,i}h^i,
	\qquad k\in\{0,1\},\ell\in\{-1,0,1\},
\]
	and $\{c^{\B_1}_{k,\ell,i}\}$ is any non-trivial solution to the linear system induced by
	\be
	\begin{split}\label{CB1:EQ}
	&\sum_{k=0}^1 \sum_{\ell=-1}^1 C^{\B_1}_{k,\ell} \left( G^{V}_{5,0,n}(kh, \ell h) +
	\sum_{i=n}^5  {i\choose n}  {\alpha}^{(i-n)} G^{V}_{5,1,i}(kh, \ell h) (1-\delta_{n,6}) \right)\\
	&=\bo(h^{7}),\quad \mbox{for all } \quad n=0,1,\dots,6.
	\end{split}
    \ee
%
%
  Moreover, the maximum accuracy order of a 6-point finite difference scheme for $\B_1u=\frac{\partial u}{\partial \nv} +\alpha u=g_1$  at the point $(x_0,y_j) \in \Gamma_1$ with two smooth functions $\alpha(y)$ and $a(x,y)$   is six.
\end{theorem}

In our numerical experiments in \cref{sec:numerical}, we use the unique solution $\{c^{\B_1}_{k,\ell,i}\}$ to \eqref{CB1:EQ} with the normalization condition
$c^{\B_1}_{1,1,0}=1$, where 
all $c^{\B_1}_{0,0,6},  c^{\B_1}_{0,1,5},
c^{\B_1}_{0,1,6},
c^{\B_1}_{1,-1,i_1}, c^{\B_1}_{1,0,i_2}, c^{\B_1}_{1,1,i_3}$ for $i_1=1,4,5,6$, $i_2=3,4,5,6$, and $i_3=2,3,4,5,6$, are set to zero.
In particular, if $a$ in \eqref{Qeques2} is a discontinuous constant coefficient and $\B_1u=\frac{\partial u}{\partial \nv}+\alpha u=g_1$ with a constant $\alpha$, then the coefficients
in \eqref{stencil:regular:Robin:1} are
{\small{
		\be \label{stencil:CB1:R}
		\begin{aligned}
			&C^{\mathcal{B}_1}_{0,1} =\frac{1}{75}{\alpha}^2h^2+\frac{1}{5}{\alpha}h+2,\ \  C^{\mathcal{B}_1}_{0,0} =\frac{8}{675}{\alpha}^5h^5-\frac{16}{675}{\alpha}^4h^4+\frac{16}{225}{\alpha}^3h^3-\frac{8}{25}{\alpha}^2h^2-\frac{34}{5}{\alpha}h-10,\\
			& C^{\mathcal{B}_1}_{1,1} =1,\ \ \  C^{\mathcal{B}_1}_{1,0} = -\frac{8}{675}{\alpha}^4h^4+\frac{8}{225}{\alpha}^3h^3-\frac{8}{75}{\alpha}^2h^2+\frac{2}{5}{\alpha}h+4,\ \ C^{\mathcal{B}_1}_{0,-1}=C^{\mathcal{B}_1}_{0,1},\ \  C^{\mathcal{B}_1}_{1,-1}=C^{\mathcal{B}_1}_{1,1}.
		\end{aligned}
		\ee
}}
Similarly, we could obtain the following \cref{hybrid:thm:regular:Neu:3,hybrid:thm:regular:Robin:4}.

\begin{theorem}\label{hybrid:thm:regular:Neu:3}
	Let  $(u_{h})_{i,j}$ be the numerical approximation of the exact solution $u$ of the elliptic interface problem \eqref{Qeques2} at the point $(x_i, y_j)$. Then the following discretization stencil  centered at $(x_i,y_0)\in \Gamma_3$:
	\be	 \label{stencil:regular:H:Neu:3}
	\mathcal{L}^{\mathcal{B}_3}_h u_h  :=
	\begin{aligned}	
		&C^{\mathcal{B}_3}_{-1,0}(u_{h})_{i-1,0}&             + &C^{\mathcal{B}_3}_{0,0}(u_{h})_{i,0}&             +&C^{\mathcal{B}_3}_{1,0}(u_{h})_{i+1,0}\\
	+	&C^{\mathcal{B}_3}_{-1,1}(u_{h})_{i-1,1}& 	 +&C^{\mathcal{B}_3}_{0,1}(u_{h})_{i,1}&           +&C^{\mathcal{B}_3}_{1,1}(u_{h})_{i+1,1}\\
	\end{aligned} =	\sum_{(m,n)\in \ind_{4}} f^{(m,n)}C^{\mathcal{B}_3}_{f,m,n}+\sum_{n=0}^{5}g_{3}^{(n)}C_{g_{3},n}^{\B_3},
	\ee
		achieves sixth order of accuracy for $\B_3u=\frac{\partial u}{\partial \nv}=g_3$ at the point $(x_i,y_0) \in \Gamma_3$, where
			\[
		C^{\mathcal{B}_3}_{f,m,n} := \sum\limits_{k=-1}^1 \sum\limits_{\ell=0}^1 C^{\mathcal{B}_3}_{k,\ell} Q^{H}_{5,m,n}(kh, \ell h), \quad \mbox{for all} \quad (m,n)  \in \ind_4,
		\]
		\[
		C_{g_3,n}^{\B_3} :=  -\sum\limits_{k=-1}^1 \sum\limits_{\ell=0}^1 C^{\mathcal{B}_3}_{k,\ell} G^{H}_{5,n,1}(kh, \ell h), \quad \mbox{for all} \quad n=0, \dots, 5,
		\]
		\[
		C^{\B_3}_{k,\ell}(h):=\sum_{i=0}^{6} c^{\B_3}_{k,\ell,i}h^i,
		\qquad k\in\{-1,0,1\},\ell\in\{0,1\},
		\]
and $\{c^{\B_3}_{k,\ell,i}\}$ is any non-trivial solution to the linear system induced by
\begin{equation}\label{CB3:EQ}
	\sum_{k=-1}^1 \sum_{\ell=0}^1 C^{\B_3}_{k,\ell} G^{H}_{5,n,0}(kh, \ell h)
	=\bo(h^{7}),\quad \mbox{for all } \quad n=0,1,\dots,6,
\end{equation}
Moreover, the maximum accuracy order of a 6-point finite difference scheme for $\B_3u=\frac{\partial u}{\partial \nv} =g_3$  at the point $(x_i,y_0) \in \Gamma_3$ with a smooth function $a(x,y)$ is six.
\end{theorem}

For our numerical experiments in \cref{sec:numerical}, we use the unique solution $\{c^{\B_3}_{k,\ell,i}\}$ to
\eqref{CB3:EQ} with the normalization condition
$c^{\B_3}_{1,1,0}=1$, presetting to zero all
$c^{\B_3}_{0,0,6}, c^{\B_3}_{-1,1,i_1}, c^{\B_3}_{0,1,i_2}, c^{\B_3}_{1,0,i_3}, c^{\B_3}_{1,1,i_4}$ for
$i_1=1,5,6$,  $i_2=4,5,6$, $i_3=3,4,5,6$, and $i_4=2,3,4,5,6$.
In particular, if $a$ is a discontinuous constant coefficient in \eqref{Qeques2}, then the coefficients
in \eqref{stencil:regular:H:Neu:3} are
\be \label{stencil:CB3:N}
C^{\mathcal{B}_3}_{1,0} =2,\ \  C^{\mathcal{B}_3}_{1,1} =1,\ \
C^{\mathcal{B}_3}_{0,0} =-10,\ \  C^{\mathcal{B}_3}_{0,1} =4,\ \ C^{\mathcal{B}_3}_{-1,0}=C^{\mathcal{B}_3}_{1,0},\ \ C^{\mathcal{B}_3}_{-1,1}=C^{\mathcal{B}_3}_{1,1}.
\ee

\begin{theorem}\label{hybrid:thm:regular:Robin:4}
	Let  $(u_{h})_{i,j}$ be the numerical approximation of the exact solution $u$ of the elliptic interface problem \eqref{Qeques2} at the point $(x_i, y_j)$. Then the following discretization stencil  centered at $(x_i,y_{N_2})\in \Gamma_4$:
	\be	 \label{stencil:regular:H:Robin:4}
	\mathcal{L}^{\mathcal{B}_4}_h u_h  :=
	\begin{aligned}	
		 &C^{\mathcal{B}_4}_{-1,-1}(u_{h})_{i-1,-1}& 	 +&C^{\mathcal{B}_4}_{0,-1}(u_{h})_{i,-1}&           +&C^{\mathcal{B}_4}_{1,-1}(u_{h})_{i+1,-1}\\
		 +&C^{\mathcal{B}_4}_{-1,0}(u_{h})_{i-1,0}&             + &C^{\mathcal{B}_4}_{0,0}(u_{h})_{i,0}&             +&C^{\mathcal{B}_4}_{1,0}(u_{h})_{i+1,0}\\	
	\end{aligned} =	\sum_{(m,n)\in \ind_{4}} f^{(m,n)}C^{\mathcal{B}_4}_{f,m,n}+\sum_{n=0}^{5}g_{4}^{(n)}C_{g_{4},n}^{\B_4},
	\ee
	achieves sixth order of accuracy for $\B_4u=\frac{\partial u}{\partial \nv}+\beta u=g_4$ at the point $(x_i,y_{N_2}) \in \Gamma_4$,
	where
	\[
	C^{\mathcal{B}_4}_{f,m,n} := \sum\limits_{k=-1}^1 \sum\limits_{\ell=-1}^0 C^{\mathcal{B}_4}_{k,\ell} Q^{H}_{5,m,n}(kh, \ell h), \quad \mbox{for all} \quad (m,n) \in \ind_4,
	\]
	\[
	C_{g_4,n}^{\B_4} :=  \sum\limits_{k=-1}^1 \sum\limits_{\ell=-1}^0 C^{\mathcal{B}_4}_{k,\ell} G^{H}_{5,n,1}(kh, \ell h), \quad \mbox{for all} \quad n=0, \dots, 5,
	\]
		\[
C^{\B_4}_{k,\ell}(h):=\sum_{i=0}^{6} c^{\B_4}_{k,\ell,i}h^i,
\qquad k\in\{-1,0,1\},\ell\in\{-1,0\},
\]	
and $\{c^{\B_4}_{k,\ell,i}\}$ is any non-trivial solution to the linear system induced by
		\be
	\begin{split}\label{CB4:EQ}
		&\sum_{k=-1}^1 \sum_{\ell=-1}^0 C^{\B_4}_{k,\ell} \left( G^{H}_{5,n,0}(kh, \ell h) -
		\sum_{i=n}^5  {i\choose n}  {\beta}^{(i-n)} G^{H}_{5,i,1}(kh, \ell h) (1-\delta_{n,6}) \right)\\
		&=\bo(h^{7}),\quad \mbox{for all } \quad n=0,1,\dots,6.
	\end{split}
	\ee
	 Moreover, the maximum accuracy order of a 6-point finite difference scheme for $\B_4u=\frac{\partial u}{\partial \nv} +\beta u=g_4$  at the point $(x_i,y_{N_2}) \in \Gamma_4$ with two smooth functions $\beta(x)$ and $a(x,y)$   is six.
\end{theorem}

For our numerical experiments in \cref{sec:numerical}, we use the unique solution $\{c^{\B_4}_{k,\ell,i}\}$ to
\eqref{CB4:EQ} with the normalization condition
$c^{\B_4}_{1,-1,0}=1$, presetting to zero all
$c^{\B_4}_{0,-1,6}, c^{\B_4}_{-1,0,5}, c^{\B_4}_{-1,0,6}, c^{\B_4}_{0,0,i_1}, c^{\B_4}_{1,-1,i_2}, c^{\B_4}_{1,0,i_3}$
with $i_1=4,5,6$, $i_2=2,3,4,5,6$, $i_3=1,3,4,5,6$.
In particular, if $a$ is a discontinuous piecewise constant coefficient  and $\B_4u=\frac{\partial u}{\partial \nv}+\beta u=g_4$ with a constant $\beta$, then the coefficients
in \eqref{stencil:regular:H:Robin:4} are
\be \label{stencil:CB4:R}
\begin{aligned}
	&C^{\mathcal{B}_4}_{1,-1} =1,\ \  C^{\mathcal{B}_4}_{0,-1} =-\frac{8}{675}{\beta}^4h^4+\frac{8}{225}{\beta}^3h^3-\frac{8}{75}{\beta}^2h^2+\frac{2}{5}{\beta}h+4,\\
	& C^{\mathcal{B}_4}_{1,0} =\frac{1}{75}{\beta}^2h^2+\frac{1}{5}{\beta}h+2
	,\ \ \  C^{\mathcal{B}_4}_{0,0} =\frac{8}{675}{\beta}^5h^5-\frac{16}{675}{\beta}^4h^4+\frac{16}{225}{\beta}^3h^3-\frac{8}{25}{\beta}^2h^2-\frac{34}{5}{\beta}h-10,\\
	& C^{\mathcal{B}_4}_{-1,-1}=C^{\mathcal{B}_4}_{1,-1},\ \ \  C^{\mathcal{B}_4}_{-1,0}=C^{\mathcal{B}_4}_{1,0}.
\end{aligned}
\ee

\subsubsection{Stencils for corner points}\label{Cornerpoints}
\begin{theorem}  \label{thm:corner:1}
	Let  $(u_{h})_{i,j}$ be the numerical approximation of the exact solution $u$ of the elliptic interface problem \eqref{Qeques2} at the point $(x_i, y_j)$. Then the following discretization on a stencil centered at the corner point $(x_0,y_0)$:
	\be \label{stencil:corner:1}
			\begin{aligned}
				\mathcal{L}^{\rr_1}_h u_h  :=
				\begin{aligned}	
					 &C^{\mathcal{R}_1}_{0,0}(u_{h})_{0,0}&             +&C^{\mathcal{R}_1}_{1,0}(u_{h})_{1,0}\\
					 +&C^{\mathcal{R}_1}_{0,1}(u_{h})_{0,1}&           +&C^{\mathcal{R}_1}_{1,1}(u_{h})_{1,1}
				\end{aligned}
				= \sum_{(m,n)\in \ind_{4}} f^{(m,n)}C^{\rr_1}_{f,m,n} + \sum_{n=0}^{5}g_{1}^{(n)}C_{g_{1},n}^{\rr_1} + \sum_{n=0}^{5}g_{3}^{(n)}C_{g_{3},n}^{\rr_1},
			\end{aligned}	
	\ee
	achieves sixth order of accuracy for $\B_1u=\frac{\partial u}{\partial \nv}+\alpha u=g_1$ and $\B_3u=\frac{\partial u}{\partial \nv}=g_3$ at the point $(x_0,y_0)$, where $\{C^{\rr_1}_{k,\ell}\}_{k,\ell \in \{0,1\}}$, $\{C^{\rr_1}_{f,m,n}\}_{(m,n) \in \ind_4}$, $\{C^{\rr_1}_{g_1,n}\}_{n=0}^5$ and  $\{C^{\rr_1}_{g_3,n}\}_{n=0}^5$ can be calculated by replacing $\B_1u=\frac{\partial u}{\partial \nv}- \ia \ka u=g_1$ by $\B_1u=\frac{\partial u}{\partial \nv}+\alpha u=g_1$
	in  \cite[Theorem 2.4]{FHM21Helmholtz} with $M=M_f=M_{g_1}=M_{g_3}=5$,  and replacing  $G^{V}_{M,m,n}$, $Q^{V}_{M,m,n}$, $G^{H}_{M,m,n}$ and $Q^{H}_{M,m,n}$ in \cite{FHM21Helmholtz} by \eqref{GVmn}, \eqref{QVmn}, \eqref{GHmn} and \eqref{QHmn}, respectively. 	
	Moreover, the maximum accuracy order of a 4-point finite difference scheme for $\B_1u=\frac{\partial u}{\partial \nv} +\alpha u=g_1$ and $\B_3u=\frac{\partial u}{\partial \nv}=g_3$  at the point $(x_0,y_0)$ with two smooth functions $\alpha(y)$ and $a(x,y)$   is six.
	\end{theorem}
	
 In particular, if   $a$ in \eqref{Qeques2} is a discontinuous piecewise constant coefficient, and $\B_1u=\frac{\partial u}{\partial \nv}+\alpha u=g_1$ with a constant $\alpha$, then the coefficients
 in \eqref{stencil:corner:1} are
\be \label{CR1:Corner:1}
		\begin{aligned}
			&C^{\mathcal{R}_1}_{0,0} =  \frac{4}{675}{\alpha}^5h^5-\frac{8}{675}{\alpha}^4h^4+\frac{8}{225}{\alpha}^3h^3-\frac{4}{25}{\alpha}^2h^2-\frac{17}{5}{\alpha}h-5,\ \ C^{\mathcal{R}_1}_{0,1} = \frac{1}{75}{\alpha}^2h^2+\frac{1}{5}{\alpha}h+2,\\	 
			&C^{\mathcal{R}_1}_{1,0} =  -\frac{4}{675}{\alpha}^4h^4+\frac{4}{225}{\alpha}^3h^3-\frac{4}{75}{\alpha}^2h^2+\frac{1}{5}{\alpha}h+2,\ \ C^{\mathcal{R}_1}_{1,1} =  1.\\			
	\end{aligned}
\ee		


\begin{theorem}  \label{thm:corner:2}
	Let  $(u_{h})_{i,j}$ be the numerical approximation of the exact solution $u$ of the elliptic interface problem \eqref{Qeques2} at the point $(x_i, y_j)$. Then the following discretization on a stencil centered at the corner point $(x_0,y_{N_2})$:
	\be \label{stencil:corner:2}
	{\footnotesize{
			\begin{aligned}
				\mathcal{L}^{\rr_2}_h u_h :=
				\begin{aligned}	
					 &C^{\mathcal{R}_2}_{0,-1}(u_{h})_{0,N_2-1}&             +&C^{\mathcal{R}_2}_{1,-1}(u_{h})_{1,N_2-1}\\
					 +&C^{\mathcal{R}_2}_{0,0}(u_{h})_{0,N_2}&           +&C^{\mathcal{R}_2}_{1,0}(u_{h})_{1,N_2}
				\end{aligned}
				= \sum_{(m,n)\in \ind_{4}} f^{(m,n)}C^{\rr_2}_{f,m,n} + \sum_{n=0}^{5}g_{1}^{(n)}C_{g_{1},n}^{\rr_2} + \sum_{n=0}^{5}g_{4}^{(n)}C_{g_{4},n}^{\rr_2},
			\end{aligned}
		}
	}
	\ee
	achieves sixth order of accuracy for $\B_1u=\frac{\partial u}{\partial \nv}+\alpha u=g_1$ and $\B_4u=\frac{\partial u}{\partial \nv}+\beta u=g_4$ at the point $(x_0,y_{N_2})$, where
 $\{C^{\rr_2}_{k,\ell}\}_{k\in \{0,1\},\ell \in \{-1,0\}}$, $\{C^{\rr_2}_{f,m,n}\}_{(m,n) \in \ind_4}$, $\{C^{\rr_2}_{g_1,n}\}_{n=0}^5$ and  $\{C^{\rr_2}_{g_4,n}\}_{n=0}^5$ can be  calculated by replacing $\B_1u=\frac{\partial u}{\partial \nv}- \ia \ka u=g_1$ and $\B_4u=\frac{\partial u}{\partial \nv}-\ia \ka u=g_4$ by $\B_1u=\frac{\partial u}{\partial \nv}+\alpha u=g_1$ and $\B_4u=\frac{\partial u}{\partial \nv}+\beta u=g_4$ respectively
	in \cite[Theorem 2.5]{FHM21Helmholtz} with $M=M_f=M_{g_1}=M_{g_4}=5$  and replacing  $G^{V}_{M,m,n}$, $Q^{V}_{M,m,n}$, $G^{H}_{M,m,n}$ and $Q^{H}_{M,m,n}$ in \cite{FHM21Helmholtz} by \eqref{GVmn}, \eqref{QVmn}, \eqref{GHmn} and \eqref{QHmn}, respectively.
	Moreover, the maximum accuracy order of a 4-point finite difference scheme for $\B_1u=\frac{\partial u}{\partial \nv} +\alpha u=g_1$ and $\B_4u=\frac{\partial u}{\partial \nv} +\beta u=g_4$  at the point $(x_0,y_{N_2})$ with three smooth functions $\alpha(y)$, $\beta(x)$ and $a(x,y)$   is six, where $\alpha(y_{N_2})\ne \beta(x_0)$.
	
	\end{theorem}
	Again, if  $a$ in \eqref{Qeques2} is a discontinuous constant coefficient, $\B_1u=\frac{\partial u}{\partial \nv}+\alpha u=g_1$ and $\B_4u=\frac{\partial u}{\partial \nv}+\beta u=g_4$ with  $\alpha$ and $\beta$ being constant, then the coefficients on the left hand side
	 in \eqref{stencil:corner:2} are
	\be \label{CR1:Corner:2}
	\begin{aligned}
		C^{\mathcal{R}_2}_{0,-1} &=\frac{1}{675}(4{\alpha}^5-6{\alpha}^4{\beta}+6{\alpha}^3{\beta}^2-4{\alpha}^2{\beta}^3)h^5+\frac{1}{675}(4{\alpha}^4-6{\alpha}^3{\beta}+6{\alpha}^2{\beta}^2-4{\alpha}{\beta}^3)h^4\\
		 &+\frac{1}{675}(9{\alpha}^2+63{\alpha}{\beta}-36{\beta}^2)h^2+\frac{1}{675}(135{\beta}+135{\alpha})h+2,\\
		C^{\mathcal{R}_2}_{0,0} &= \frac{1}{225}(-4{\alpha}^4+6{\alpha}^3{\beta}-6{\alpha}^2{\beta}^2+4{\alpha}{\beta}^3)h^4+\frac{1}{225}(8{\alpha}^3-18{\alpha}^2{\beta}-30{\alpha}{\beta}^2+16{\beta}^3)h^3\\
		 &+\frac{1}{225}(-36{\alpha}^2-357{\alpha}{\beta}-36{\beta}^2)h^2+\frac{1}{225}(-765{\alpha}-765{\beta})h-5,\\
		C^{\mathcal{R}_2}_{1,-1} &=\frac{1}{675}(-4{\alpha}^4+6{\alpha}^3{\beta}-6{\alpha}^2{\beta}^2+4{\alpha}{\beta}^3)h^4+1, \\
		C^{\mathcal{R}_2}_{1,0} &= \frac{1}{225}(4{\alpha}^3-6{\alpha}^2{\beta}+6{\alpha}{\beta}^2-4{\beta}^3)h^3+\frac{1}{225}(-12{\alpha}^2+21{\alpha}{\beta}+3{\beta}^2)h^2\\
		&+\frac{1}{225}(45{\beta}+45{\alpha})h+2.					 
	\end{aligned}
	\ee	
	When $\alpha=\beta$, we further have
	 $C^{\mathcal{R}_2}_{0,-1}=C^{\mathcal{R}_2}_{1,0}=\frac{4}{75}{\beta}^2h^2+\frac{2}{5}{\beta}h+2$ and $C^{\mathcal{R}_2}_{1,-1}=1$ in \eqref{CR1:Corner:2}.

%
	
	\subsection{Stencils for irregular points}\label{hybrid:Irregular:points}
	Let $(x_i,y_j)$ be an irregular point (i.e., both $d_{i,j}^+$ and $d_{i,j}^-$ are nonempty, see \cref{Extend the compact} for an example) and choose the base point $(x^*_i,y^*_j)\in \Gamma_I \cap (x_i-h,x_i+h)\times (y_j-h,y_j+h)$.
	By \eqref{base:pt}, we have
	\begin{equation} \label{base:pt:gamma}
		x_i^*=x_i-v_0h  \quad \mbox{and}\quad y_j^*=y_j-w_0h  \quad \mbox{with}\quad
		-1<v_0, w_0<1 \quad \mbox{and}\quad (x_i^*,y_j^*)\in \Gamma_I.
	\end{equation}
	Let $a_{\pm}$, $u_{\pm}$ and $f_{\pm}$ represent the coefficient function $a$, the solution $u$ and source term $f$ in $\Omega^{\pm}$.
	Similar to \eqref{ufmn}, we define that
	\begin{align*}
		& a_{\pm}^{(m,n)}:=\frac{\partial^{m+n} a_{\pm}}{ \partial^m x \partial^n y}(x^*_i,y^*_j),\qquad u_{\pm}^{(m,n)}:=\frac{\partial^{m+n} u_{\pm}}{ \partial^m x \partial^n y}(x^*_i,y^*_j),\qquad f_{\pm}^{(m,n)}:=\frac{\partial^{m+n} f_{\pm}}{ \partial^m x \partial^n y}(x^*_i,y^*_j), \\
		& \gd^{(m,n)}:=\frac{\partial^{m+n} \gd}{ \partial^m x \partial^n y}(x^*_i,y^*_j),\qquad
		\gn^{(m,n)}:=\frac{\partial^{m+n} \gn}{ \partial^m x \partial^n y}(x^*_i,y^*_j).
	\end{align*}
	Similar to \eqref{u:approx:key:V}, we have
	\begin{align*}
		\label{u:approx:ir:key:2}
		u_\pm (x+x_i^*,y+y_j^*)
		& =\sum_{(m,n)\in \ind_{M+1}^{V,1}}
		u_\pm^{(m,n)} G^{\pm,V}_{M,m,n}(x,y) +\sum_{(m,n)\in \ind_{M-1}}
		f_\pm ^{(m,n)} Q^{\pm,V}_{M,m,n}(x,y)+\bo(h^{M+2}),
	\end{align*}
	for $x,y\in (-2h,2h)$, where $\ind_{M-1}$ and $\ind_{M+1}^{V,1}$ are defined in \eqref{Sk} and  \eqref{indV12} respectively, $G^{\pm,V}_{M,m,n}(x,y)$ and $Q^{\pm,V}_{M,m,n}(x,y)$ are obtained by  replacing $\{a^{(m,n)}: (m,n) \in \ind_{M}\}$ by $\{a_{\pm}^{(m,n)}: (m,n) \in \ind_{M}\}$  in \eqref{GVmn} and \eqref{QVmn}.
	As in \cite{FHM21a,FHM21b,FHM21Helmholtz}, near the point $(x_i^*,y_j^*)$, the parametric equation of $\Gamma_I$ can be written as:
	\be \label{parametric}
	x=r(t)+x_i^*,\quad y=s(t)+y_j^*,\quad
	(r'(t))^2+(s'(t))^2>0
	\quad \mbox{for}\;\; t\in (-\epsilon,\epsilon) \quad \mbox{with}\quad \epsilon>0,
	\ee
	where $r$ and $s$ are smooth functions. Similarly to the definition of the 9-point compact stencil in \eqref{Compact:Set}, we define the following 4-point set for the 13-point scheme:
	 \be\label{13point:Set}
	\begin{split}
		& e_{i,j}^+:=\{(k,\ell) \; : \;
		(k,\ell)\in \{(-2,0),(0,-2),(0,2),(2,0)\}, \psi(x_i+kh, y_j+\ell h)\ge 0\}, \quad \mbox{and}\\
		& e_{i,j}^-:=\{(k,\ell) \; : \;
		(k,\ell)\in \{(-2,0),(0,-2),(0,2),(2,0)\}, \psi(x_i+kh, y_j+\ell h)<0\}.
	\end{split}
	\ee
	
	In the next theorem we present a simplified version of  \cite[Theorem 3.2]{FHM21b}, adapted to the aim of developing of a fifth order hybrid 13-point scheme for irregular points.
	\begin{theorem}\label{thm:interface}
		Let $u$ be the solution to the elliptic interface problem in \eqref{Qeques2} and let  $\Gamma_I$ be parameterized near $(x_i^*,y_j^*)$ by \eqref{parametric}.
		Then
\be \label{tranmiss:cond}
\begin{split}
	u_-^{(m',n')}&=\sum_{  \substack{ (m,n)\in \ind_{M+1}^{V,1} \\  m+n \le m'+n'} }T^{u_{+}}_{m',n',m,n}u_+^{(m,n)}+\sum_{(m,n)\in \ind_{M-1}} \left(T^+_{m',n',m,n} f_+^{(m,n)}
	+ T^-_{m',n',m,n} f_{-}^{(m,n)}\right)\\
	&+\sum_{(m,n)\in \ind_{M+1}} T^{g_D}_{m',n',m,n} g_{D}^{(m,n)}
	+\sum_{(m,n)\in \ind_{M}} T^{g_N}_{m',n',m,n} g_N^{(m,n)},\qquad \forall\; (m',n')\in \ind_{M+1}^{V,1},
\end{split}
\ee		
		where all the transmission coefficients $T^{u_+}, T^{\pm}, T^{\gd}, T^{\gn}$ are uniquely determined by
		$r^{(k)}(0)$, $s^{(k)}(0)$  for $k=0,\ldots,M+1$ and $\{a_{\pm}^{(m,n)}: (m,n)\in \ind_{M}\}$.	Moreover, let  $T^{u_+}_{m',n',m,n}$ be the transmission coefficient of $u_+^{(m,n)}$ in \eqref{tranmiss:cond} with $(m,n)\in \ind_{M+1}^{V,1}$, $m+n=m'+n'$ and $(m',n')\in \ind_{M+1}^{V,1}$. Then $T^{u_+}_{m',n',m,n}$ only depends on $r^{(k)}(0)$, $s^{(k)}(0)$  for $k=0,\ldots,M+1$ of \eqref{parametric} and $a_{\pm}^{(0,0)}$. Particularly, 	
		\be\label{T0000}
T^{u_+}_{0,0,0,0}=1
\quad \mbox{and}\quad
T^{u_+}_{m',n',0,0}=0 \quad \mbox{if }  (m',n')\ne (0,0).
		\ee
		
	\end{theorem}


	Next, we provide the 13-point finite difference scheme for interior irregular points.
	
	\begin{theorem}\label{13point:Inter}
		Let $(u_{h})_{i,j}$ be the numerical approximation to the solution of \eqref{Qeques2} at an interior irregular point $(x_i, y_j)$. Pick a base point
		$(x_i^*,y_j^*)$ as in \eqref{base:pt:gamma}.
		Then
		the following 13-point scheme centered at the interior irregular point $(x_i,y_j)$:
		\begin{align}\label{13point:interface}
			\mathcal{L}^{\Gamma_I}_h & :=
			\begin{aligned}
				 & & & &	 &C_{0,-2}(u_{h})_{i,j-2}&\\
			& &+&C_{-1,-1}(u_{h})_{i-1,j-1}&
			+&C_{0,-1}(u_{h})_{i,j-1}&
			+&C_{1,-1}(u_{h})_{i+1,j-1}& \\
			 +&C_{-2,0}(u_{h})_{i-2,j}&+&C_{-1,0}(u_{h})_{i-1,j}&
			+&C_{0,0}(u_{h})_{i,j}&
			+&C_{1,0}(u_{h})_{i+1,j} &+&C_{2,0}(u_{h})_{i+2,j}&\\
				& &+&C_{-1,1}(u_{h})_{i-1,j+1}&
			+&C_{0,1}(u_{h})_{i,j+1}&
			+&C_{1,1}(u_{h})_{i+1,j+1}& \\
			& & & &	+&C_{0,2}(u_{h})_{i,j+2}&
			\end{aligned}\\
			\nonumber
			&=\sum_{(m,n)\in \ind_{3} } f_+^{(m,n)}J^{+}_{m,n} + \sum_{(m,n)\in \ind_{3}} f_-^{(m,n)}J^{-}_{m,n} +\sum_{(m,n)\in \ind_{5}} \gd^{(m,n)}J^{\gd}_{m,n} + \sum_{(m,n)\in \ind_{4}} \gn^{(m,n)}J^{\gn}_{m,n},
		\end{align}
		achieves fifth order accuracy, where all $\{C_{k,\ell}\}$ in \eqref{13point:interface} are calculated by \eqref{AX0}, $J^{\pm}_{m,n}:=
		J_{m,n}^{\pm,0}+J^{\pm,T}_{m,n}$ for all $(m,n)\in \ind_{3}$,
		{\footnotesize
			\begin{align*}
				& J^{\pm,0}_{m,n}:=\sum_{(k,\ell)\in d_{i,j}^\pm \cup e_{i,j}^\pm} C_{k,\ell} Q^{\pm,V}_{4,m,n}((v_0 +k)h,(w_0+\ell) h), \quad J^{\pm,T}_{m,n}:=
				\sum_{  \substack{ (m',n')\in \ind_{5}^{V,1} }} I^{-}_{m',n'} T^{\pm}_{m',n',m,n}, \quad \forall (m,n) \in \ind_{3},\\
				& J^{\gd}_{m,n}:=
				\sum_{ \substack{ (m',n')\in \ind_{5}^{V,1} }} I^{-}_{m',n'} T^{\gd}_{m',n',m,n},
				\quad \forall (m,n) \in \ind_{5}, \quad
				J^{\gn}_{m,n}:=
				\sum_{ \substack{ (m',n')\in \ind_{5}^{V,1} }} I^{-}_{m',n'} T^{\gn}_{m',n',m,n}, \quad \forall (m,n) \in \ind_{4},\\
				& I^{-}_{m,n}:=\sum_{(k,\ell)\in d_{i,j}^- \cup e_{i,j}^-}
				C_{k,\ell} G^{-,V}_{4,m,n}((v_0+k)h,(w_0+\ell) h), \quad \forall (m,n) \in \ind_{5}^{V,1}.
			\end{align*}
		}
		Moreover, the maximum accuracy order of a 13-point finite difference stencil for \eqref{Qeques2} at an interior irregular point $(x_i, y_j)$ is five.
	\end{theorem}

For the $13$-point scheme in \cref{13point:Inter}, if only one point in the set $\{(x_i-h,y_j-h),(x_i-h,y_j+h),(x_i+h,y_j-h),(x_i+h,y_j+h)\}$ belongs to $\Omega^{-}$ and the other 12 points all belong to $\Omega^{+}$, we can set $C_{k,\ell}=0$ for $(x_i+k h, y_j+\ell h)\in \Omega^{-}$, $x^*_i=x_i$, $y^*_i=y_i$ to achieve sixth order accuracy in $(x_i,y_j)$.

Finally, we provide a way of achieving an efficient implementation for the $13$-point scheme in irregular points in \cref{13point:Inter}.
	
\textbf{Efficient implementation details:}\\

By \cref{thm:interface},
a simpler $J^{u_+,T}_{m,n}(h)$ in \cite[(3.26)]{FHM21b} can be written as:
\be\label{newJuT}
J^{u_+,T}_{m,n}(h):=
\sum_{  \substack{ (m',n')\in \ind_{M+1}^{V,1} \\  m'+n' \ge m+n}} I^-_{m',n'}(h) T^{u_+}_{m',n',m,n}.
\ee
Replacing  $\ind_{M+1}^{1}$ by $\ind_{M+1}^{V,1}$ for \cite[(3.28) and (3.29)]{FHM21b}, we have
\be\label{IJT}
I^+_{m,n}(h)+J^{u_+,T}_{m,n}(h)=\bo(h^{M+2}),  \  h\to 0, \; \mbox{ for all }\; (m,n)\in \ind_{M+1}^{V,1}.
\ee
Replacing  $G^{\pm}_{m,n}$, $H^{\pm}_{m,n}$  and $d_{i,j}^\pm$  by $G^{\pm,V}_{M,m,n}$, $Q^{\pm,V}_{M,m,n}$ and $d_{i,j}^\pm \cup e_{i,j}^\pm$ for  \cite[(3.25) and (3.26)]{FHM21b}, we obtain
\[
\sum_{(k,\ell)\in d_{i,j}^+\cup e_{i,j}^+}
C_{k,\ell}(h) G^{+,V}_{M,m,n}(v_0h+kh,w_0h+\ell h)+\sum_{  \substack{ (m',n')\in \ind_{M+1}^{V,1} \\  m'+n' \ge m+n}} I^-_{m',n'}(h) T^{u_+}_{m',n',m,n}=\bo(h^{M+2}),
\]	
and
\[
\begin{split}
	&\sum_{(k,\ell)\in d_{i,j}^+\cup e_{i,j}^+}
	C_{k,\ell}(h) G^{+,V}_{M,m,n}(v_0h+kh,w_0h+\ell h)\\
	&+\sum_{  \substack{ (m',n')\in \ind_{M+1}^{V,1} \\  m'+n' \ge m+n}} \sum_{(k,\ell)\in d_{i,j}^-\cup e_{i,j}^-}
	C_{k,\ell}(h) G^{-,V}_{M,m',n'}(v_0h+kh,w_0h+\ell h) T^{u_+}_{m',n',m,n}=\bo(h^{M+2}).
\end{split}
\]
So, \eqref{IJT} is equivalent to
\be\label{eq:interface:sum}
\begin{split}
	&\sum_{(k,\ell)\in d_{i,j}^-\cup e_{i,j}^-}
	C_{k,\ell}(h) \sum_{  \substack{ (m',n')\in \ind_{M+1}^{V,1} \\  m'+n' \ge m+n}}  G^{-,V}_{M,m',n'}(v_0h+kh,w_0h+\ell h) T^{u_+}_{m',n',m,n}\\
	&+\sum_{(k,\ell)\in d_{i,j}^+\cup e_{i,j}^+}
	C_{k,\ell}(h) G^{+,V}_{M,m,n}(v_0h+kh,w_0h+\ell h)=\bo(h^{M+2}), \qquad \mbox{for all } \ (m,n)\in \ind_{M+1}^{V,1}.
\end{split}
\ee		
Let
\be\label{ChandX}
C_{k,\ell}(h):=\sum_{i=0}^{M+1} c_{k,\ell,i}h^i,  \qquad
X_{k,\ell}:=( c_{k,\ell,0}, c_{k,\ell,1},\dots,c_{k,\ell,M+1})^T.
\ee
Since $G^{\pm,V}_{M,m,n}((k+v_0)h,(\ell+w_0)h)$ is the polynomial of $h$ and the degree of $h$ of every term in $G^{\pm,V}_{M,m,n}((k+v_0)h,(\ell+w_0)h)$ is non-negative, we deduce that
\be\label{eq:plus}
C_{k,\ell}(h) G^{+,V}_{M,m,n}((k+v_0)h,(\ell+w_0)h)=DA^{+,m,n}_{k,\ell}X_{k,\ell}+\bo(h^{M+2}),
\ee
\be\label{eq:minus}
C_{k,\ell}(h) \sum_{  \substack{ (m',n')\in \ind_{M+1}^{V,1} \\  m'+n' \ge m+n}}  G^{-,V}_{M,m',n'}((k+v_0)h,(\ell+w_0)h) T^{u_+}_{m',n',m,n}=DA^{-,m,n}_{k,\ell}X_{k,\ell}+\bo(h^{M+2}),
\ee		
where
\[
D=(h^0,h^1,\dots,h^{M+1}),
\]
and $A^{\pm,m,n}_{k,\ell}$ is independent for $h$ for all $(m,n)\in \ind_{M+1}^{V,1}$.
So \eqref{eq:interface:sum} is equivalent to

\be\label{eq:interface:sum:matrix}
\begin{split}
	\sum_{(k,\ell)\in d_{i,j}^+\cup e_{i,j}^+}
	 DA^{+,m,n}_{k,\ell}X_{k,\ell}+\sum_{(k,\ell)\in d_{i,j}^-\cup e_{i,j}^-}
DA^{-,m,n}_{k,\ell}X_{k,\ell}=\bo(h^{M+2}), \qquad \mbox{for all } \ (m,n)\in \ind_{M+1}^{V,1}.
\end{split}
\ee	

Define
\be\label{Amnkl}
A^{m,n}_{k,\ell}:=
\begin{cases}
	A^{+,m,n}_{k,\ell}, &\text{if } (k,\ell)\in d_{i,j}^+\cup e_{i,j}^+,\\
	A^{-,m,n}_{k,\ell}, &\text{if } (k,\ell)\in d_{i,j}^-\cup e_{i,j}^-.
\end{cases}
\ee
Then \eqref{eq:interface:sum:matrix} is equivalent to
\[
A^{m,n}X=0, \qquad \mbox{for all } \ (m,n)\in \ind_{M+1}^{V,1},
\]
where
\be\label{Amn:Totla}
A^{m,n}=(A^{m,n}_{-1,-1},A^{m,n}_{-1,0},A^{m,n}_{-1,1},A^{m,n}_{0,-1},A^{m,n}_{0,0},A^{m,n}_{0,1},A^{m,n}_{1,-1},A^{m,n}_{1,0},A^{m,n}_{1,1},A^{m,n}_{-2,0},A^{m,n}_{2,0},A^{m,n}_{0,-2},A^{m,n}_{0,2}),
\ee
and
\be\label{X:Totla}
X=(X_{-1,-1}^T,X_{-1,0}^T,X_{-1,1}^T,X_{0,-1}^T,X_{0,0}^T,X_{0,1}^T,X_{1,-1}^T,X_{1,0}^T,X_{1,1}^T,X_{-2,0}^T,X_{2,0}^T,X_{0,-2}^T,X_{0,2}^T)^T.
\ee
Let
\be\label{TotalA}
A=\left(
(A^{0,0})^T,
(A^{0,1})^T,
\dots,
(A^{0,M+1})^T,
(A^{1,0})^T,
(A^{1,1})^T,
\dots,
(A^{1,M})^T\right)^T.
\ee
Finally, \eqref{eq:interface:sum} is equivalent to
\be\label{AX0}
AX=0.
\ee
Since we use 13-point scheme for the irregular points, we have 13 components in  \eqref{Amn:Totla} and  \eqref{X:Totla}. If we use 9-point compact scheme for the irregular points, we only need to delete the last four components in \eqref{Amn:Totla} and  \eqref{X:Totla}. For the 25-point or 36-point schemes for the irregular points, the only change is to add more $A^{m,n}_{k,\ell}$ and $X_{k,\ell}$ in \eqref{Amn:Totla} and  \eqref{X:Totla}. Even
there are many different cases for the $13$-point schemes for the irregular points depending on how the interface curve $\Gamma_I$ partitions the $13$ points in it, we can repeatedly use $A^{\pm,m,n}_{k,\ell}$ in \eqref{eq:plus}, \eqref{eq:minus} and \eqref{Amnkl} to cover all the cases which significantly reduce the computation cost and make the implementation very effective and flexible. Furthermore, if we want to obtain the lower or higher finite schemes for irregular points, we only need to delete or add some $A^{0,n+1}$ and $A^{1,n}$ in \eqref{TotalA}.

After the above simplification, we find that the $A$ in \eqref{AX0} is a 36 by 78 matrix for the 13-point scheme with fifth order accuracy while $A$ is a 16 by 36 matrix  and the 9-point scheme with third order accuracy.
Observing the following identity (whose proof is given in \cref{hybrid:sec:proofs})
\be\label{ckln:irregular}		 c_{0,-2,i}+c_{-2,0,i}+c_{2,0,i}+c_{0,2,i}+\sum\limits_{k=-1}^{1} \sum\limits_{\ell=-1}^{1}
		c_{k,\ell,i}=0, \quad \mbox{for} \quad i=0,1,\dots, M+1,
		\ee	
we can further reduce the size of the matrix $A$ in \eqref{AX0} to $30$ by $72$ for the $13$-point scheme.

	\section{Numerical experiments}
	\label{sec:numerical}

Let $\Omega=(l_1,l_2)\times(l_3,l_4)$ with
$l_4-l_3=N_0(l_2-l_1)$ for some positive integer $N_0$. For a given $J\in \NN$, we define
$h:=(l_2-l_1)/N_1$ with $N_1:=2^J$ and let
$x_i=l_1+ih$ and
$y_j=l_3+jh$ for $i=0,1,\dots,N_1$ and $j=0,1,\dots,N_2$ with $N_2:=N_0N_1$.
Let
$u(x,y)$ be the exact solution of \eqref{Qeques2} and $(u_{h})_{i,j}$ be a numerical solution at $(x_i, y_j)$ using the mesh size $h$.
We  measure the consistency of the proposed scheme in the $l_2$ norm by the relative error
$\frac{\|u_{h}-u\|_{2}}{\|u\|_{2}}$, if the exact solution $u$ is available.
If it is not, then we quantify the consistency error
by  ${\|u_{h}-u_{h/2}\|_{2}}$, where
\begin{align*}
	\|u_{h}-u\|_{2}^2:= h^2&\sum_{i=0}^{N_1}\sum_{j=0}^{N_2} \left((u_h)_{i,j}-u(x_i,y_j)\right)^2, \ \ \|u\|_{2}^2:=h^2 \sum_{i=0}^{N_1}\sum_{j=0}^{N_2} \left(u(x_i,y_j)\right)^2,\\
	&\|u_{h}-u_{h/2}\|_{2}^2:= h^2\sum_{i=0}^{N_1}\sum_{j=0}^{N_2} \left((u_{h})_{i,j}-(u_{h/2})_{2i,2j}\right)^2.
\end{align*}
In addition we also provide results for the infinity norm of the errors given by:
\[
\|u_h-u\|_\infty
:=\max_{0\le i\le N_1, 0\le j\le N_2} \left|(u_h)_{i,j}-u(x_i,y_j)\right|,
\quad
\|u_{h}-u_{h/2}\|_\infty:=\max_{0\le i\le N_1,0\le j\le N_2} \left|(u_{h})_{i,j}-(u_{h/2})_{2i,2j}\right|.
\]

\subsection{Numerical examples with known $u$}

In this subsection, we provide five numerical examples with a known solution $u$ of \eqref{Qeques2}. Note that the maximum accuracy order for the compact 9-point finite difference scheme in irregular and regular points, for elliptic interface problems with discontinuous coefficients, is three  and six, respectively. So, in \cref{hybrid:ex3,hybrid:ex5} we compare the proposed hybrid scheme with the compact 9-point scheme of a  sixth order of accuracy at  regular points and third order of accuracy at irregular points. That is, both uses the same compact $9$-point stencils with accuracy order six at all regular points, and they only differ at irregular points such that the proposed hybrid scheme uses $13$-point stencils having fifth order accuracy, while the compact $9$-point scheme uses $9$-point stencils having third order accuracy. Their computational costs are comparable,
because the percentage of the number of irregular points over all the grid points 
decays exponentially to $0$ at the rate $\bo(2^{-J})$, e.g., this percentage is less than or around $1\%$ at the level $J=9$ for all our numerical examples.

The  five numerical examples can be characterized as follows:

\begin{itemize}	
	\item \cref{hybrid:ex3,hybrid:ex5} compare the proposed hybrid scheme and the $9$-point compact scheme.
\item In all examples, either $a_+/a_-$ or $a_-/a_+$ is very large on $\Gamma_I$ for high contrast coefficients $a$.
	\item 4-side Dirichlet boundary conditions are demonstrated in  \cref{hybrid:ex3,hybrid:ex4,hybrid:ex5}.
	\item 1-side Dirichlet, 1-side Neumann and 2-side Robin boundary conditions
are considered
in \cref{hybrid:ex1,hybrid:ex2}.
	\item Results for smooth interface curves $\Gamma_I$ are presented in \cref{hybrid:ex1,hybrid:ex2,hybrid:ex3,hybrid:ex4}.
	\item Results for a  sharp-edged interface curve $\Gamma_I$ are demonstrated in \cref{hybrid:ex5}.
	\item Results for two constant jump functions $g_D$ and $g_N$  are shown in \cref{hybrid:ex1,hybrid:ex2,hybrid:ex3,hybrid:ex4}.
	\item Results for two non-constant jump functions $g_D$ and $g_N$ are presented in \cref{hybrid:ex5}.
	
\end{itemize}

\begin{figure}[htbp]
	\centering
	\begin{subfigure}[b]{0.20\textwidth}
		 \includegraphics[width=3.0cm,height=3cm]{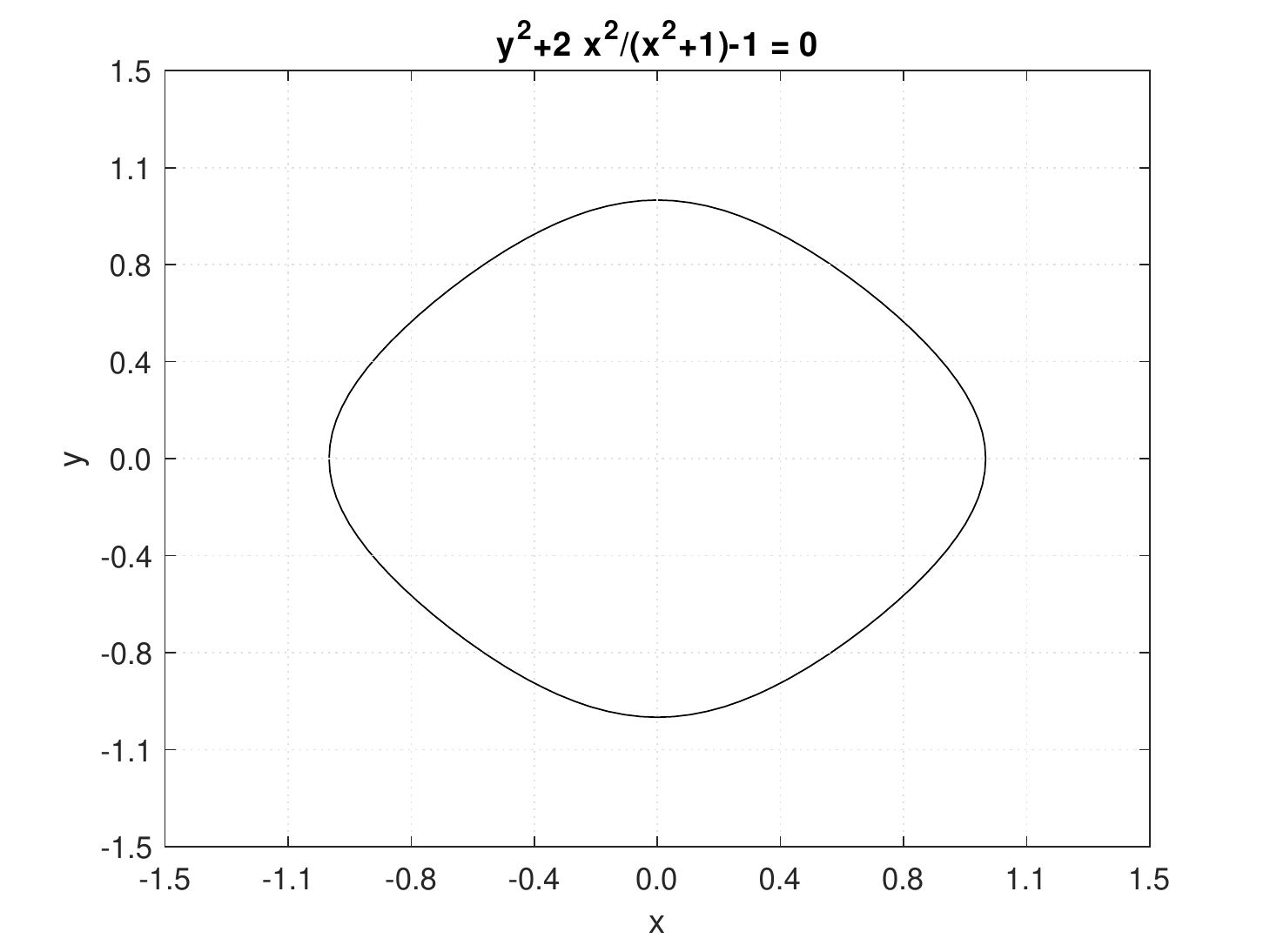}
	\end{subfigure}
	\begin{subfigure}[b]{0.25\textwidth}
		 \includegraphics[width=3.5cm,height=3.5cm]{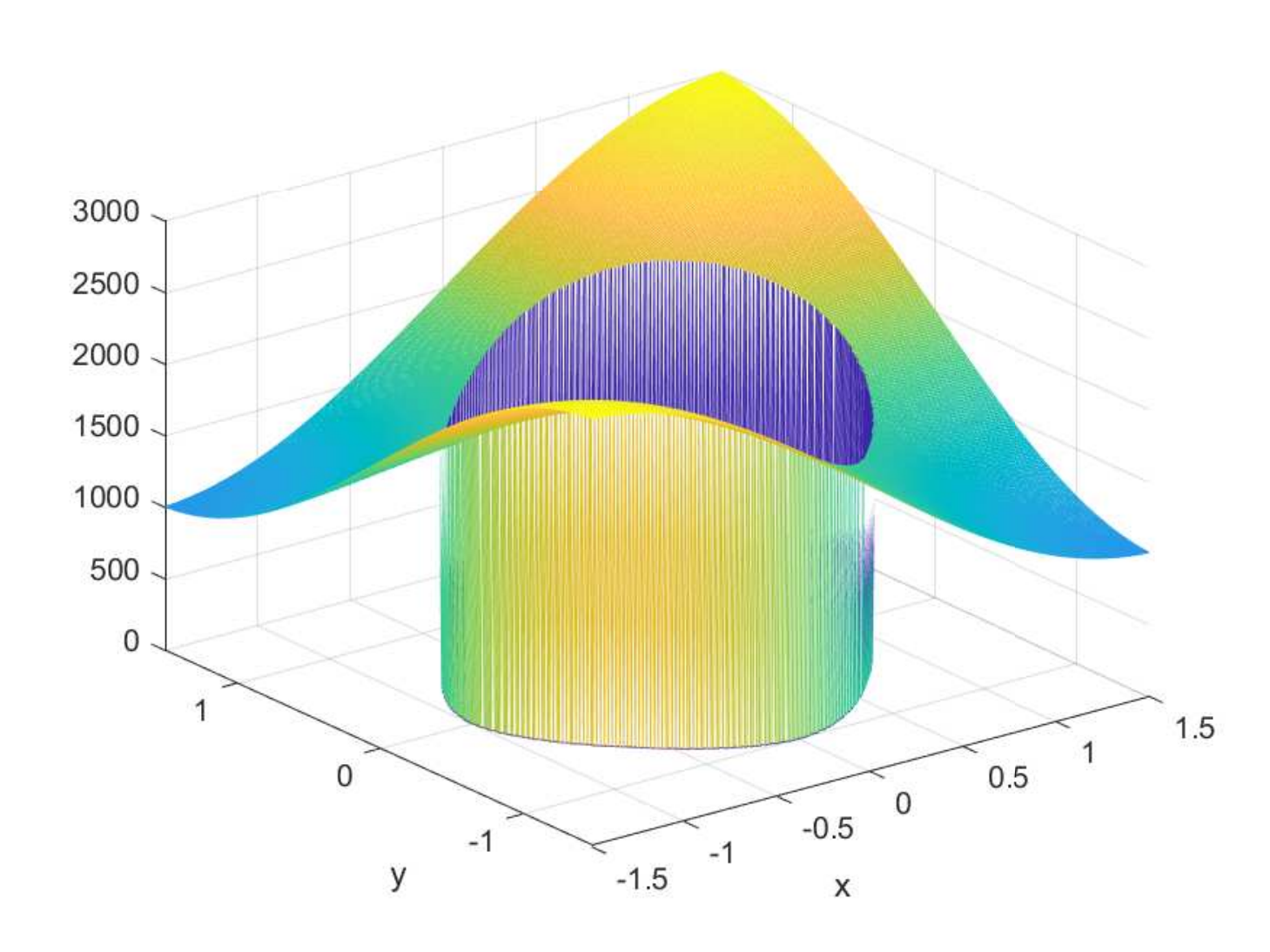}
	\end{subfigure}
	\begin{subfigure}[b]{0.25\textwidth}
		 \includegraphics[width=3.5cm,height=3.5cm]{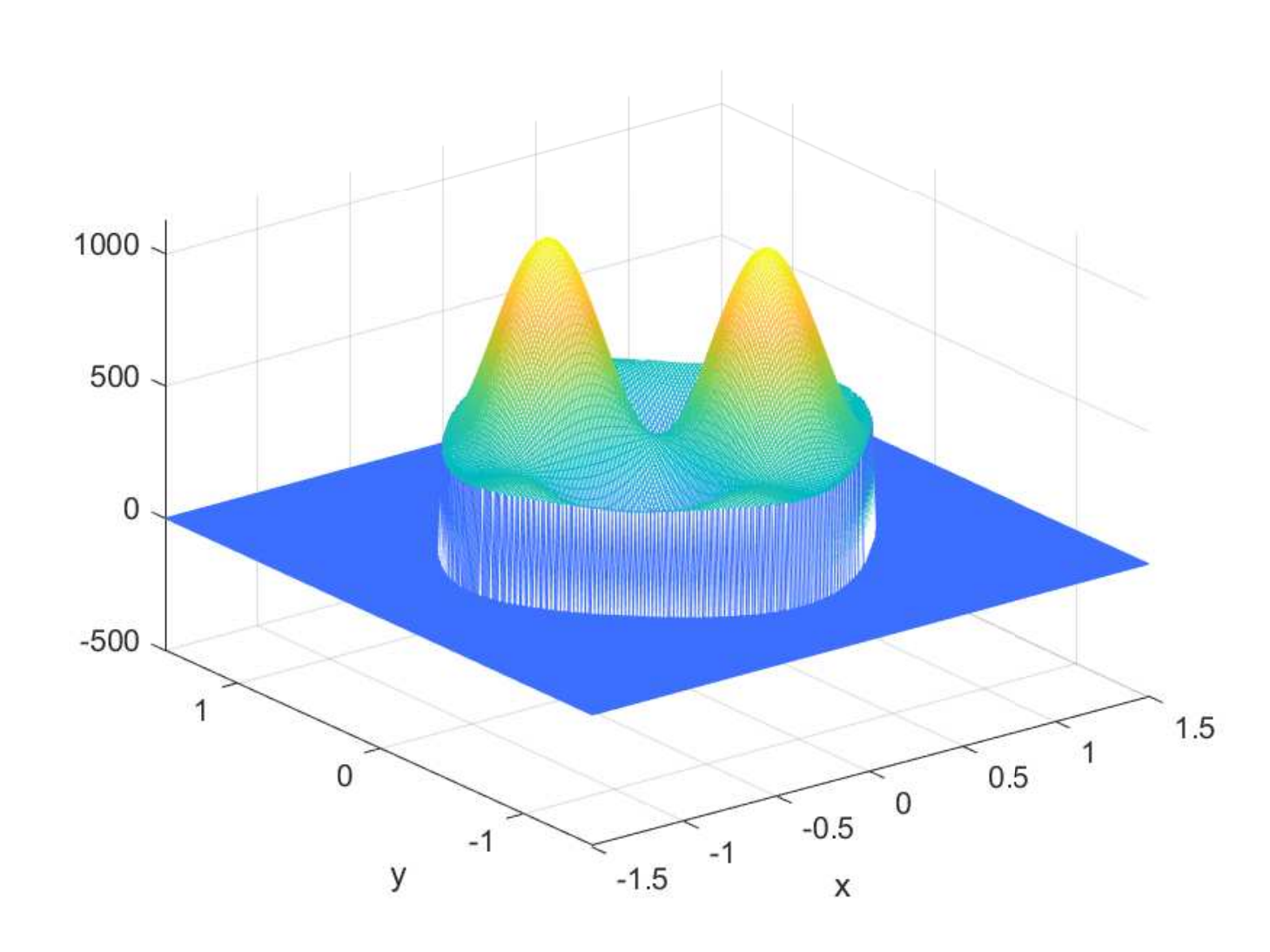}
	\end{subfigure}
	\begin{subfigure}[b]{0.25\textwidth}
		 \includegraphics[width=3.5cm,height=3.5cm]{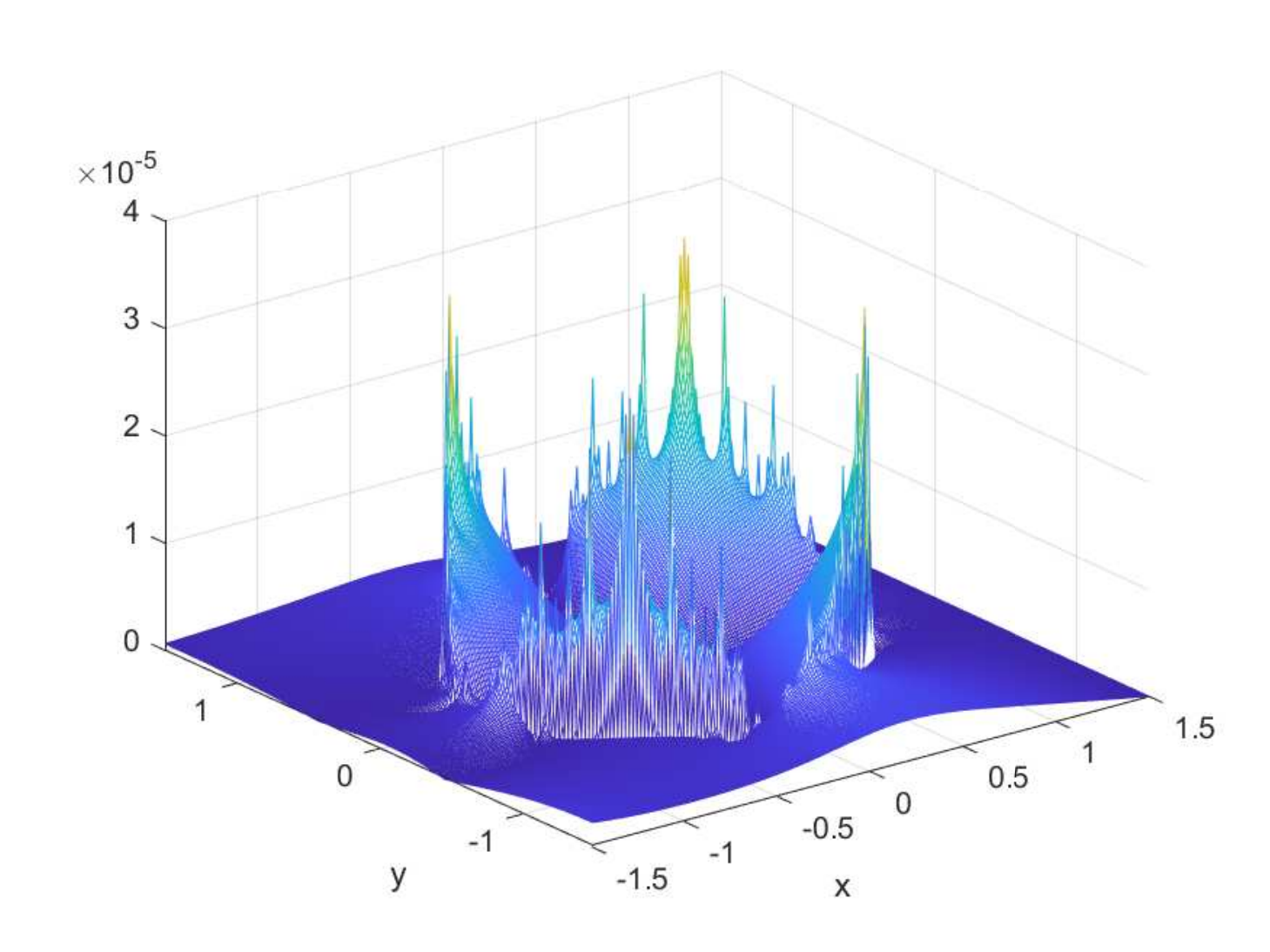}
	\end{subfigure}
	\caption
	{\cref{hybrid:ex3}: the interface curve $\Gamma_I$ (first panel), the coefficient $a(x,y)$ (second panel),  the numerical solution $u_h$ (third panel), and the error $|u_h-u|$ (fourth panel) with $h=2^{-8}\times 3$,  where $u_h$ is computed by our proposed hybrid finite difference scheme.}
	\label{hybrid:fig:QSp3}
\end{figure}

\begin{example}\label{hybrid:ex3}
	\normalfont
	Let $\Omega=(-1.5,1.5)^2$ and
	the interface curve be given by
	$\Gamma_I:=\{(x,y)\in \Omega:\; \psi(x,y)=0\}$ with
	$\psi (x,y)=y^2+\frac{2x^2}{x^2+1}-1$.
	The functions  in \eqref{Qeques2} are given by
	\begin{align*}
		&a_{+}=10^3(2+\sin(x)\sin(y)),
		\qquad a_{-}=10^{-3}(2+\sin(x)\sin(y)), \qquad g_D=-200, \qquad g_N=0,\\
		&u_{+}=10^{-3}\sin(4 x)\sin(4 y)(y^2(x^2+1)+x^2-1),\\
		&u_{-}=10^{3}\sin(4 x)\sin(4 y)(y^2(x^2+1)+x^2-1)+200,\\
		&  u(-1.5,y)=g_1,
		\qquad u(1.5,y)=g_2, \qquad \mbox{for} \qquad y\in(-1.5,1.5),\\
		& u(x,-1.5)=g_3,
		\qquad  u(x,1.5)=g_4,  \qquad \mbox{for} \qquad x\in(-1.5,1.5),
	\end{align*}
	the other functions $f^{\pm}$, $g_1, \ldots,g_4$ in \eqref{Qeques2} can be obtained by plugging the above functions into \eqref{Qeques2}.
Note the high contrast $a_+/a_-=10^6$ on $\Gamma_I$.
The numerical results are presented in \cref{hybrid:table:QSp3} and \cref{hybrid:fig:QSp3}.	
\end{example}

\begin{table}[htbp]
	\caption{Performance in \cref{hybrid:ex3}  of our proposed hybrid finite difference scheme and compact 9-point scheme on uniform Cartesian meshes with $h=2^{-J}\times 3$. $\kappa$ is the condition number of the coefficient matrix.}
	\centering
	\setlength{\tabcolsep}{0.5mm}{
		\begin{tabular}{c|c|c|c|c|c|c|c|c|c|c}
			\hline
\multicolumn{1}{c|}{} &
\multicolumn{5}{c|}{Our proposed hybrid scheme} &
\multicolumn{5}{c}{Compact 9-point scheme} \\
\cline{1-11}			
			$J$
			& $\frac{\|u_{h}-u\|_{2}}{\|u\|_{2}}$
			
			&order & $\|u_{h}-u\|_{\infty}$
			
			&order &  $\kappa$ & $\frac{\|u_{h}-u\|_{2}}{\|u\|_{2}}$
			
			&order & $\|u_{h}-u\|_{\infty}$
			
			&order &  $\kappa$  \\
			\hline
4    &1.493E-01    &0    &1.362E+02    &0    &2.136E+02    &5.465E-01    &0    &4.515E+02    &0    &8.685E+01\\
5    &3.124E-03    &5.6    &3.872E+00    &5.1    &4.262E+02    &4.751E-02    &3.5    &4.453E+01    &3.3    &4.896E+02\\
6    &6.081E-05    &5.7    &7.168E-02    &5.8    &6.261E+03    &2.464E-03    &4.3    &2.890E+00    &3.9    &2.069E+03\\
7    &1.238E-06    &5.6    &1.490E-03    &5.6    &1.701E+04    &2.745E-04    &3.2    &3.318E-01    &3.1    &9.171E+03\\
8    &1.803E-08    &6.1    &3.305E-05    &5.5    &1.169E+05    &1.557E-05    &4.1    &1.894E-02    &4.1    &4.054E+04\\
9    &    &    &    &    &    &9.053E-07    &4.1    &1.185E-03    &4.0    &1.648E+05\\			
			\hline
	\end{tabular}}
	\label{hybrid:table:QSp3}
\end{table}

\begin{figure}[htbp]
	\centering
	\begin{subfigure}[b]{0.2\textwidth}
		 \includegraphics[width=3.0cm,height=3.0cm]{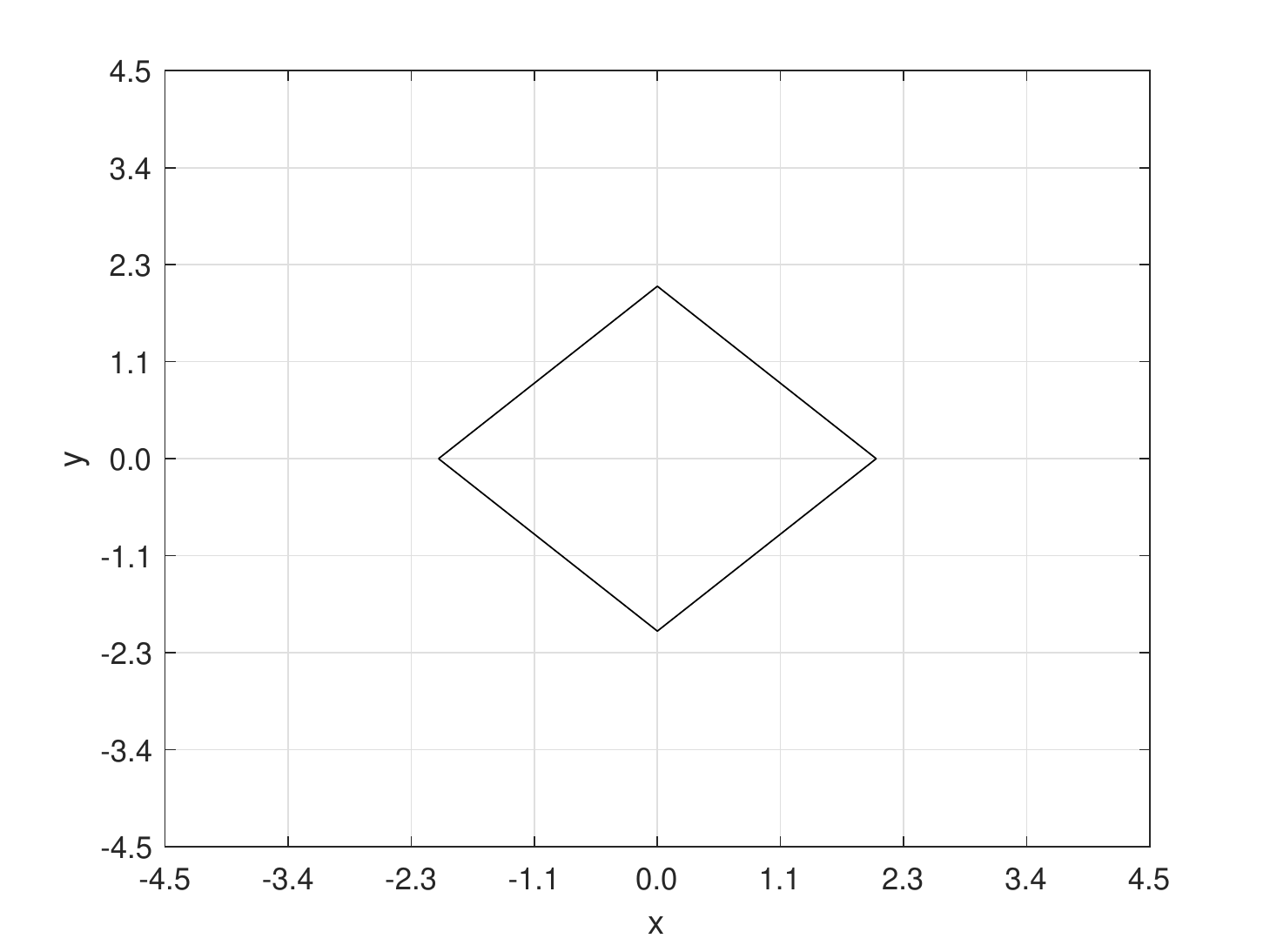}
	\end{subfigure}
	\begin{subfigure}[b]{0.25\textwidth}
		 \includegraphics[width=3.5cm,height=3.5cm]{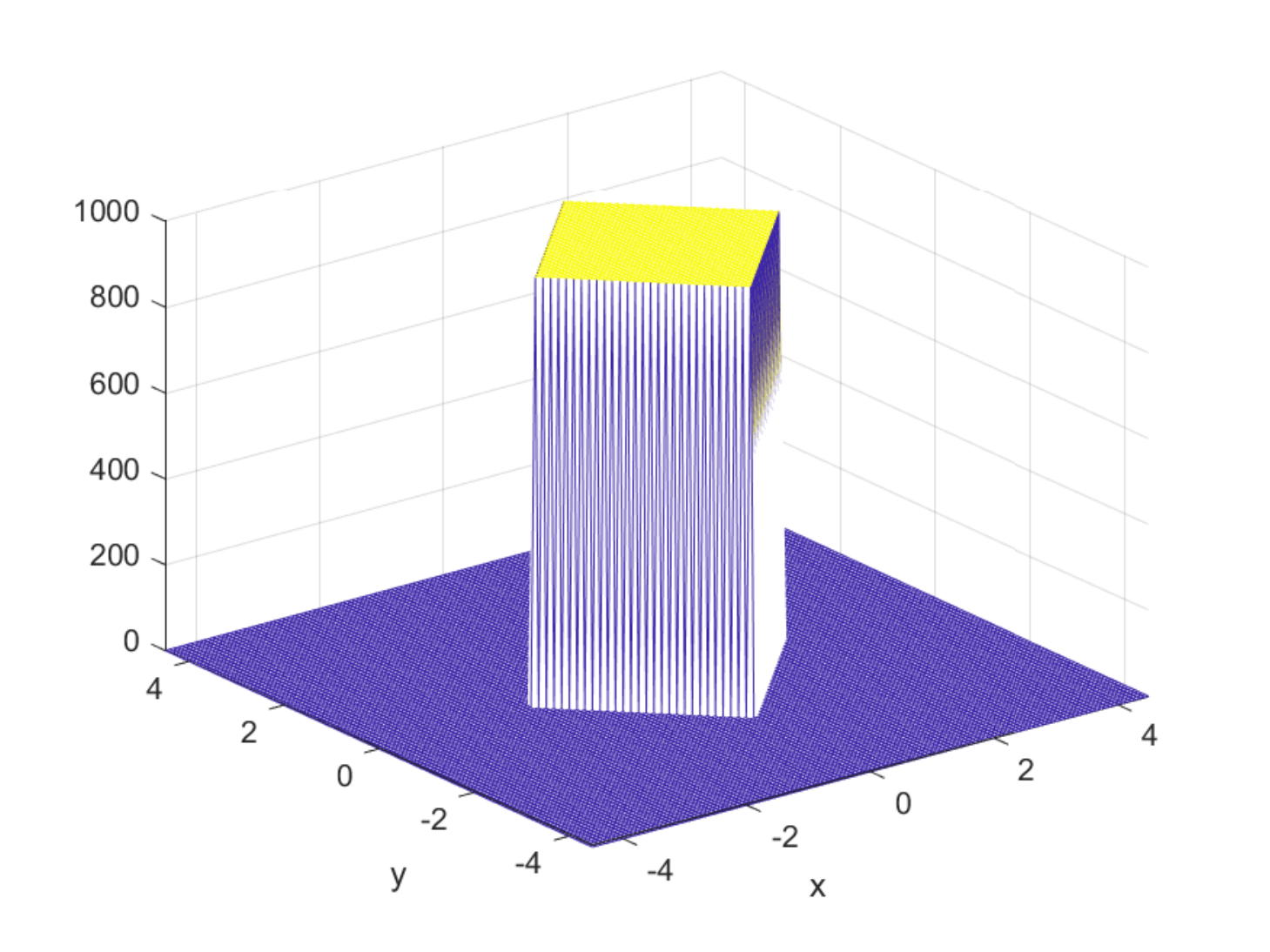}
	\end{subfigure}
	\begin{subfigure}[b]{0.25\textwidth}
		 \includegraphics[width=3.5cm,height=3.5cm]{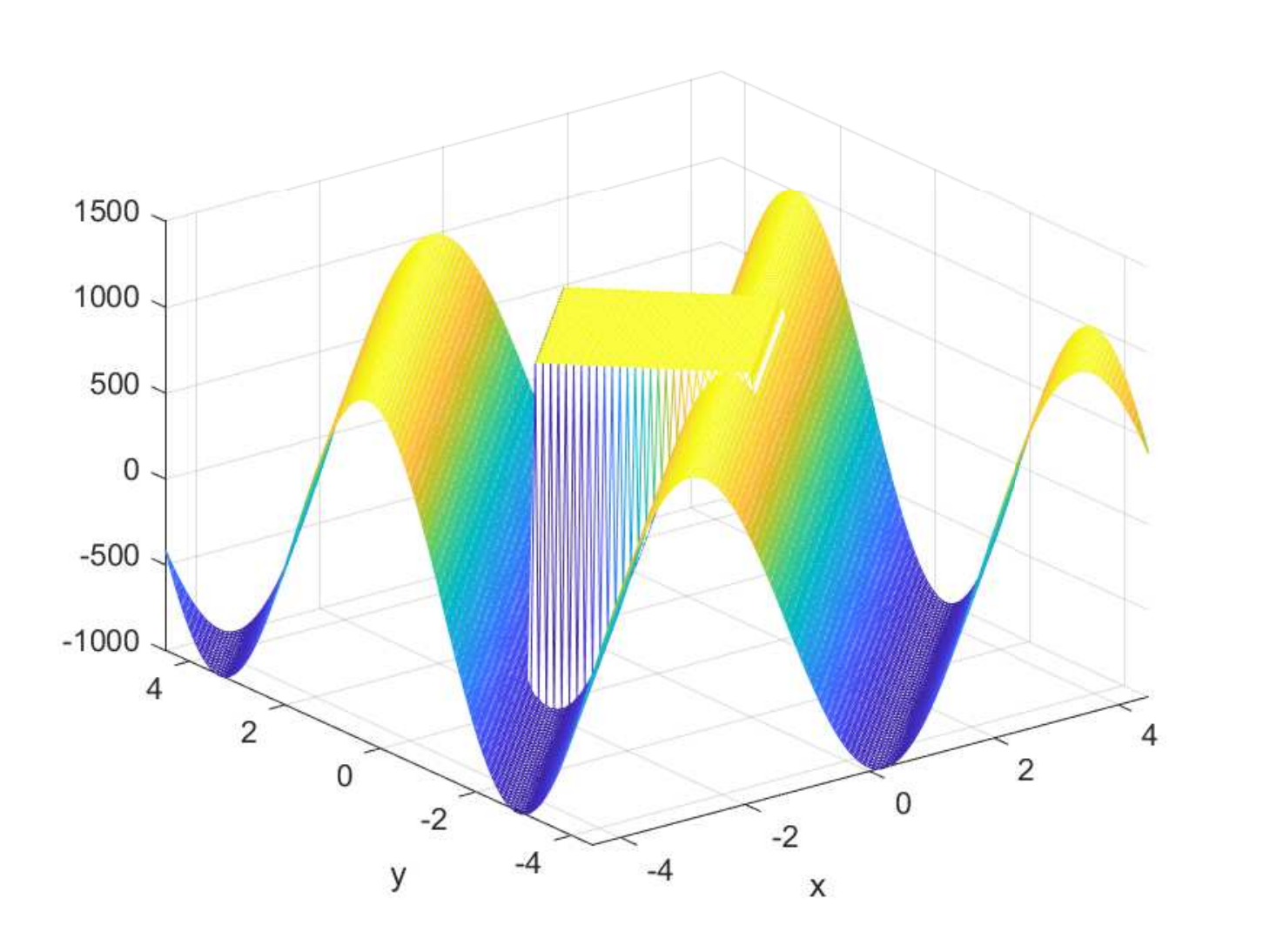}
	\end{subfigure}
	\begin{subfigure}[b]{0.25\textwidth}
		 \includegraphics[width=3.5cm,height=3.5cm]{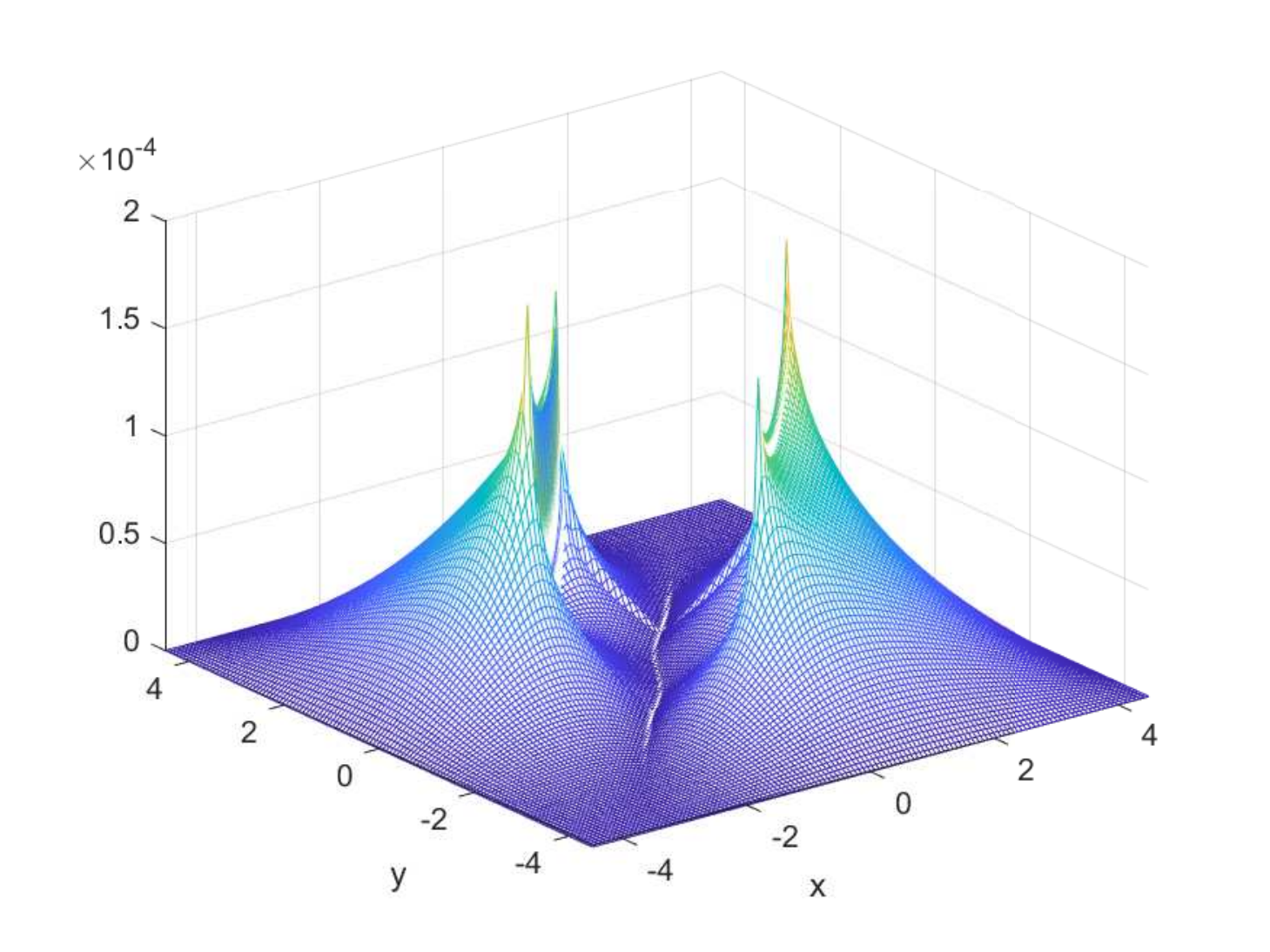}
	\end{subfigure}
	\caption
	{\cref{hybrid:ex5}: the interface curve $\Gamma_I$ (first panel), the coefficient $a(x,y)$ (second panel),  the numerical solution $u_h$ (third panel), and the error $|u_h-u|$ (fourth panel) with $h=2^{-7}\times 9$,  where $u_h$ is computed by our proposed hybrid finite difference scheme.}
	\label{hybrid:fig:QSp5}
\end{figure}

\begin{example}\label{hybrid:ex5}
	\normalfont
	Let $\Omega=(-4.5,4.5)^2$ and
	the interface curve be given by
	$\Gamma_I:=\{(x,y)\in \Omega:\; \psi(x,y)=0\}$ which is shown in \cref{hybrid:fig:QSp5}. Precisely, the sharp-edged interface is a square with 4  corner points $(-2,0)$, $(0,2)$, $(2,0)$ and $(0,-2)$.
	The functions in \eqref{Qeques2} are given by
	\begin{align*}
		&a_{+}=10^{-3},
		\qquad a_{-}=10^{3}, \qquad u_{+}=10^{3}\sin(x-y),
		\quad u_{-}=10^{-3}\cos(x)\cos(y)+1000,\\
		&  u(-4.5,y)= g_1, \qquad
		\qquad u(4.5,y)= g_2,\qquad \mbox{for} \qquad y\in(-4.5,4.5),\\
		& u(x,-4.5)= g_3, \qquad
		\qquad  u(x,4.5)= g_4, \qquad \mbox{for} \qquad x\in(-4.5,4.5),
	\end{align*}
	the other functions $f^{\pm}$, $g_D$, $g_N$, $g_1, \ldots,g_4$ in \eqref{Qeques2} can be obtained by plugging the above functions into \eqref{Qeques2}.  Clearly,
	$g_D$ and $g_N$ are not constants.
Note the high contrast $a_-/a_+=10^6$ on $\Gamma_I$.
	The numerical results are presented in \cref{hybrid:table:QSp5} and \cref{hybrid:fig:QSp5}.	
\end{example}

\begin{table}[htbp]
	\caption{Performance in \cref{hybrid:ex5}  of our proposed hybrid finite difference scheme and compact 9-point scheme on uniform Cartesian meshes with $h=2^{-J}\times 9$. $\kappa$ is the condition number of the coefficient matrix.}
\centering
\setlength{\tabcolsep}{0.5mm}{
	\begin{tabular}{c|c|c|c|c|c|c|c|c|c|c}
		\hline
		\multicolumn{1}{c|}{} &
		\multicolumn{5}{c|}{Our proposed hybrid scheme} &
		\multicolumn{5}{c}{Compact 9-point scheme} \\
		\cline{1-11}			
		$J$
		& $\frac{\|u_{h}-u\|_{2}}{\|u\|_{2}}$
		
		&order & $\|u_{h}-u\|_{\infty}$
		
		&order &  $\kappa$ & $\frac{\|u_{h}-u\|_{2}}{\|u\|_{2}}$
		
		&order & $\|u_{h}-u\|_{\infty}$
		
		&order &  $\kappa$  \\
		\hline
4    &7.431E-03    &0    &2.062E+01    &0    &1.337E+03    &6.254E-02    &0    &1.574E+02    &0    &1.238E+03\\
5    &4.505E-04    &4.0    &1.322E+00    &4.0    &1.020E+04    &1.110E-02    &2.5    &2.837E+01    &2.5    &6.529E+03\\
6    &5.701E-06    &6.3    &1.778E-02    &6.2    &6.394E+04    &6.953E-04    &4.0    &1.929E+00    &3.9    &4.152E+04\\
7    &4.937E-08    &6.9    &1.869E-04    &6.6    &3.920E+05    &2.993E-05    &4.5    &1.059E-01    &4.2    &3.286E+05\\
8    &6.087E-10    &6.3    &2.942E-06    &6.0    &2.132E+07    &1.155E-06    &4.7    &4.177E-03    &4.7    &1.474E+06\\
9    &    &    &    &    &    &8.390E-08    &3.8    &3.391E-04    &3.6    &1.006E+07\\			
		\hline
\end{tabular}}
	\label{hybrid:table:QSp5}
\end{table}

\begin{figure}[htbp]
	\centering
	\begin{subfigure}[b]{0.20\textwidth}
		 \includegraphics[width=3.0cm,height=3.0cm]{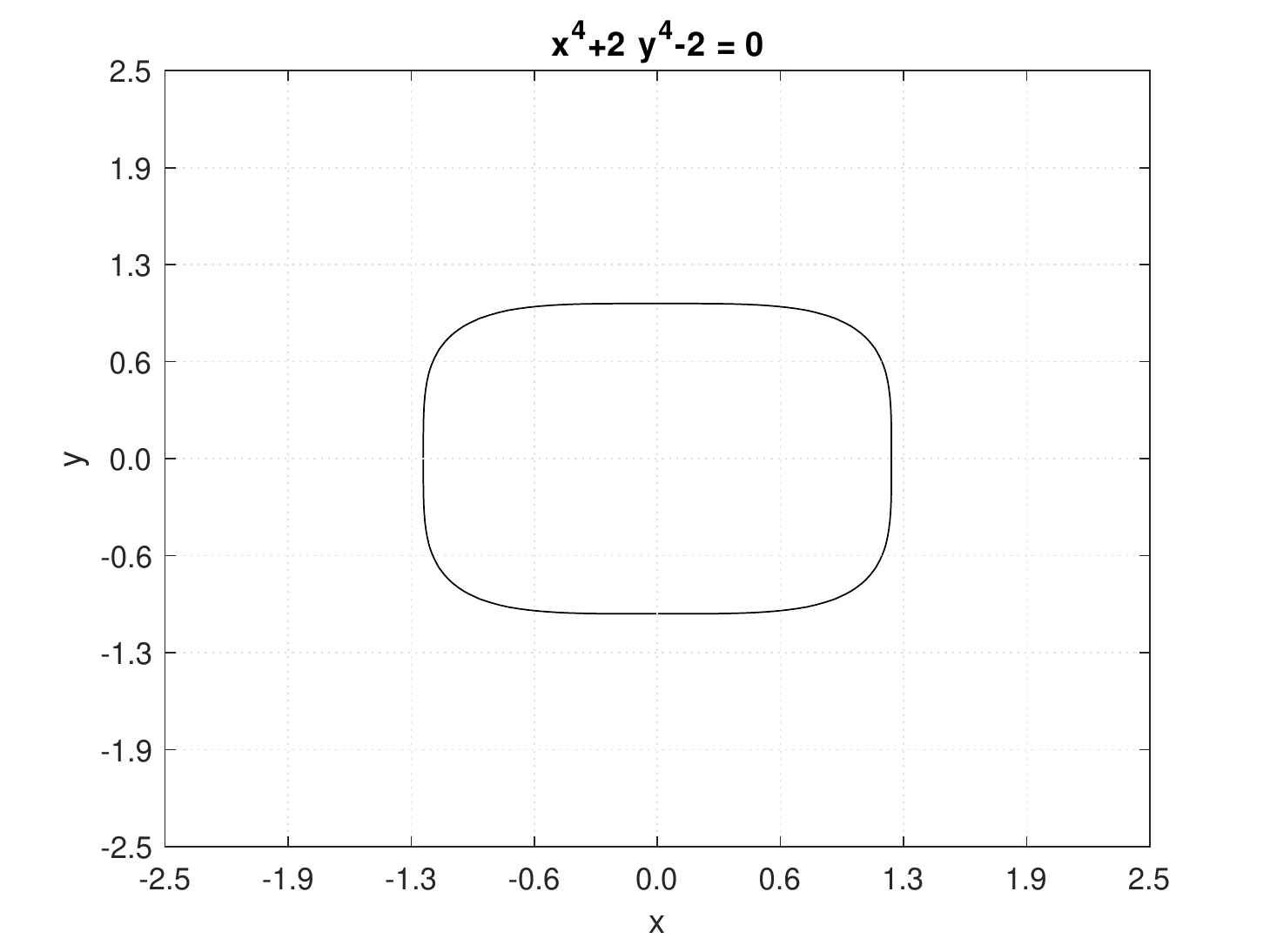}
	\end{subfigure}
	\begin{subfigure}[b]{0.25\textwidth}
		 \includegraphics[width=3.5cm,height=3.5cm]{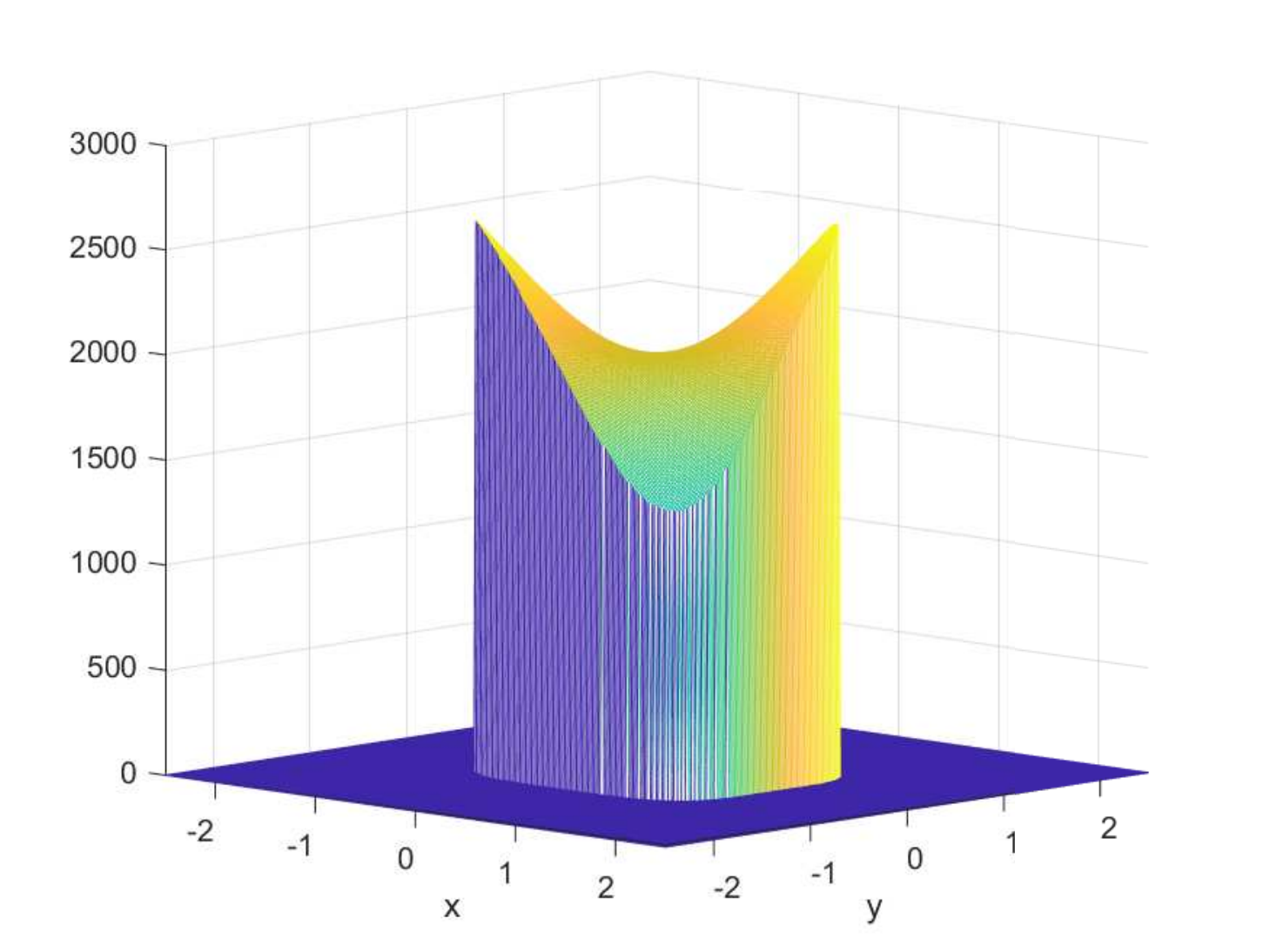}
	\end{subfigure}
	\begin{subfigure}[b]{0.25\textwidth}
		 \includegraphics[width=3.5cm,height=3.5cm]{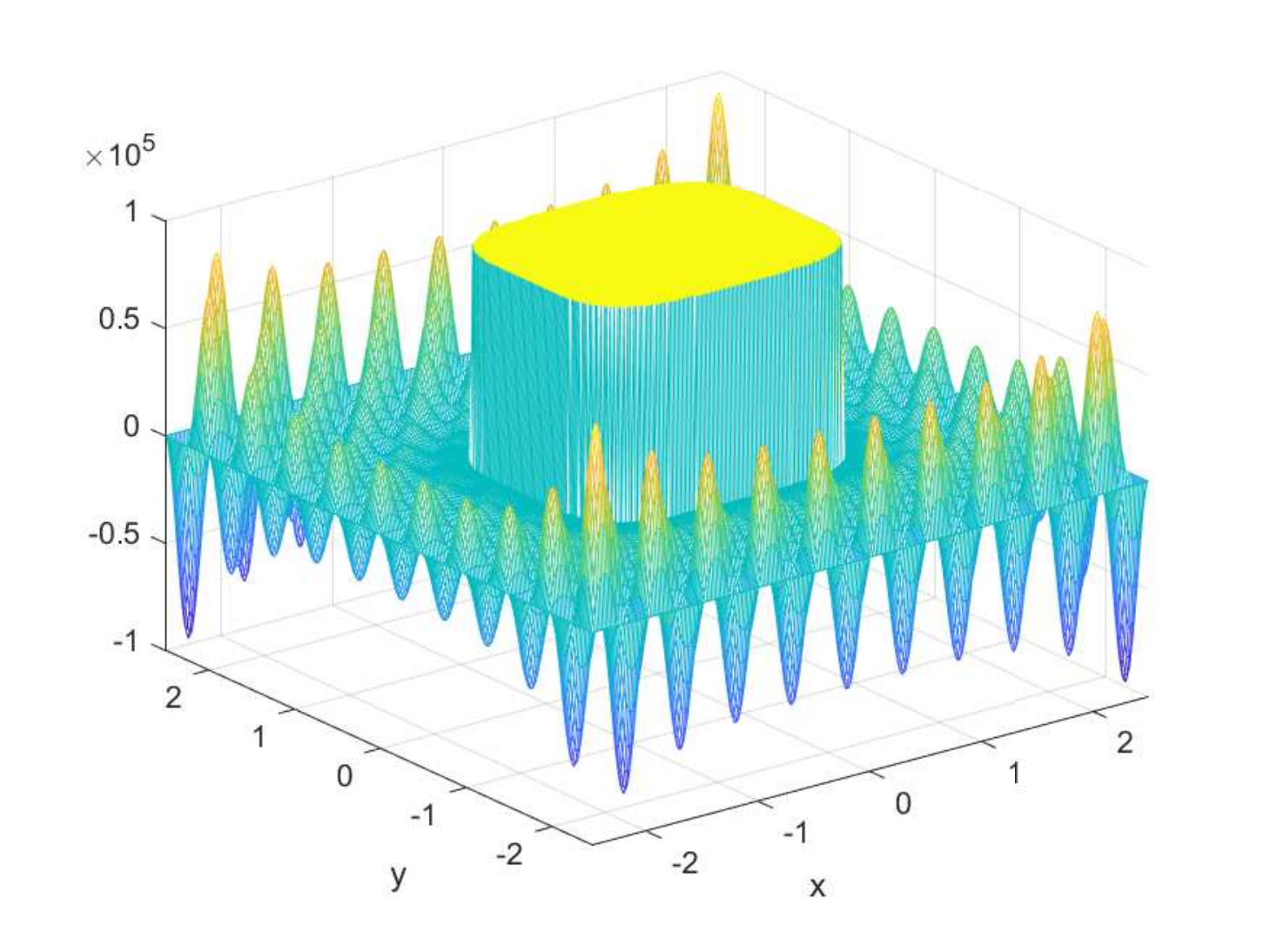}
	\end{subfigure}
	\begin{subfigure}[b]{0.25\textwidth}
		 \includegraphics[width=3.5cm,height=3.5cm]{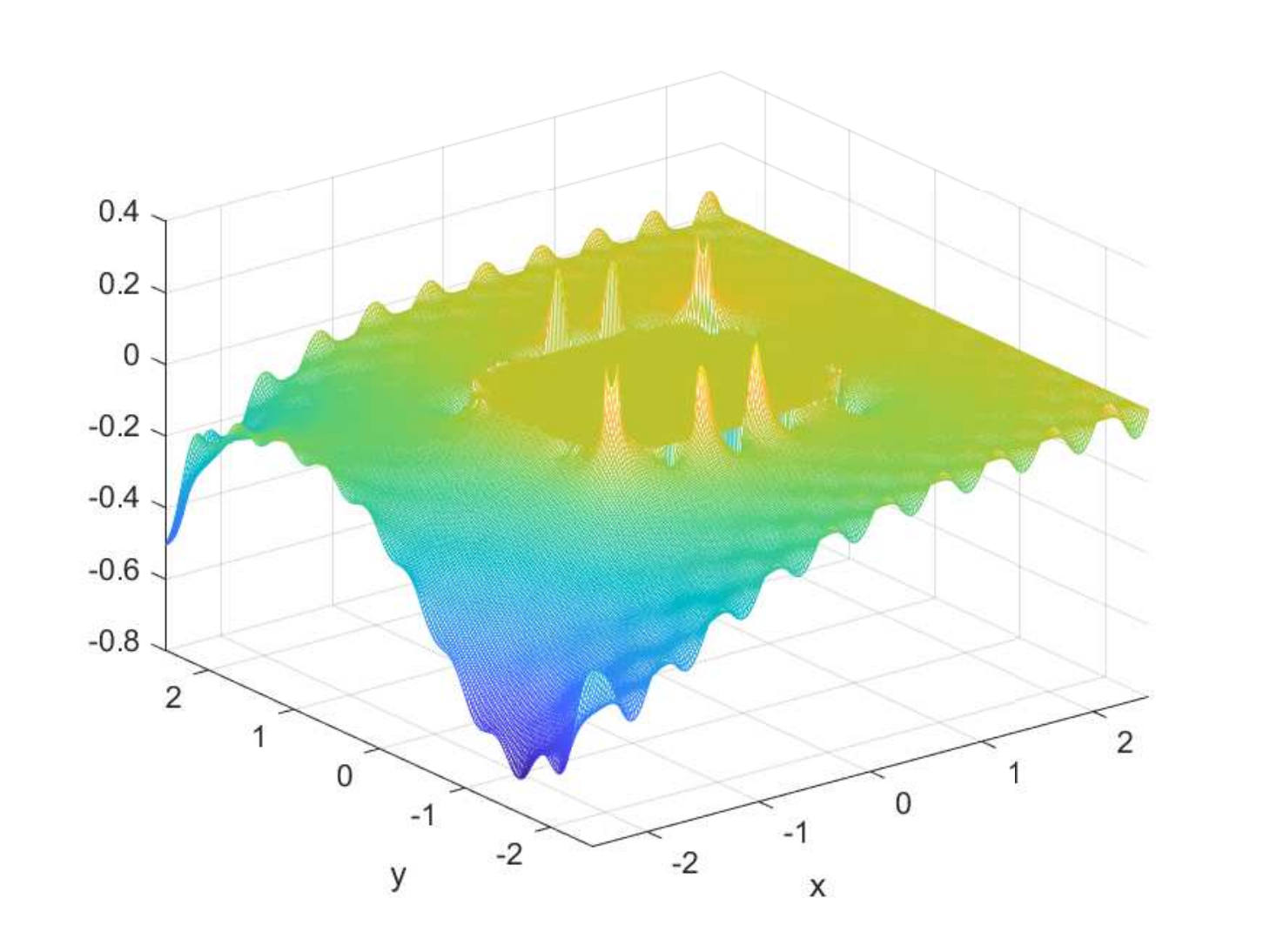}
	\end{subfigure}
	\caption
	{\cref{hybrid:ex1}: the interface curve $\Gamma_I$ (first panel), the coefficient $a(x,y)$ (second panel),  the numerical solution $u_h$ (third panel), and the error $u-u_h$ (fourth panel) with $h=2^{-8}\times 5$,  where $u_h$ is computed by our proposed hybrid finite difference scheme.}
	\label{hybrid:fig:QSp1}
\end{figure}

\begin{example}\label{hybrid:ex1}
	\normalfont
	Let $\Omega=(-2.5,2.5)^2$ and
	the interface curve be given by
	$\Gamma_I:=\{(x,y)\in \Omega:\; \psi(x,y)=0\}$ with
	$\psi (x,y)=x^4+2y^4-2$.
	The functions in \eqref{Qeques2} are given by
		\begin{align*}
		&a_{+}=10^{-3}(2+\sin(x)\sin(y)),
		\qquad a_{-}=10^3(2+\sin(x)\sin(y)), \qquad g_D=-10^5, \qquad g_N=0,\\
		&u_{+}=10^3\sin(4\pi x)\sin(4\pi y)(x^4+2y^4-2),
		\qquad u_{-}=10^{-3}\sin(4\pi x)\sin(4\pi y)(x^4+2y^4-2)+10^5,\\
		& -u_x(-2.5,y)+\alpha u(-2.5,y)= g_1, \qquad
		\qquad u(2.5,y)= g_2, \qquad \alpha=\sin(y),\qquad \mbox{for} \qquad y\in(-2.5,2.5),\\
		& -u_y(x,-2.5)= g_3, \qquad
		\qquad u_y(x,2.5)+\beta u(x,2.5)= g_4, \qquad \beta=\cos(x), \qquad \mbox{for} \qquad x\in(-2.5,2.5),
	\end{align*}
 the other functions $f^{\pm}$, $g_1, \ldots,g_4$ in \eqref{Qeques2} can be obtained by plugging the above functions into \eqref{Qeques2}.
Note the high contrast $a_-/a_+=10^6$ on $\Gamma_I$.
	The numerical results are presented in \cref{hybrid:table:QSp1} and \cref{hybrid:fig:QSp1}.	
\end{example}

\begin{table}[htbp]
	\caption{Performance in \cref{hybrid:ex1}  of our proposed hybrid finite difference scheme on uniform Cartesian meshes with $h=2^{-J}\times 5$. }
	\centering
	\setlength{\tabcolsep}{1.5mm}{
		\begin{tabular}{c|c|c|c|c|c|c|c|c}
			\hline
			$J$
			& $\frac{\|u_{h}-u\|_{2}}{\|u\|_{2}}$ 
			
			&order & $\|u_{h}-u\|_{\infty}$
			
			&order &  ${\|u_{h}-u_{h/2}\|_{2}}$
			&order &  $\|u_{h}-u_{h/2}\|_{\infty}$
			
			&order \\
			\hline
5   &8.167E-01   & 0  &1.758E+05   & 0  &1.811E+05   & 0  &1.734E+05   & 0  \\
6   &1.123E-02   &6.2   &2.488E+03   &6.1   &2.471E+03   &6.2   &2.441E+03   &6.2   \\
7   &2.059E-04   &5.8   &4.711E+01   &5.7   &4.550E+01   &5.8   &4.640E+01   &5.7   \\
8   &3.035E-06   &6.1   &7.028E-01   &6.1   &6.701E-01   &6.1   &6.919E-01   &6.1   \\
9   &4.632E-08   &6.0   &1.087E-02   &6.0   &9.946E-03   &6.1   &1.037E-02   &6.1   \\
			\hline
	\end{tabular}}
	\label{hybrid:table:QSp1}
\end{table}

\begin{figure}[htbp]
	\centering
	\begin{subfigure}[b]{0.20\textwidth}
		 \includegraphics[width=3.0cm,height=3.0cm]{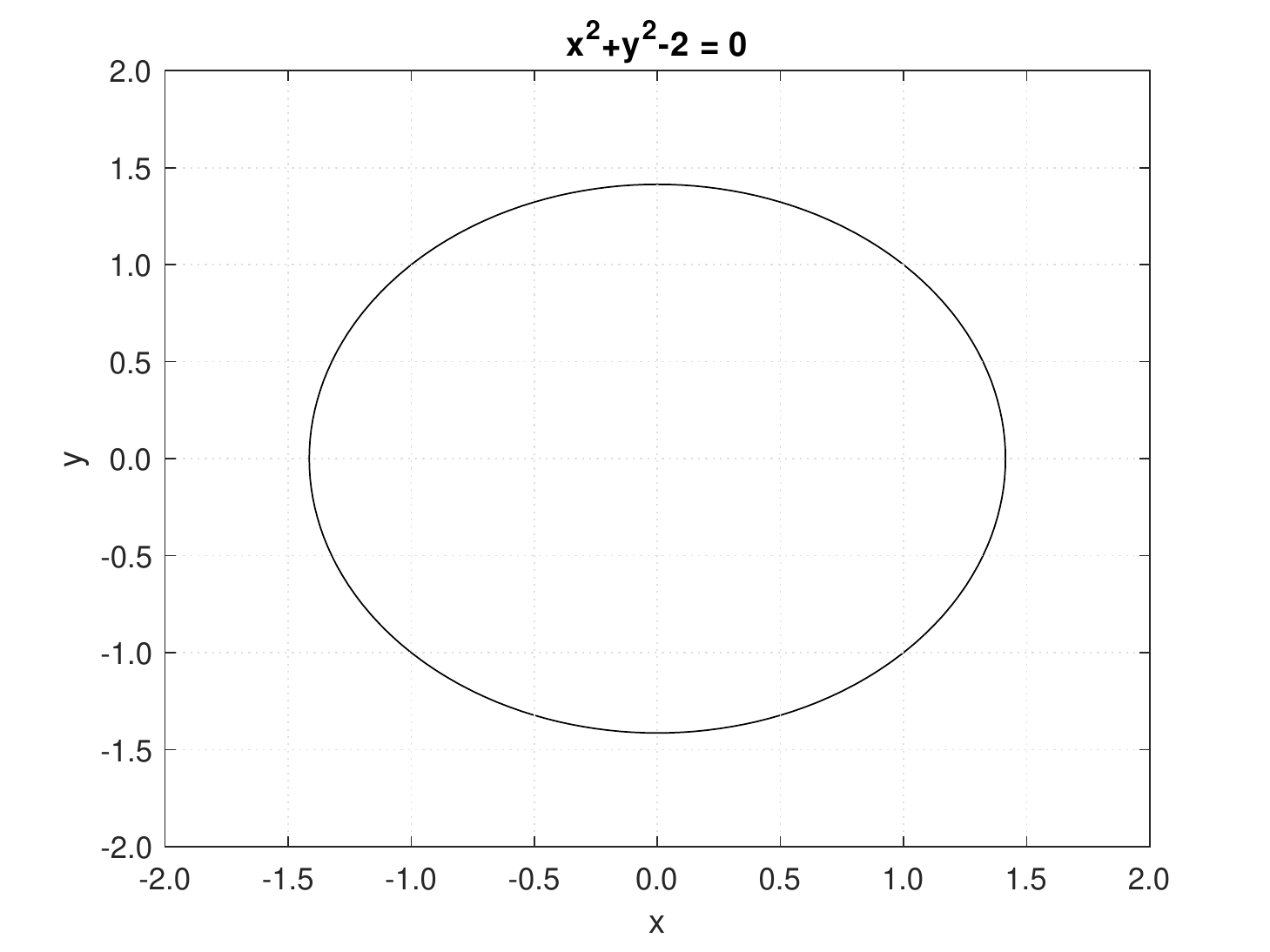}
	\end{subfigure}
	\begin{subfigure}[b]{0.25\textwidth}
		 \includegraphics[width=3.5cm,height=3.5cm]{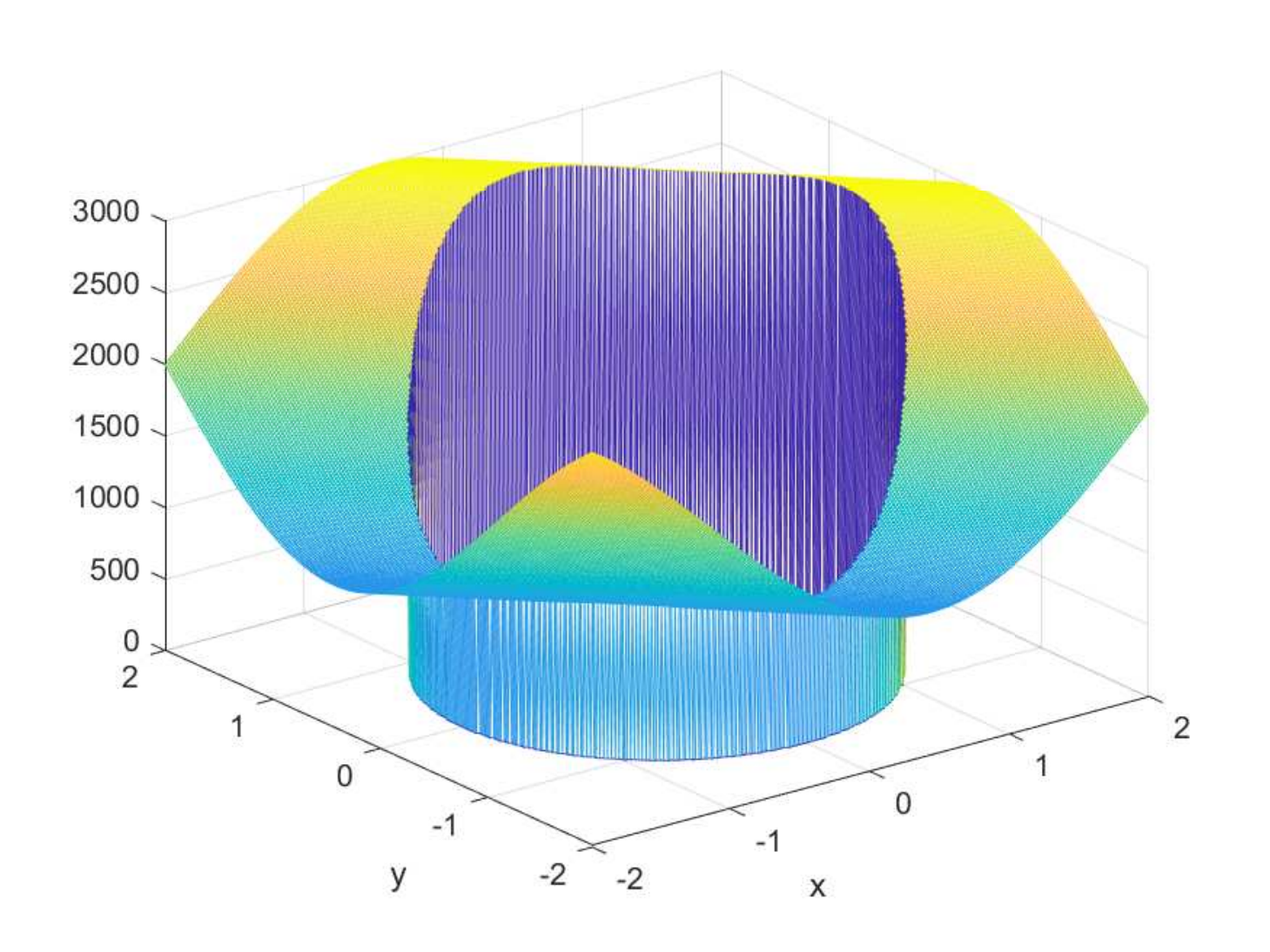}
	\end{subfigure}
	\begin{subfigure}[b]{0.25\textwidth}
		 \includegraphics[width=3.5cm,height=3.5cm]{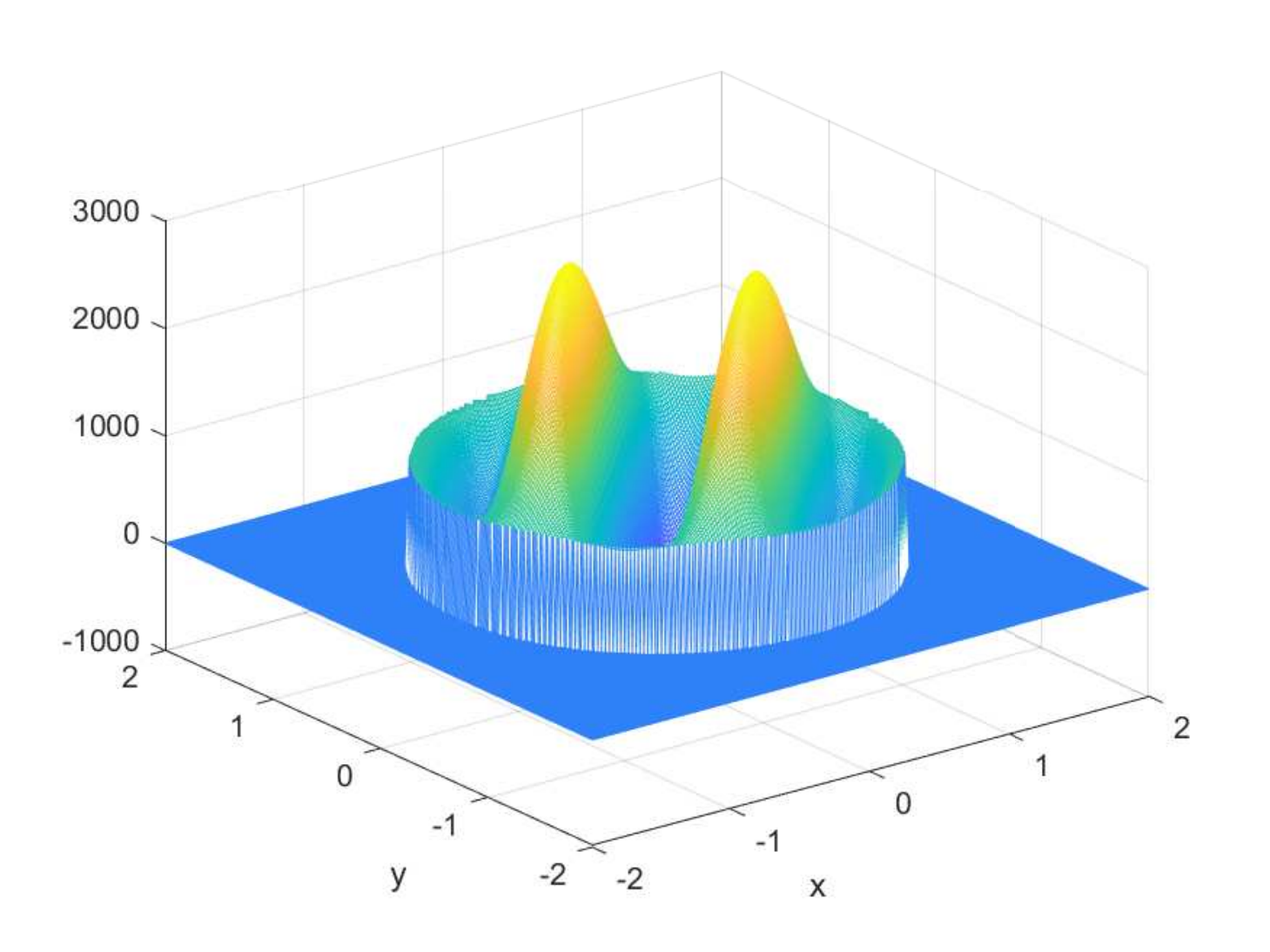}
	\end{subfigure}
	\begin{subfigure}[b]{0.25\textwidth}
		 \includegraphics[width=3.5cm,height=3.5cm]{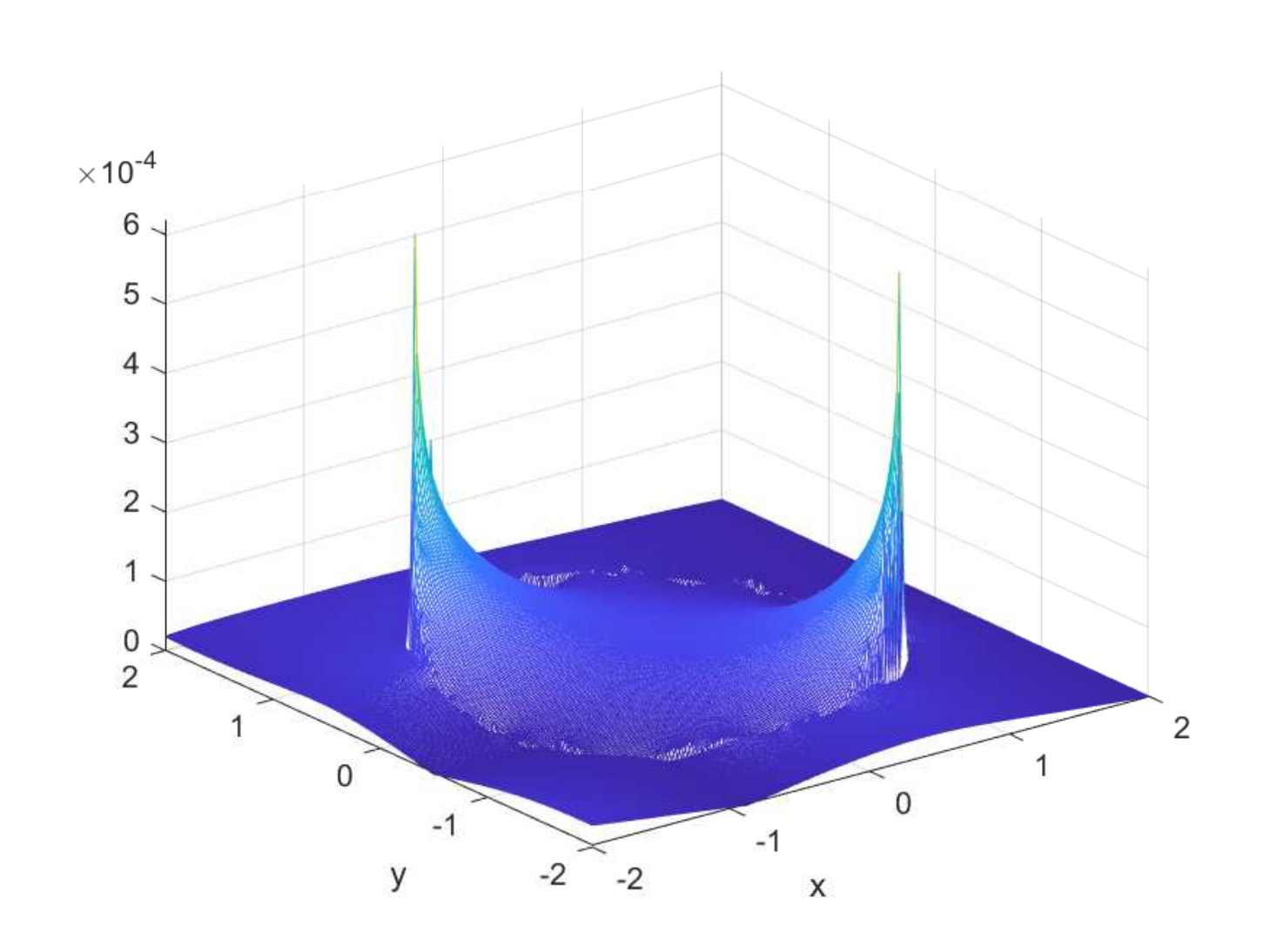}
	\end{subfigure}
	\caption
	{\cref{hybrid:ex2}: the interface curve $\Gamma_I$ (first panel), the coefficient $a(x,y)$ (second panel),  the numerical solution $u_h$ (third panel), and the error $|u_h-u|$ (fourth panel) with $h=2^{-8}\times 4$,  where $u_h$ is computed by our proposed hybrid finite difference scheme.}
	\label{hybrid:fig:QSp2}
\end{figure}

\begin{example}\label{hybrid:ex2}
	\normalfont
	Let $\Omega=(-2,2)^2$ and
	the interface curve be given by
	$\Gamma_I:=\{(x,y)\in \Omega:\; \psi(x,y)=0\}$ with
	$\psi (x,y)=x^2+y^2-2$.
	The functions in \eqref{Qeques2} are given by
	\begin{align*}
		&a_{+}=10^3(2+\sin(x+y)),
		\qquad a_{-}=10^{-3}(2+\sin(x+y)), \qquad g_D=-10^3, \qquad g_N=0,\\
		&u_{+}=10^{-3}\cos(4 (x-y))(x^2+y^2-2),
		\qquad u_{-}=10^{3}\cos(4 (x-y))(x^2+y^2-2)+10^3,\\
		& -u_x(-2,y)+\alpha u(-2,y)= g_1, \qquad
		\qquad u(2,y)= g_2, \qquad \alpha=\sin(y),\qquad \mbox{for} \qquad y\in(-2,2),\\
		& -u_y(x,-2)= g_3, \qquad
		\qquad u_y(x,2)+\beta u(x,2)= g_4, \qquad \beta=\cos(x), \qquad \mbox{for} \qquad x\in(-2,2),
	\end{align*}
	the other functions $f^{\pm}$, $g_1, \ldots,g_4$ in \eqref{Qeques2} can be obtained by plugging the above functions into \eqref{Qeques2}.
Note the high contrast $a_+/a_-=10^6$ on $\Gamma_I$.
	The numerical results are presented in \cref{hybrid:table:QSp2} and \cref{hybrid:fig:QSp2}.	
\end{example}

\begin{table}[htbp]
	\caption{Performance in \cref{hybrid:ex2}  of our proposed hybrid finite difference scheme  on uniform Cartesian meshes with $h=2^{-J}\times 4$. }
	\centering
	\setlength{\tabcolsep}{1.5mm}{
		\begin{tabular}{c|c|c|c|c|c|c|c|c}
			\hline
			$J$
			& $\frac{\|u_{h}-u\|_{2}}{\|u\|_{2}}$ 
			
			&order & $\|u_{h}-u\|_{\infty}$
			
			&order &  ${\|u_{h}-u_{h/2}\|_{2}}$
			&order &  $\|u_{h}-u_{h/2}\|_{\infty}$
			
			&order  \\
			\hline
4   &8.087E-01   &0   &4.191E+03   &0   &2.568E+03   &0   &4.141E+03   &0   \\
5   &1.443E-02   &5.8   &1.061E+02   &5.3   &4.623E+01   &5.8   &1.048E+02   &5.3   \\
6   &2.679E-04   &5.8   &2.154E+00   &5.6   &8.629E-01   &5.7   &2.132E+00   &5.6   \\
7   &3.432E-06   &6.3   &3.518E-02   &5.9   &1.100E-02   &6.3   &3.477E-02   &5.9   \\
8   &6.625E-08   &5.7   &6.192E-04   &5.8   &2.120E-04   &5.7   &6.118E-04   &5.8   \\
			\hline
	\end{tabular}}
	\label{hybrid:table:QSp2}
\end{table}

\begin{figure}[htbp]
	\centering
	\begin{subfigure}[b]{0.20\textwidth}
		 \includegraphics[width=3.0cm,height=3.0cm]{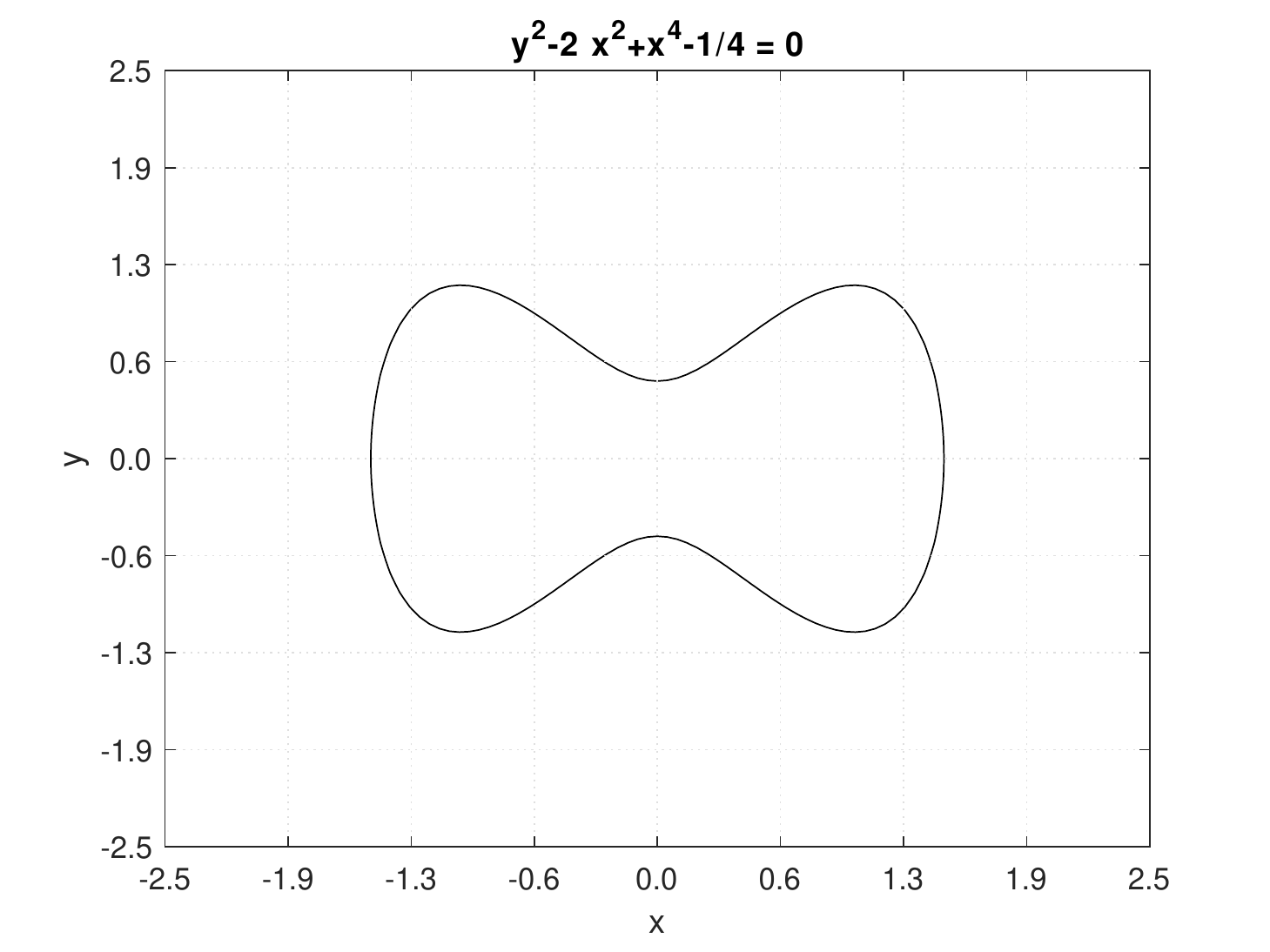}
	\end{subfigure}
	\begin{subfigure}[b]{0.25\textwidth}
		 \includegraphics[width=3.5cm,height=3.5cm]{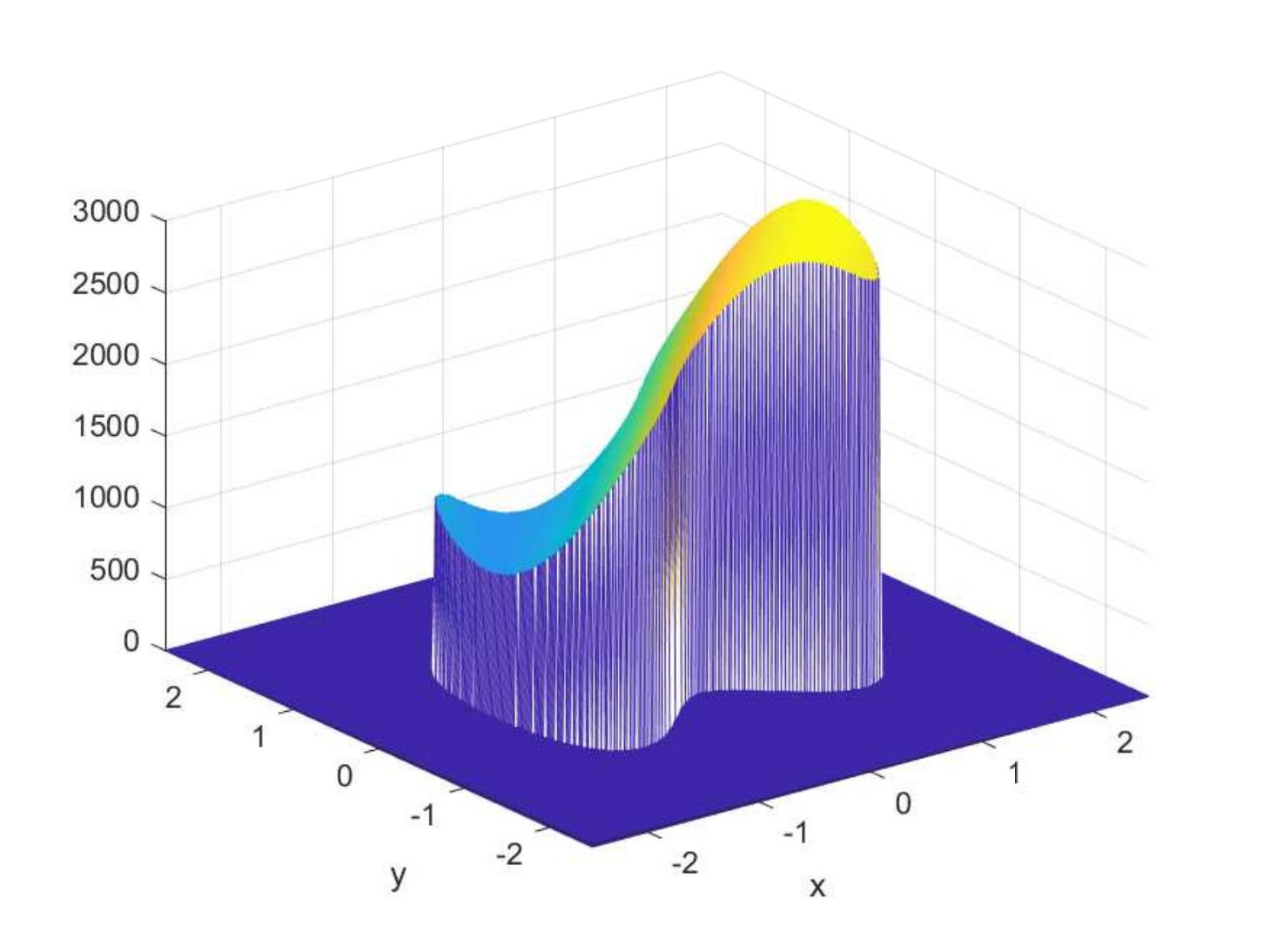}
	\end{subfigure}
	\begin{subfigure}[b]{0.25\textwidth}
		 \includegraphics[width=3.5cm,height=3.5cm]{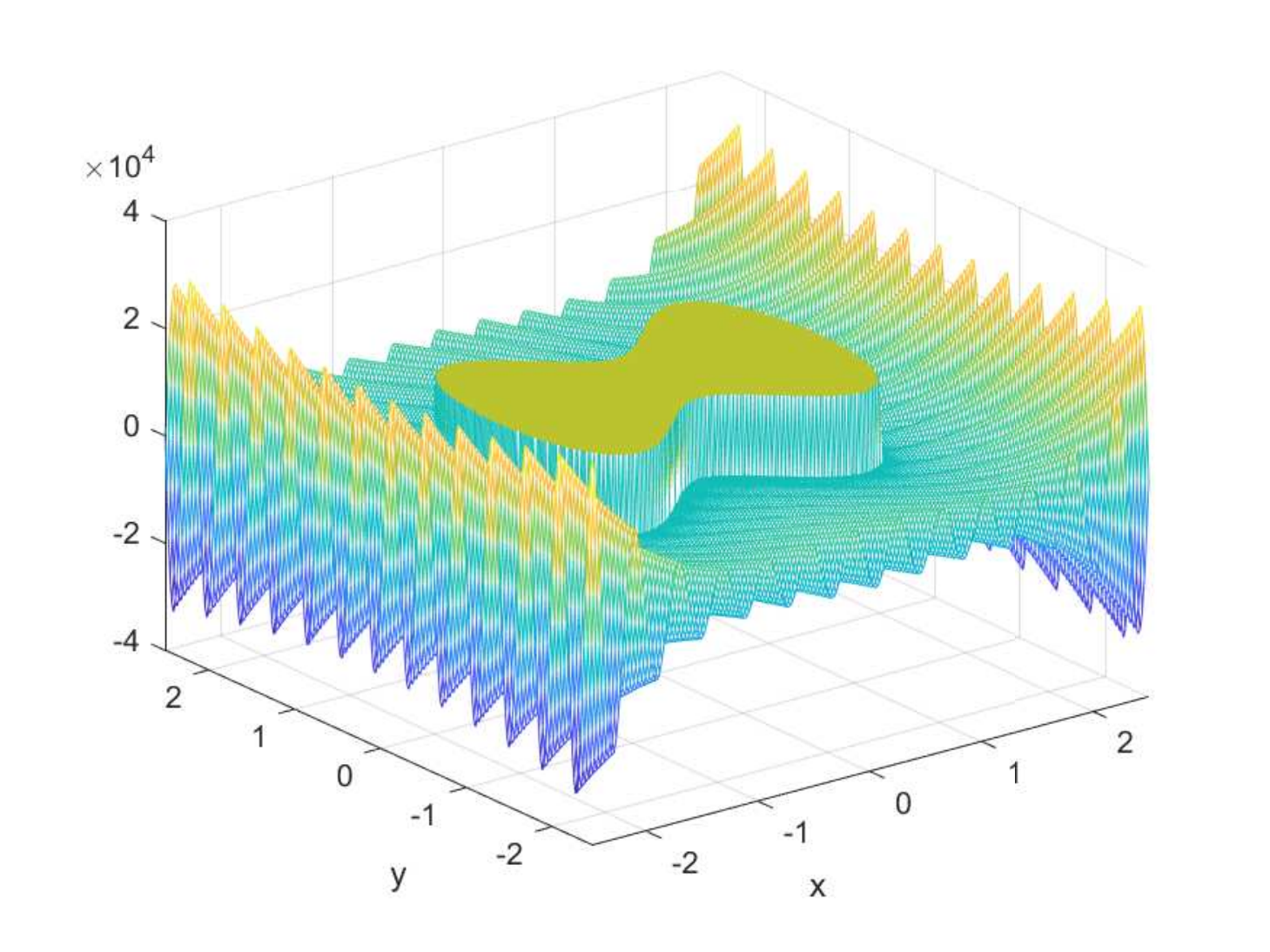}
	\end{subfigure}
	\begin{subfigure}[b]{0.25\textwidth}
		 \includegraphics[width=3.5cm,height=3.5cm]{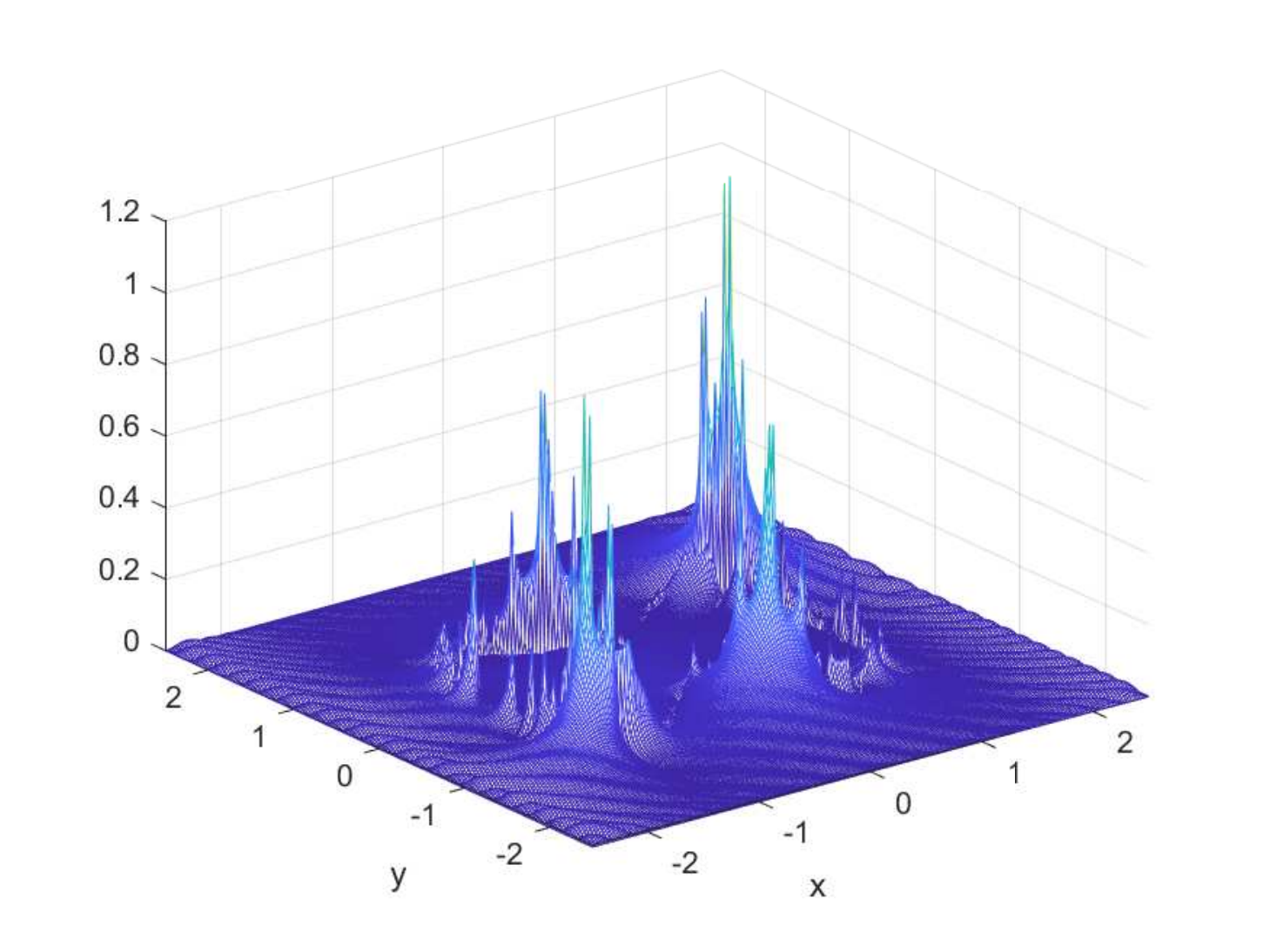}
	\end{subfigure}
	\caption
	{\cref{hybrid:ex4}: the interface curve $\Gamma_I$ (first panel), the coefficient $a(x,y)$ (second panel),  the numerical solution $u_h$ (third panel), and the error $|u_h-u|$ (fourth panel) with $h=2^{-8}\times 5$,  where $u_h$ is computed by our proposed hybrid finite difference scheme.}
	\label{hybrid:fig:QSp4}
\end{figure}

\begin{example}\label{hybrid:ex4}
	\normalfont
	Let $\Omega=(-2.5,2.5)^2$ and
	the interface curve be given by
	$\Gamma_I:=\{(x,y)\in \Omega:\; \psi(x,y)=0\}$ with
	$\psi (x,y)=y^2-2x^2+x^4-\frac{1}{4}$.
	The functions in \eqref{Qeques2} are given by
	\begin{align*}
		&a_{+}=10^{-3}(2+\sin(x-y)),
\qquad a_{-}=10^3(2+\sin(x-y)), \qquad g_D=-1.5\times 10^4, \qquad g_N=0,\\
&u_{+}=10^{3}\sin(16 (x+y))(y^2-2x^2+x^4-1/4),\\
&u_{-}=10^{-3}\sin(16 (x+y))(y^2-2x^2+x^4-1/4)+1.5\times 10^4,\\
&  u(-2.5,y)= g_1, \qquad
\qquad u(2.5,y)= g_2,\qquad \mbox{for} \qquad y\in(-2.5,2.5),\\
& u(x,-2.5)= g_3, \qquad
\qquad  u(x,2.5)= g_4, \qquad \mbox{for} \qquad x\in(-2.5,2.5),
	\end{align*}
	the other functions $f^{\pm}$, $g_1, \ldots,g_4$ in \eqref{Qeques2} can be obtained by plugging the above functions into \eqref{Qeques2}. Note the high contrast $a_-/a_+=10^6$ on $\Gamma_I$.
	The numerical results are presented in \cref{hybrid:table:QSp4} and \cref{hybrid:fig:QSp4}.	
\end{example}

\begin{table}[htbp]
	\caption{Performance in \cref{hybrid:ex4}  of our proposed hybrid finite difference scheme on uniform Cartesian meshes with $h=2^{-J}\times 5$.}
	\centering
	\setlength{\tabcolsep}{1.5mm}{
		\begin{tabular}{c|c|c|c|c|c|c|c|c}
			\hline
			$J$
			& $\frac{\|u_{h}-u\|_{2}}{\|u\|_{2}}$ 
			
			&order & $\|u_{h}-u\|_{\infty}$
			
			&order &  ${\|u_{h}-u_{h/2}\|_{2}}$
			&order &  $\|u_{h}-u_{h/2}\|_{\infty}$
			
			&order  \\
			\hline
5   &8.627E-01   &0   &9.480E+04   &0   &4.284E+04   &0   &9.338E+04   &0   \\
6   &2.854E-02   &4.9   &2.758E+03   &5.1   &1.360E+03   &5.0   &2.736E+03   &5.1   \\
7   &4.543E-04   &6.0   &5.673E+01   &5.6   &2.128E+01   &6.0   &5.658E+01   &5.6   \\
8   &6.195E-06   &6.2   &1.184E+00   &5.6   &2.856E-01   &6.2   &1.177E+00   &5.6   \\
9   &8.902E-08   &6.1   &1.738E-02   &6.1   &4.441E-03   &6.0   &1.788E-02   &6.0   \\
			\hline
	\end{tabular}}
	\label{hybrid:table:QSp4}
\end{table}

%
%
%

\subsection{Numerical examples with unknown $u$}

In this subsection, we provide five numerical examples with unknown $u$ of \eqref{Qeques2}. They can be characterized as follows.

\begin{itemize}	
\item In all examples, either $a_+/a_-$ or $a_-/a_+$ is very large on $\Gamma_I$ for high-contrast coefficients $a$.
	\item 4-side Dirichlet boundary conditions are demonstrated in  \cref{hybrid:unknown:ex1,hybrid:unknown:ex4}.
	\item 3-side Dirichlet and 1-side Robin boundary conditions
in \cref{hybrid:unknown:ex2,hybrid:unknown:ex3}.
	\item 1-side Dirichlet, 1-side Neumann and 2-side Robin boundary conditions
in \cref{hybrid:unknown:ex5}.
	
	\item All the interface curves $\Gamma_I$ are smooth and
	all the jump functions $g_D$ and $g_N$  are non-constant.
\end{itemize}

\begin{figure}[htbp]
	\centering
	\begin{subfigure}[b]{0.3\textwidth}
		 \includegraphics[width=5.5cm,height=3.5cm]{HyCUR1.pdf}
	\end{subfigure}
	\begin{subfigure}[b]{0.3\textwidth}
		 \includegraphics[width=5.5cm,height=3.5cm]{HyAA1.pdf}
	\end{subfigure}
	\begin{subfigure}[b]{0.3\textwidth}
		 \includegraphics[width=5.5cm,height=3.5cm]{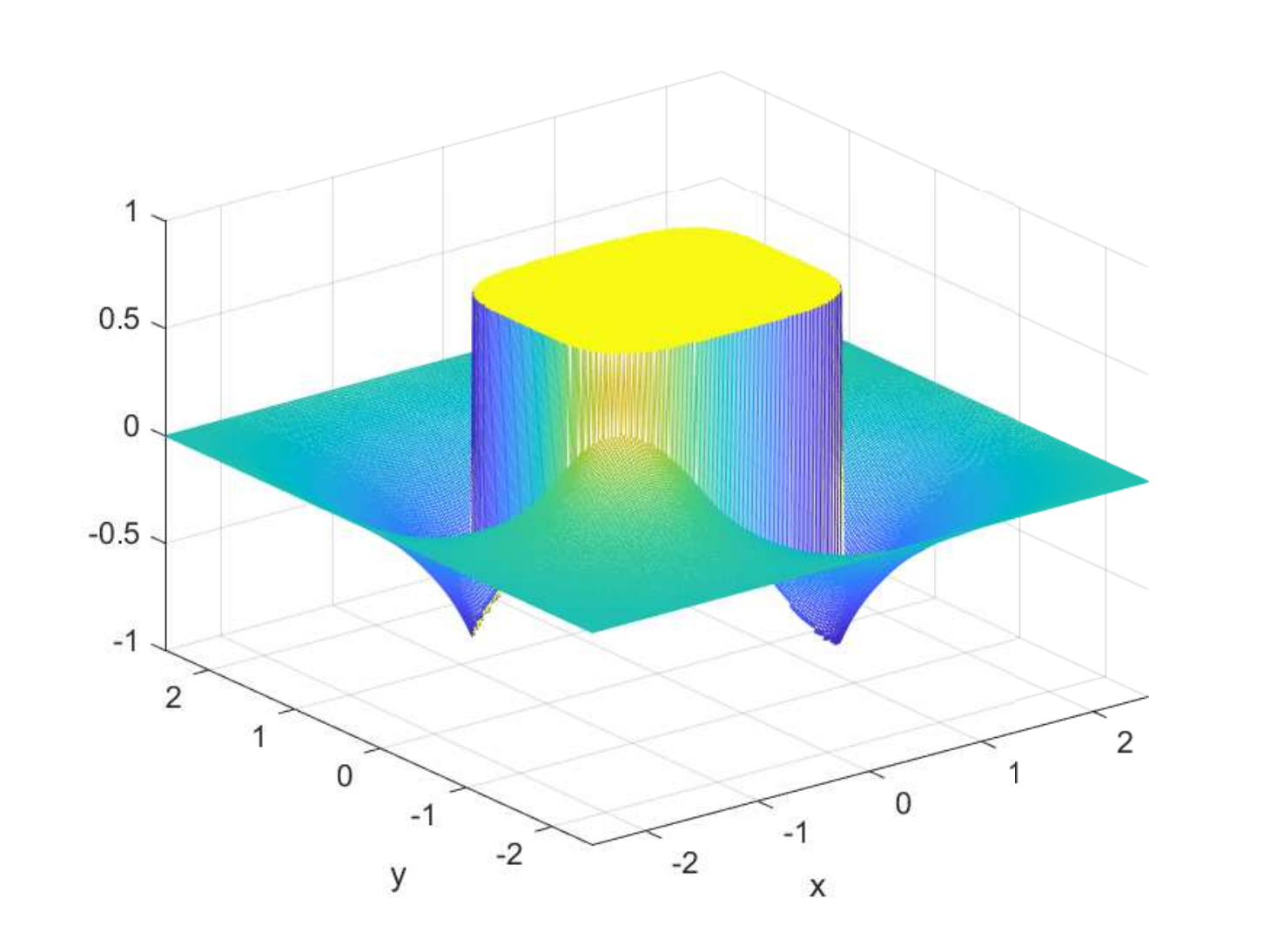}
	\end{subfigure}
	\caption
	{\cref{hybrid:unknown:ex1}: the interface curve $\Gamma_I$ (left), the coefficient $a(x,y)$ (middle) and the numerical solution $u_h$ (right) with $h=2^{-8}\times 5$,  where $u_h$ is computed by our proposed hybrid finite difference scheme. In order  to show the graph of $a(x,y)$ clearly, we rotate the graph of $a(x,y)$ by $\pi/2$ in this figure.}
	\label{hybrid:fig:unknown:QSp1}
\end{figure}

\begin{example}\label{hybrid:unknown:ex1}
	\normalfont
	Let $\Omega=(-2.5,2.5)^2$ and
	the interface curve be given by
	$\Gamma_I:=\{(x,y)\in \Omega:\; \psi(x,y)=0\}$ with
	$\psi (x,y)=x^4+2y^4-2$.
	The functions in \eqref{Qeques2} are given by
	\begin{align*}
		&a_{+}=2+\cos(x)\cos(y),
		\qquad a_{-}=10^3(2+\sin(x)\sin(y)), \qquad g_D=\sin(x)\sin(y)-1, \\
		&f_{+}=\sin(4\pi x)\sin(4\pi y),
		\qquad f_{-}=\cos(4\pi x)\cos(4\pi y), \qquad g_N=\cos(x)\cos(y),\\
		&  u(-2.5,y)= 0, \qquad
		\qquad u(2.5,y)= 0, \qquad \mbox{for} \qquad y\in(-2.5,2.5),\\
		& u(x,-2.5)= 0, \qquad
		\qquad  u(x,2.5)= 0,  \qquad \mbox{for} \qquad x\in(-2.5,2.5).
	\end{align*}
Note the high contrast $a_-/a_+\approx 10^3$ on $\Gamma_I$.
	The numerical results are presented in \cref{hybrid:table:unknown:QSp1} and \cref{hybrid:fig:unknown:QSp1}.	
\end{example}

\begin{table}[htbp]
	\caption{Performance in \cref{hybrid:unknown:ex1}  of our proposed hybrid finite difference scheme on uniform Cartesian meshes with $h=2^{-J}\times 5$. }
	\centering
	\setlength{\tabcolsep}{1.5mm}{
		\begin{tabular}{c|c|c|c|c}
			\hline
			$J$
		    &  ${\|u_{h}-u_{h/2}\|_{2}}$
			&order &  $\|u_{h}-u_{h/2}\|_{\infty}$

			&order  \\
			\hline
4   &9.83385E+02   &0   &3.29078E+02   &0  \\
5   &1.93678E+01   &5.7   &6.50631E+00   &5.7   \\
6   &3.13024E-01   &6.0   &1.04785E-01   &6.0   \\
8   &9.47776E-05   &5.8   &3.20754E-05   &5.8   \\
			\hline
	\end{tabular}}
	\label{hybrid:table:unknown:QSp1}
\end{table}

\begin{figure}[htbp]
	\centering
	\begin{subfigure}[b]{0.3\textwidth}
		 \includegraphics[width=5.5cm,height=3.5cm]{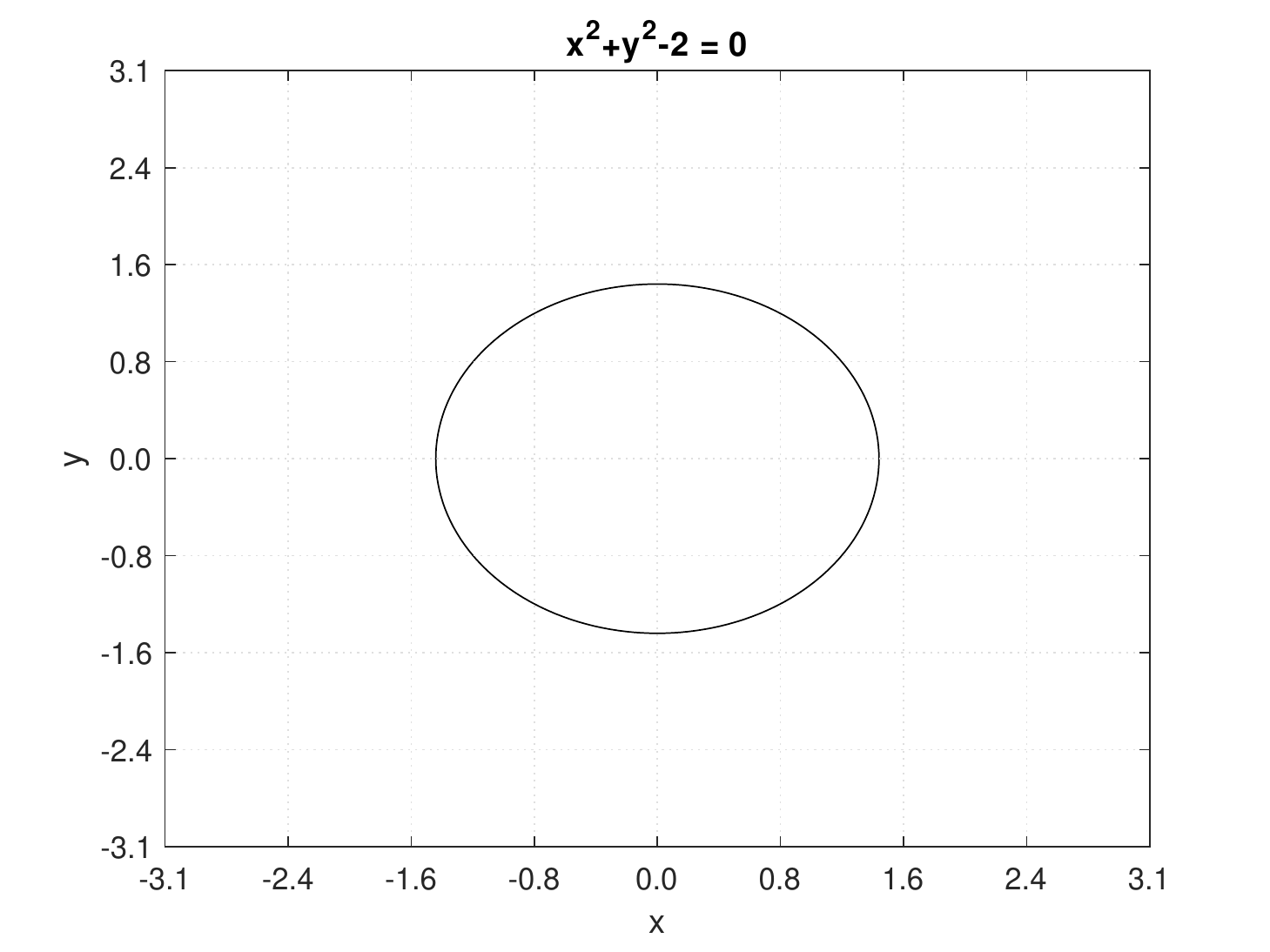}
	\end{subfigure}
	\begin{subfigure}[b]{0.3\textwidth}
		 \includegraphics[width=5.5cm,height=3.5cm]{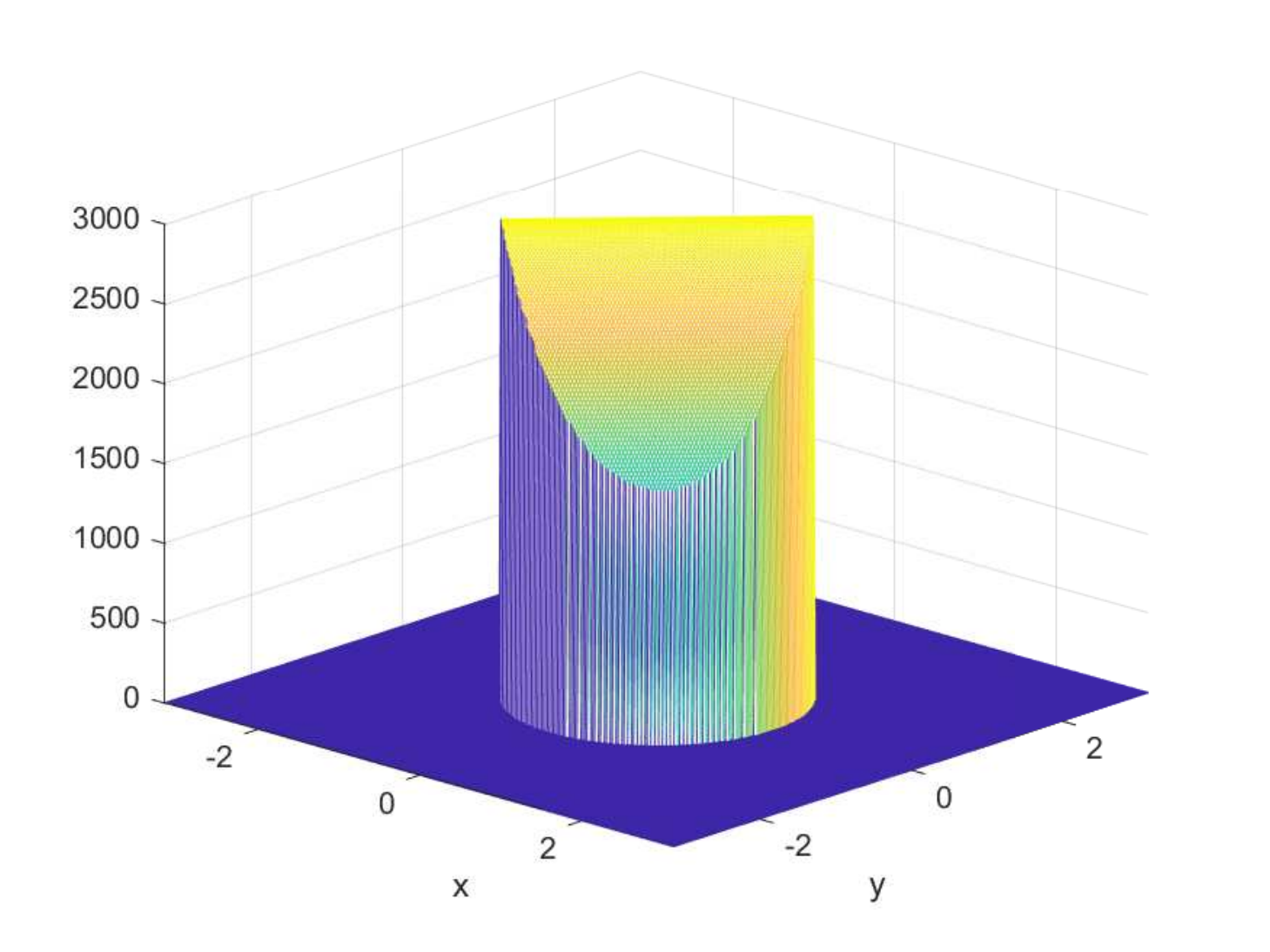}
	\end{subfigure}
	\begin{subfigure}[b]{0.3\textwidth}
		 \includegraphics[width=5.5cm,height=3.5cm]{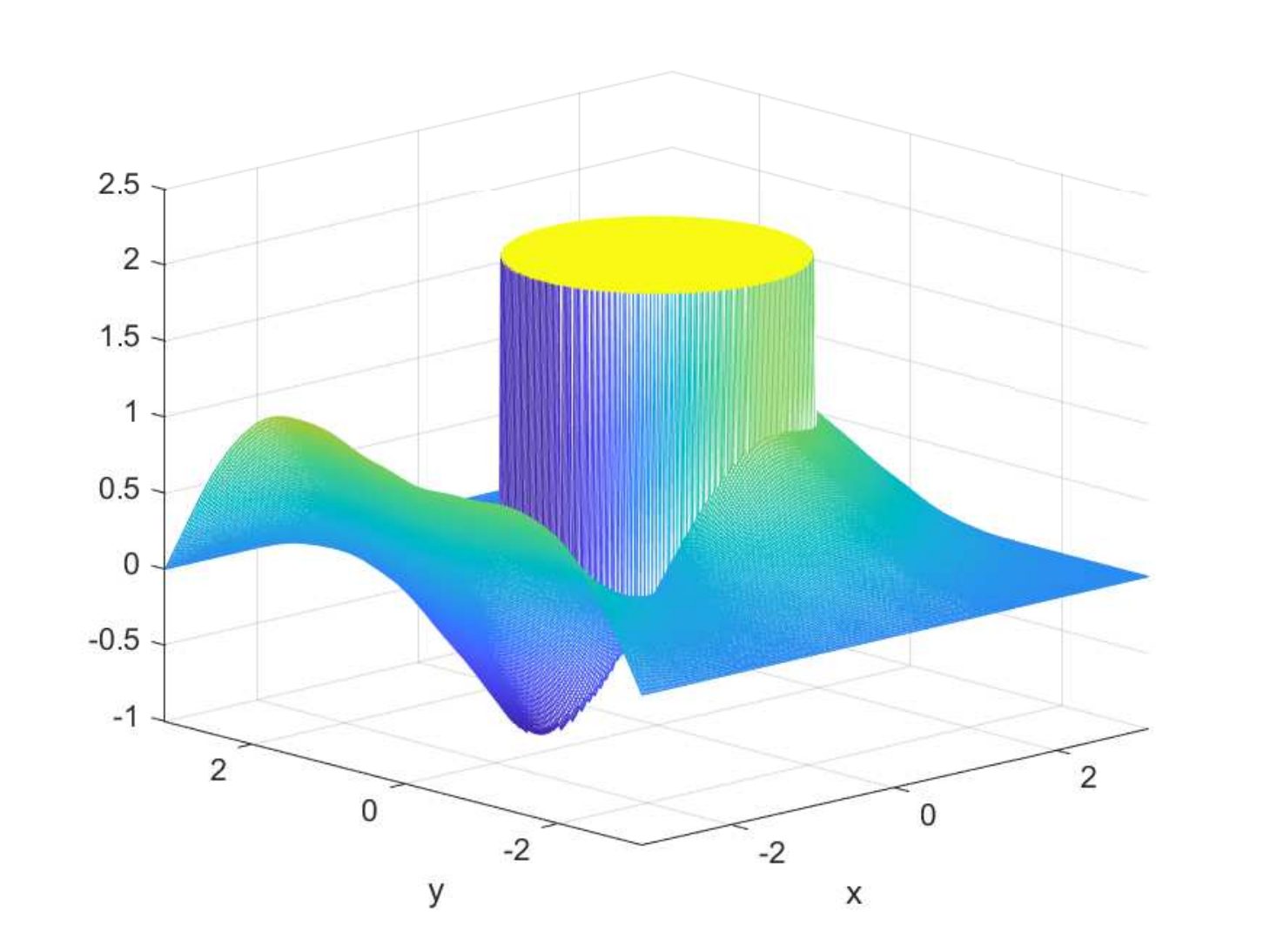}
	\end{subfigure}
	\caption
	{\cref{hybrid:unknown:ex2}: the interface curve $\Gamma_I$ (left), the coefficient $a(x,y)$ (middle) and the numerical solution $u_h$ (right) with $h=2^{-8}\times 2\pi$,  where $u_h$ is computed by our proposed hybrid finite difference scheme. In order  to show the graph of $a(x,y)$ clearly, we rotate the graph of $a(x,y)$ by $\pi/2$ in this figure.}
	\label{hybrid:fig:unknown:QSp2}
\end{figure}

\begin{example}\label{hybrid:unknown:ex2}
	\normalfont
	Let $\Omega=(-\pi,\pi)^2$ and
	the interface curve be given by
	$\Gamma_I:=\{(x,y)\in \Omega:\; \psi(x,y)=0\}$ with
	$\psi (x,y)=x^2+y^2-2$.
	The functions in \eqref{Qeques2} are given by
	\begin{align*}
		&a_{+}=2+\cos(x-y),
		\qquad a_{-}=10^3(2+\cos(x-y)), \qquad g_D=\sin(x-y)-2, \\
		&f_{+}=\sin(8 x)\sin(8 y),
		\qquad f_{-}=\cos(8 x)\cos(8 y), \qquad g_N=\cos(x+y),\\
		&  -u_x(-\pi,y)+\cos(y) u(-\pi,y)= \cos(y)+1, \qquad
		\qquad u(\pi,y)= 0, \qquad \mbox{for} \qquad y\in(-\pi,\pi),\\
		& u(x,-\pi)= 0, \qquad
		\qquad  u(x,\pi)= 0,  \qquad \mbox{for} \qquad x\in(-\pi,\pi).
	\end{align*}
Note the high contrast $a_-/a_+=10^3$ on $\Gamma_I$.
	The numerical results are presented in \cref{hybrid:table:unknown:QSp2} and \cref{hybrid:fig:unknown:QSp2}.	
\end{example}

\begin{table}[htbp]
	\caption{Performance in \cref{hybrid:unknown:ex2}  of our proposed hybrid finite difference scheme on uniform Cartesian meshes with $h=2^{-J}\times 2\pi$. }
	\centering
	\setlength{\tabcolsep}{1.5mm}{
		\begin{tabular}{c|c|c|c|c}
			\hline
			$J$
			&  ${\|u_{h}-u_{h/2}\|_{2}}$
			&order &  $\|u_{h}-u_{h/2}\|_{\infty}$

			&order \\
			\hline
4   &7.02037E+02   &0   &1.84708E+02   &0   \\
5   &9.69424E+00   &6.2   &2.54978E+00   &6.2   \\
6   &2.26556E-01   &5.4   &5.97145E-02   &5.4   \\
7   &2.57284E-03   &6.5   &6.79725E-04   &6.5   \\
8   &5.07886E-05   &5.7   &1.34801E-05   &5.7   \\
			\hline
	\end{tabular}}
	\label{hybrid:table:unknown:QSp2}
\end{table}

\begin{figure}[htbp]
	\centering
	\begin{subfigure}[b]{0.3\textwidth}
		 \includegraphics[width=5.5cm,height=3.5cm]{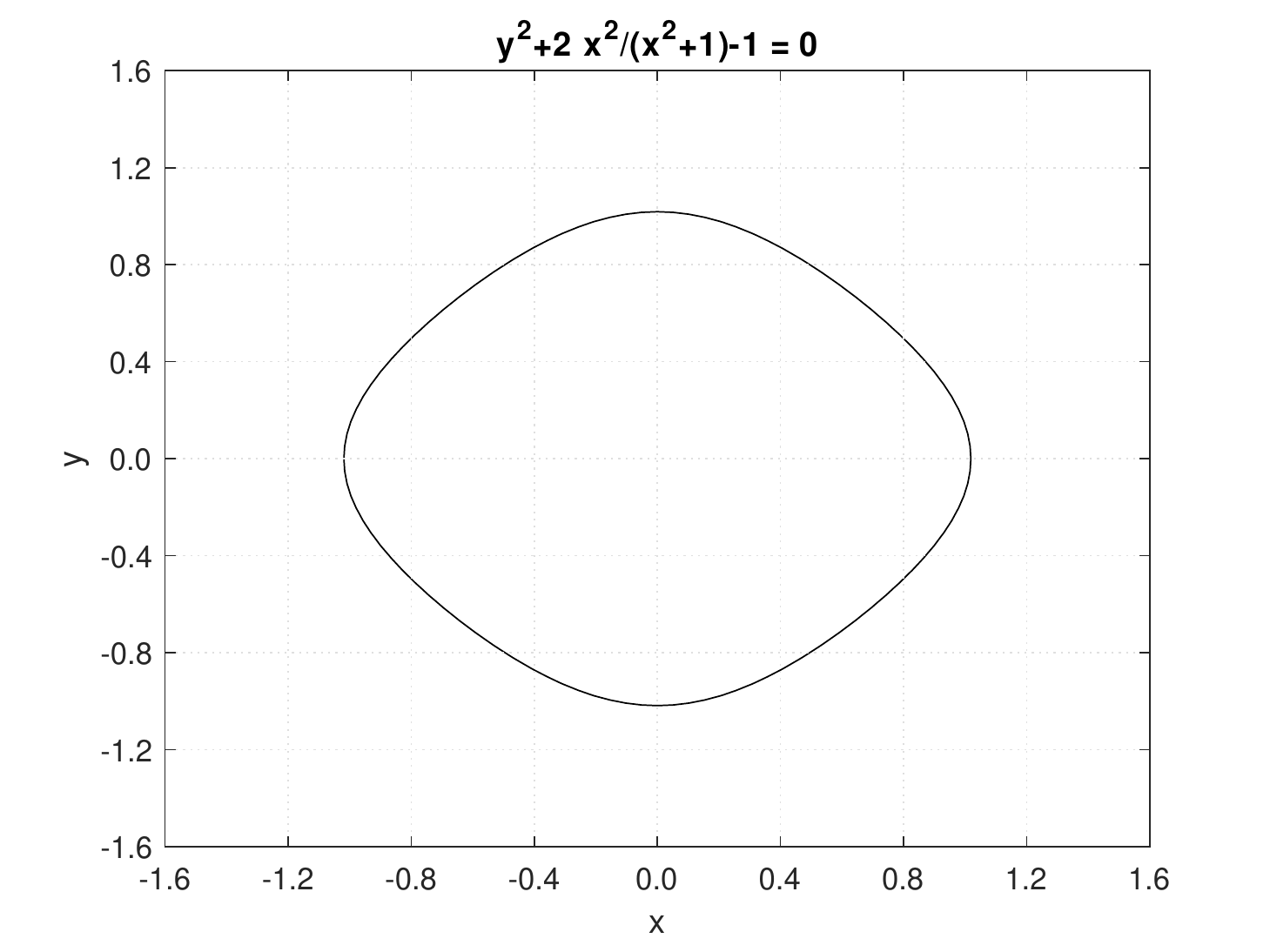}
	\end{subfigure}
	\begin{subfigure}[b]{0.3\textwidth}
		 \includegraphics[width=5.5cm,height=3.5cm]{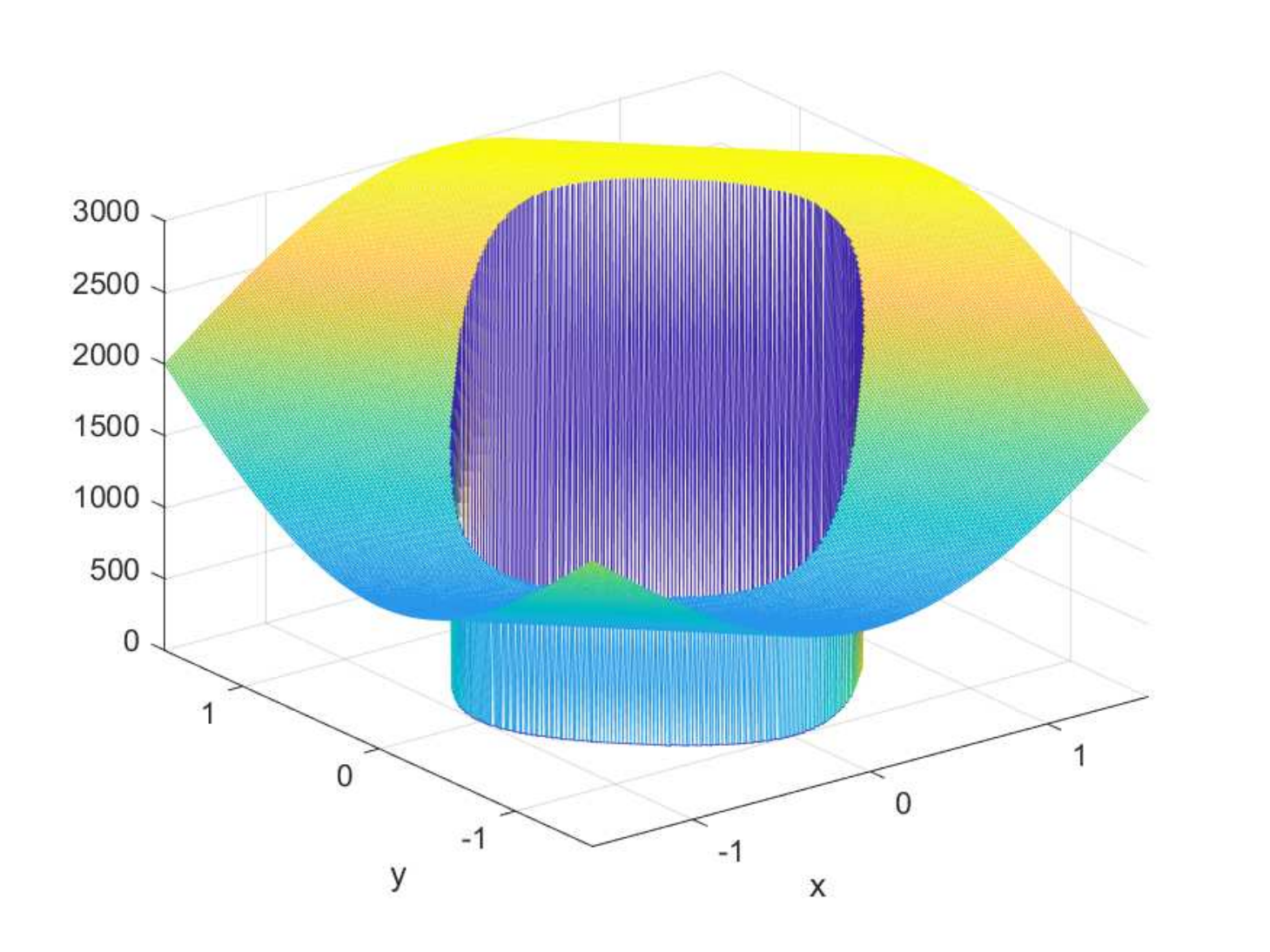}
	\end{subfigure}
	\begin{subfigure}[b]{0.3\textwidth}
		 \includegraphics[width=5.5cm,height=3.5cm]{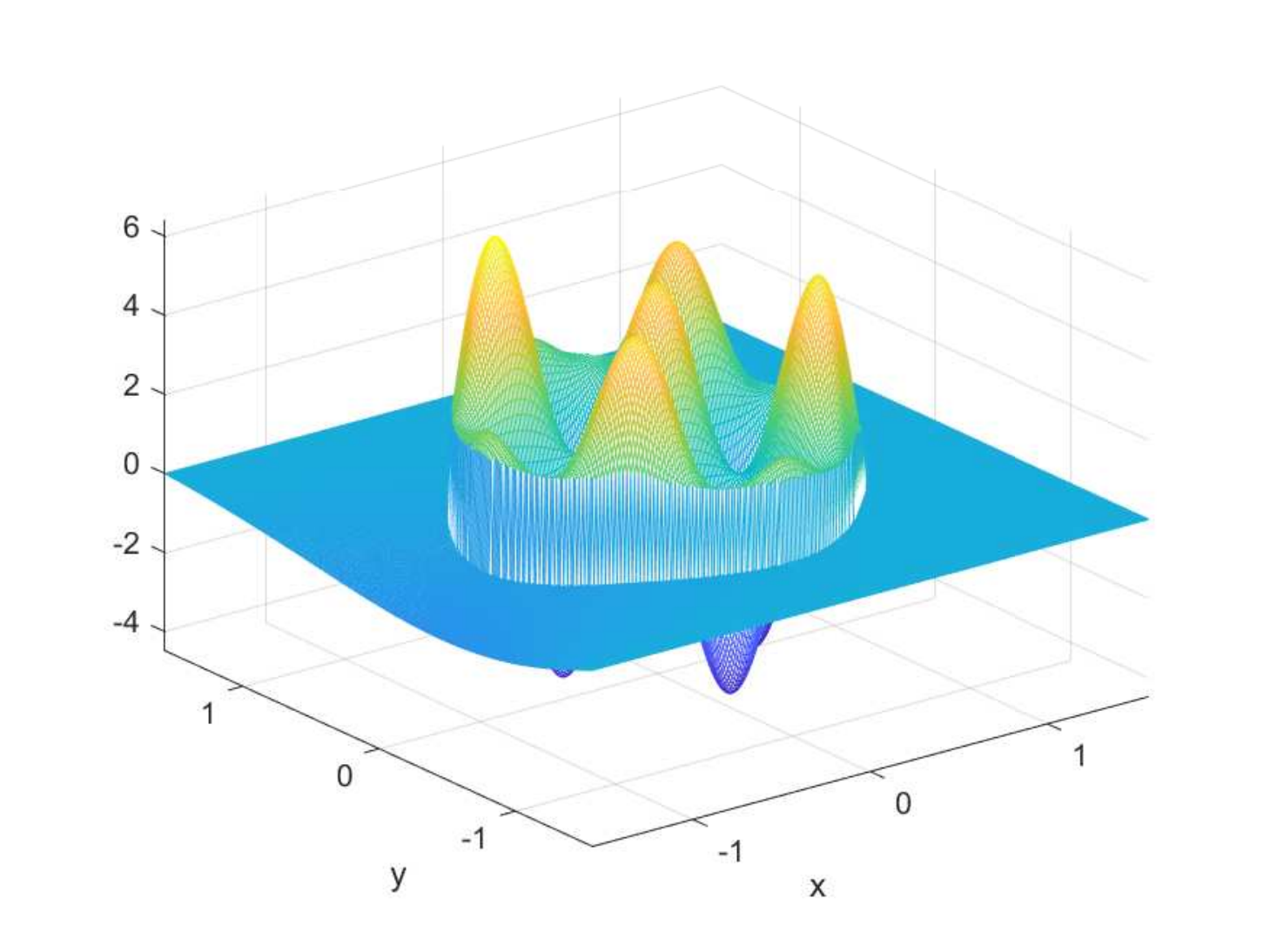}
	\end{subfigure}
	\caption
	{\cref{hybrid:unknown:ex3}: the interface curve $\Gamma_I$ (left), the coefficient $a(x,y)$ (middle) and the numerical solution $u_h$ (right) with $h=2^{-8}\times \pi$,  where $u_h$ is computed by our proposed hybrid finite difference scheme.}
	\label{hybrid:fig:unknown:QSp3}
\end{figure}

\begin{example}\label{hybrid:unknown:ex3}
	\normalfont
	Let $\Omega=(-\frac{\pi}{2},\frac{\pi}{2})^2$ and
	the interface curve be given by
	$\Gamma_I:=\{(x,y)\in \Omega:\; \psi(x,y)=0\}$ with
	$\psi (x,y)=y^2+\frac{2x^2}{x^2+1}-1$.
	The functions in \eqref{Qeques2} are given by
	\begin{align*}
		&a_{+}=10^{3}(2+\sin(x+y)),
		\qquad a_{-}=10^{-3}(2+\cos(x-y)), \qquad g_D=\sin(x)\cos(y)-2, \\
		&f_{+}=\sin(6 x)\sin(6 y),
		\qquad f_{-}=\cos(6 x)\cos(6 y), \qquad g_N=\cos(x+y),\\
		&  -u_x(-\frac{\pi}{2},y)+\cos(y) u(-\frac{\pi}{2},y)= \sin(y+\frac{\pi}{2})(y-\frac{\pi}{2}), \qquad
		\qquad u(\frac{\pi}{2},y)= 0, \qquad \mbox{for} \qquad y\in(-\frac{\pi}{2},\frac{\pi}{2}),\\
		& u(x,-\frac{\pi}{2})= 0, \qquad
		\qquad  u(x,\frac{\pi}{2})= 0,  \qquad \mbox{for} \qquad x\in(-\frac{\pi}{2},\frac{\pi}{2}).
	\end{align*}
The high contrast $a_+/a_-\approx 10^6$ on $\Gamma_I$.
	The numerical results are presented in \cref{hybrid:table:unknown:QSp3} and \cref{hybrid:fig:unknown:QSp3}.	
\end{example}

\begin{table}[htbp]
	\caption{Performance in \cref{hybrid:unknown:ex3}  of our proposed hybrid finite difference scheme on uniform Cartesian meshes with $h=2^{-J}\times \pi$.  }
	\centering
	\setlength{\tabcolsep}{1.5mm}{
		\begin{tabular}{c|c|c|c|c}
			\hline
			$J$
			&  ${\|u_{h}-u_{h/2}\|_{2}}$
			&order &  $\|u_{h}-u_{h/2}\|_{\infty}$

			&order  \\
			\hline
5   &1.17512E-01   &0   &1.95534E-01   &0   \\
6   &1.34603E-03   &6.4   &5.01334E-03   &5.3  \\
7   &2.97345E-05   &5.5   &9.62920E-05   &5.7   \\
8   &3.63705E-07   &6.4   &1.11523E-06   &6.4   \\
			\hline
	\end{tabular}}
	\label{hybrid:table:unknown:QSp3}
\end{table}

\begin{figure}[htbp]
	\centering
	\begin{subfigure}[b]{0.3\textwidth}
		 \includegraphics[width=5.5cm,height=3.5cm]{HyCUR4.pdf}
	\end{subfigure}
	\begin{subfigure}[b]{0.3\textwidth}
		 \includegraphics[width=5.5cm,height=3.5cm]{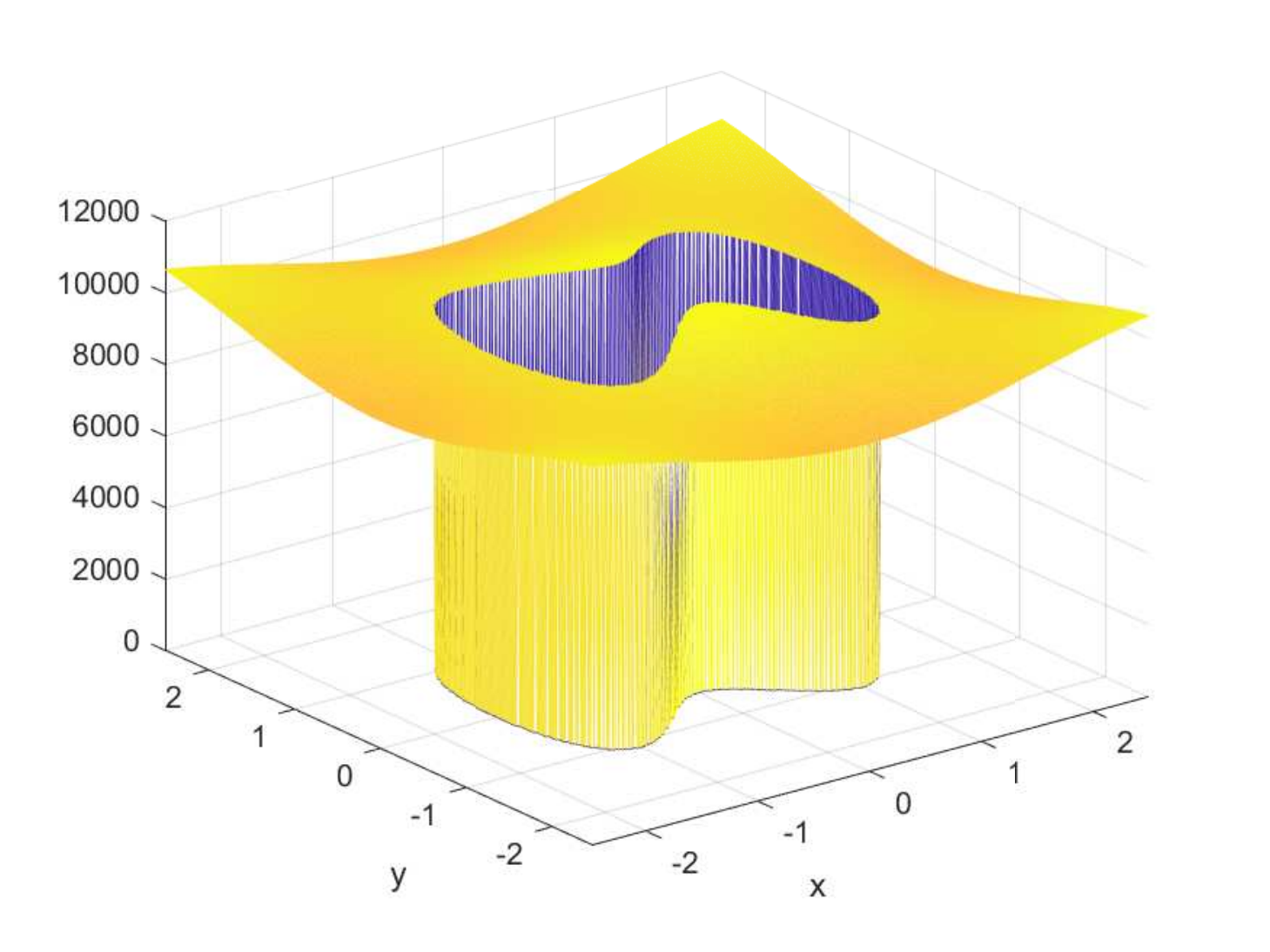}
	\end{subfigure}
	\begin{subfigure}[b]{0.3\textwidth}
		 \includegraphics[width=5.5cm,height=3.5cm]{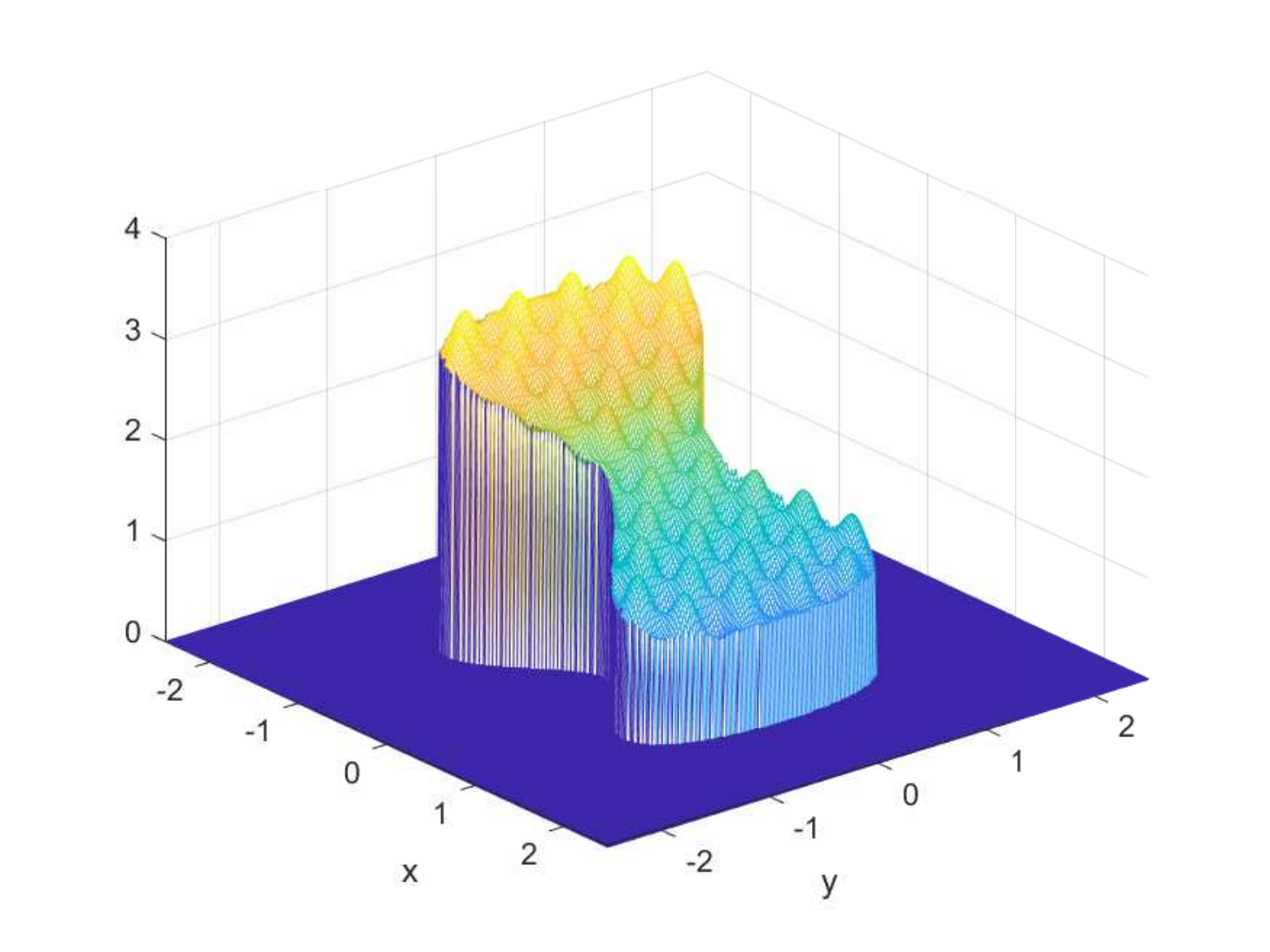}
	\end{subfigure}
	\caption
	{\cref{hybrid:unknown:ex4}: the interface curve $\Gamma_I$ (left), the coefficient $a(x,y)$ (middle) and the numerical solution $u_h$ (right) with $h=2^{-8}\times 5$. In order  to show the graph of $u_h$ clearly, we rotate the graph of $u_h$ by $\pi/2$ in this figure.}
	\label{hybrid:fig:unknown:QSp4}
\end{figure}

\begin{example}\label{hybrid:unknown:ex4}
	\normalfont
	Let $\Omega=(-2.5,2.5)^2$ and
	the interface curve be given by
	$\Gamma_I:=\{(x,y)\in \Omega:\; \psi(x,y)=0\}$ with
	$\psi (x,y)=y^2-2x^2+x^4-\frac{1}{4}$.
	The functions in \eqref{Qeques2} are given by
	\begin{align*}
		&a_{+}=10^{3}(10+\cos(x)\cos(y)),
		\qquad a_{-}=10^{-3}(10+\sin(x)\sin(y)), \qquad g_D=\sin(x)-2, \\
		&f_{+}=\sin(4\pi x)\sin(4\pi y),
		\qquad f_{-}=\cos(4\pi x)\cos(4\pi y), \qquad g_N=\cos(y),\\
		&  u(-2.5,y)= 0, \qquad
		\qquad u(2.5,y)= 0, \qquad \mbox{for} \qquad y\in(-2.5,2.5),\\
		& u(x,-2.5)= 0, \qquad
		\qquad  u(x,2.5)= 0,  \qquad \mbox{for} \qquad x\in(-2.5,2.5).
	\end{align*}
The high contrast $a_+/a_-\approx 10^6$ on $\Gamma_I$.	The numerical results are presented in \cref{hybrid:table:unknown:QSp4} and \cref{hybrid:fig:unknown:QSp4}.	
\end{example}

\begin{table}[htbp]
	\caption{Performance in \cref{hybrid:unknown:ex4}  of our proposed hybrid finite difference scheme on uniform Cartesian meshes with $h=2^{-J}\times 5$. }
	\centering
	\setlength{\tabcolsep}{1.5mm}{
		\begin{tabular}{c|c|c|c|c}
			\hline
			$J$
			&  ${\|u_{h}-u_{h/2}\|_{2}}$
			&order &  $\|u_{h}-u_{h/2}\|_{\infty}$

			&order  \\
			\hline
			5   &6.18678E+00   &0   &9.88338E+00   &0   \\
			6   &9.69535E-02   &6.0   &2.17089E-01   &5.5   \\
			7   &1.67043E-03   &5.9   &3.52407E-03   &5.9   \\
			8   &2.43148E-05   &6.1   &5.22530E-05   &6.1  \\
			\hline
	\end{tabular}}
	\label{hybrid:table:unknown:QSp4}
\end{table}

\begin{figure}[htbp]
	\centering
	\begin{subfigure}[b]{0.3\textwidth}
		 \includegraphics[width=5.5cm,height=3.5cm]{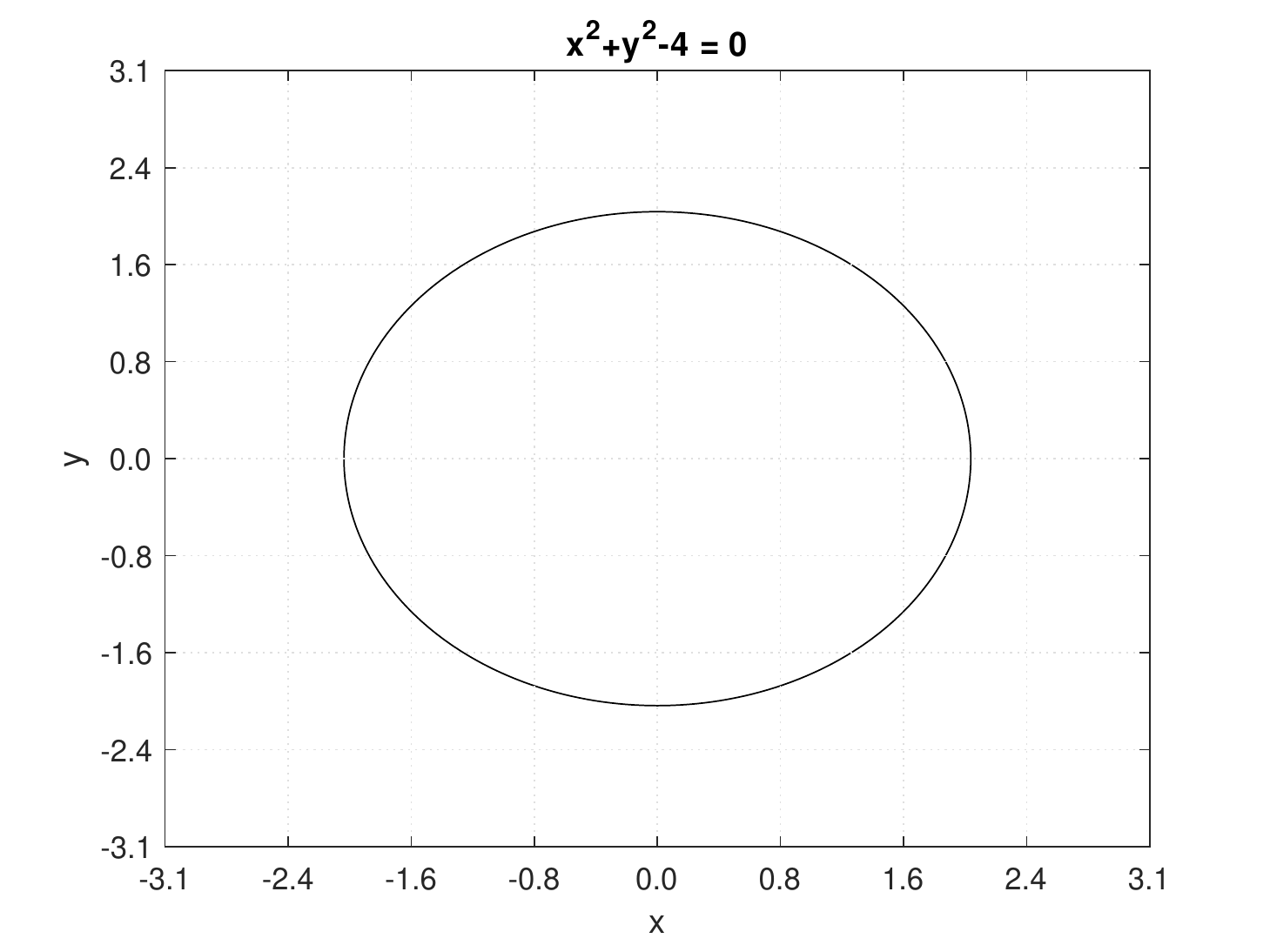}
	\end{subfigure}
	\begin{subfigure}[b]{0.3\textwidth}
		 \includegraphics[width=5.5cm,height=3.5cm]{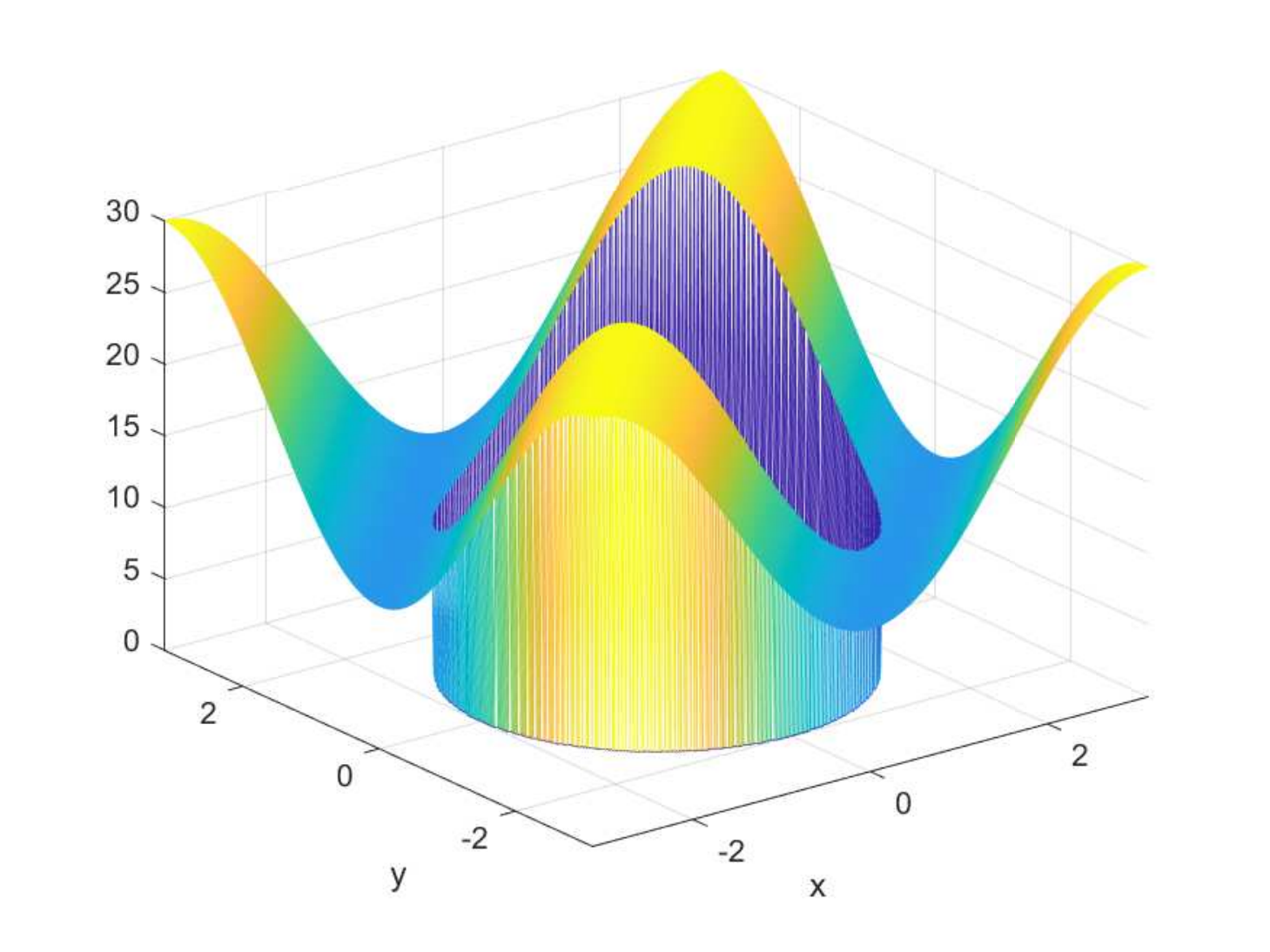}
	\end{subfigure}
	\begin{subfigure}[b]{0.3\textwidth}
		 \includegraphics[width=5.5cm,height=3.5cm]{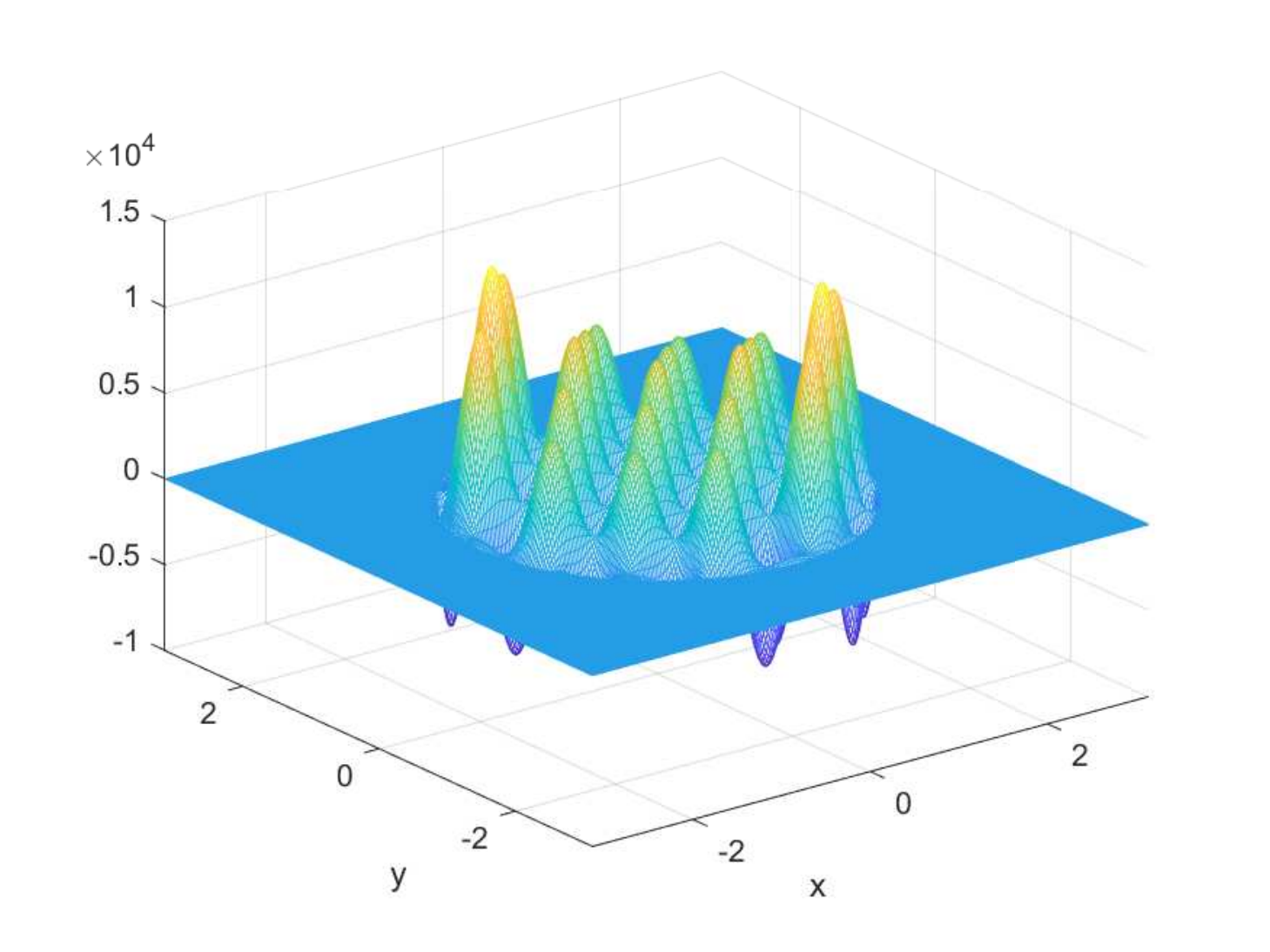}
	\end{subfigure}
	\caption
	{\cref{hybrid:unknown:ex5}: the interface curve $\Gamma_I$ (left), the coefficient $a(x,y)$ (middle) and the numerical solution $u_h$ (right) with $h=2^{-8}\times 2\pi$,  where $u_h$ is computed by our proposed hybrid finite difference scheme.}
	\label{hybrid:fig:unknown:QSp5}
\end{figure}

\begin{example}\label{hybrid:unknown:ex5}
	\normalfont
	Let $\Omega=(-\pi,\pi)^2$ and
	the interface curve be given by
	$\Gamma_I:=\{(x,y)\in \Omega:\; \psi(x,y)=0\}$ with
	$\psi (x,y)=x^2+y^2-4$.
	The functions in \eqref{Qeques2} are given by
	\begin{align*}
		&a_{+}=10(2+\cos(x-y)),
		\qquad a_{-}=10^{-6}(2+\sin(x)\sin(y)), \qquad g_D=\sin(y)-10, \\
		&f_{+}=\sin(6 x)\sin(6 y),
		\qquad f_{-}=\cos(6 x)\cos(6 y), \qquad g_N=\cos(x),\\
		& -u_x(-\pi,y)+\sin(y) u(-\pi,y)= \cos(y), \qquad
		\qquad u(\pi,y)=0, \qquad \mbox{for} \qquad y\in(-\pi,\pi),\\
		& -u_y(x,-\pi)= \sin(x-\pi), \qquad
		\qquad u_y(x,\pi)+\cos(x) u(x,\pi)=\cos(x)+1, \qquad \mbox{for} \qquad x\in(-\pi,\pi).
	\end{align*}
The high contrast $a_+/a_-\approx 10^7$ on $\Gamma_I$.		The numerical results are presented in \cref{hybrid:table:unknown:QSp5} and \cref{hybrid:fig:unknown:QSp5}.	
\end{example}

\begin{table}[htbp]
	\caption{Performance in \cref{hybrid:unknown:ex5}  of our proposed hybrid finite difference scheme on uniform Cartesian meshes with $h=2^{-J}\times 2\pi$.}
	\centering
	\setlength{\tabcolsep}{1.5mm}{
		\begin{tabular}{c|c|c|c|c}
			\hline
			$J$
			&  ${\|u_{h}-u_{h/2}\|_{2}}$
			&order &  $\|u_{h}-u_{h/2}\|_{\infty}$

			&order  \\
			\hline
5   &1.60217E+04   &0   &1.39059E+04   &0   \\
6   &2.94197E+02   &5.8   &2.79828E+02   &5.6   \\
7   &4.54676E+00   &6.0   &6.36193E+00   &5.5   \\
8   &5.82759E-02   &6.3   &1.02577E-01   &6.0   \\
			\hline
	\end{tabular}}
	\label{hybrid:table:unknown:QSp5}
\end{table}

%
%
%
%

\section{Conclusion}\label{sec:hybrid:Conclu}

To our best knowledge, so far there were no 13-point finite difference schemes for irregular points available in the literature, that can achieve fifth or sixth order
for elliptic interface problems with discontinuous coefficients.
Our contributions of this paper are as follows:
\begin{itemize}	
	\item We propose a hybrid (13-point for irregular points and compact 9-point for interior regular points) finite difference scheme, which demonstrates six order accuracy in all our numerical experiments,
for elliptic interface problems with  discontinuous, variable and high-contrast coefficients, discontinuous source terms and two non-homogeneous jump conditions.
	
	\item  The proposed hybrid scheme demonstrates a robust high-order convergence for the challenging cases of high-contrast ratios of the coefficients $a_{\pm}$: $\sup(a_+)/\inf(a_-)=10^{-3},10^{-6},10^{6},10^{7}$.

	\item Due to the flexibility and efficiency of the implementation, it is very simple to achieve the implementation for 25-point or 36-point schemes for irregular points of elliptic interface problems and Helmholtz interface equations with discontinuous wave numbers.


\item From the results in \cref{hybrid:table:QSp3,hybrid:table:QSp5}, we find that if we only replace the $13$-point scheme for irregular points by a $9$-point scheme, then the numerical errors increase significantly, while the condition number only slightly decreases. Thus, the proposed hybrid scheme could significantly improve the numerical performance with a slight increase in  the complexity of the corresponding linear system.

	\item  We  also derive a $6$-point/$4$-point schemes with a sixth order accuracy at the side/corner points for the case of smooth coefficients $\alpha$ and $\beta$ in the Robin boundary conditions   $\tfrac{\partial u}{\partial \nv}+\alpha u=g_1$ and $\tfrac{\partial u}{\partial \nv}+\beta u=g_4$.
	
	\item
	The presented numerical experiments confirm the sixth order of accuracy in the $l_2$ and $l_{\infty}$ norms of  our proposed hybrid scheme.
\end{itemize}
\section{Appendix}
\label{hybrid:sec:proofs}

Let us first present the definitions of several index sets $\ind_{M+1}, \ind_{M+1}^{V, 1}, \ind_{M+1}^{V, 2}, \ind_{M+1}^{H, 1}, \ind_{M+1}^{H, 2}$, which are employed in \cref{sec:sixord}.
	Define $\NN:=\N\cup\{0\}$, the set of all nonnegative integers.
Given $M+1\in \NN$, we use the same definitions  \cite[(2.4) and (2.7)]{FHM21Helmholtz} as follows:
\be \label{Sk}
\ind_{M+1}:=\{(m,n-m) \; : \; n=0,\ldots,M+1
\; \mbox{ and }\; m=0,\ldots,n\}, \qquad M+1\in \NN,
\ee
\be \label{indV12}
\ind_{M+1}^{V, 2}:=\ind_{M+1}\setminus \ind_{M+1}^{V, 1}\quad \mbox{with}\quad
\ind_{M+1}^{V, 1}:=\{(\ell,k-\ell) \; :   k=\ell,\ldots, M+1-\ell\;\; \mbox{and} \;\;\ell=0,1\; \},
\ee
\be \label{indH12}
\ind_{M+1}^{H, j}:=\{(n,m):(m,n) \in \ind_{M+1}^{V,j}, j =1,2\}.
\ee
%

For all $(m,n)\in \ind_{M+1}^{V,1}$, we define
\be\label{GVmn}
G^V_{M,m,n}(x,y):=\sum_{\ell=0}^{\lfloor \frac{n}{2}\rfloor}
\frac{(-1)^\ell x^{m+2\ell} y^{n-2\ell}}{(m+2\ell)!(n-2\ell)!}+\sum_{(m',n')\in \ind_{M+1}^{V,2} \setminus \ind_{m+n}^{V,2} }A^{V,u}_{m',n',m,n} \frac{ x^{m'} y^{n'}}{m'!n'!},
\ee
and for all $(m,n)\in \ind_{M-1}$,
\be\label{QVmn}
\begin{split}
	Q^V_{M,m,n}(x,y):=\sum_{\ell=1}^{1+\lfloor \frac{n}{2}\rfloor} \frac{(-1)^{\ell} x^{m+2\ell} y^{n-2\ell+2}}{(m+2\ell)!(n-2\ell+2)!}\frac{1}{a^{(0,0)}}
	+\sum_{(m',n')\in \ind_{M+1}^{V,2} \setminus \ind_{m+n+2}^{V,2} }A^{V,f}_{m',n',m,n} \frac{ x^{m'} y^{n'}}{m'!n'!},
\end{split}
\ee
where $A^{V,u}_{m',n',m,n}$ and $A^{V,f}_{m',n',m,n}$ are constants which are uniquely determined by $\{a^{(m,n)}: (m,n) \in \ind_{M}\}$,  and the floor function $\lfloor x\rfloor$ is defined to be the largest integer less than or equal to $x\in \R$.

For all $(m,n)\in \ind_{M+1}^{H,1}$, we define
\be\label{GHmn}
G^H_{M,m,n}(x,y):=\sum_{\ell=0}^{\lfloor \frac{m}{2}\rfloor}
\frac{(-1)^\ell y^{n+2\ell} x^{m-2\ell}}{(n+2\ell)!(m-2\ell)!}+\sum_{(m',n')\in \ind_{M+1}^{H,2} \setminus \ind_{m+n}^{H,2} }A^{H,u}_{m',n',m,n} \frac{ x^{m'} y^{n'}}{m'!n'!},
\ee
and for all $(m,n)\in \ind_{M-1}$,
\be\label{QHmn}
\begin{split}
	Q^H_{M,m,n}(x,y):=\sum_{\ell=1}^{1+\lfloor \frac{m}{2}\rfloor} \frac{(-1)^{\ell} y^{n+2\ell} x^{m-2\ell+2}}{(n+2\ell)!(m-2\ell+2)!}\frac{1}{a^{(0,0)}}
	+\sum_{(m',n')\in \ind_{M+1}^{H,2} \setminus \ind_{m+n+2}^{H,2} }A^{H,f}_{m',n',m,n} \frac{ x^{m'} y^{n'}}{m'!n'!},
\end{split}
\ee
where $A^{H,u}_{m',n',m,n}$ and $A^{H,f}_{m',n',m,n}$ are constants which are uniquely determined by $\{a^{(m,n)}: (m,n) \in \ind_{M}\}$,  and the floor function $\lfloor x\rfloor$ is defined to be the largest integer less than or equal to $x\in \R$.


In this appendix, we provide the  proofs to all the technical results stated in \cref{sec:sixord}.

\begin{proof}[Proof of \cref{thm:regular:interior}]
Choose $M=6$ and replace  $G_{m,n}$, $H_{m,n}$ and $\ind_{M+1}^{1}$ in \cite[]{FHM21b} by $G^V_{M,m,n}$ given in  \eqref{GVmn}, $Q^V_{M,m,n}$ in \eqref{QVmn}, and  $\ind_{M+1}^{V,1}$ in  \eqref{indV12} .
\end{proof}	


\begin{proof}[Proof of \cref{hybrid:thm:regular:Robin:1}]
	Let $M_f=M_{g_1}=M$ in  the proof of \cite[Theorem 2.3]{FHM21Helmholtz}. Then  \cite[(4.7)]{FHM21Helmholtz} implies
	\be \label{stencil:regular:EQ:V：Robin:1}
	\sum_{k=0}^1 \sum_{\ell=-1}^1
	C^{\mathcal{B}_1}_{k,\ell} u(x_i+kh,y_j+\ell h)=
	\sum_{(m,n)\in \ind_{M-1}} f^{(m,n)}C^{\mathcal{B}_1}_{f,m,n}+\sum_{n=0}^{M} g_1^{(n)}C^{\mathcal{B}_1}_{g_1,n}+\bo(h^{M+2}), \qquad h \to 0
	\ee
	Since $-u_x+\alpha u=g_1$ on $\Gamma_1$, we have $u^{(1,n)} = \sum_{i=0}^n  {n\choose i}  {\alpha}^{(n-i)}u^{(0,i)} - g_1^{(n)}$ for all $n = 0,\dots, M$. By \eqref{u:approx:key:V},
	{\footnotesize{
			\begin{align*}
				&u(x+x_i^*,y+y_j^*)=
				\sum_{n=0}^{M+1}
				u^{(0,n)} G^{V}_{M,0,n}(x,y)+\sum_{n=0}^{M}
				u^{(1,n)} G^{V}_{M,1,n}(x,y) +\sum_{(m,n)\in \ind_{M-1}}
				f^{(m,n)} Q^{V}_{M,m,n}(x,y) +\bo(h^{M+2}) \\
				&=
				\sum_{n=0}^{M+1}
				u^{(0,n)} G^{V}_{M,0,n}(x,y)+\sum_{n=0}^{M}
				u^{(1,n)} G^{V}_{M,1,n}(x,y) +\sum_{(m,n)\in \ind_{M-1}}
				f^{(m,n)} Q^{V}_{M,m,n}(x,y) +\bo(h^{M+2}) \\
				& = \sum_{n=0}^{M+1}
				u^{(0,n)} G^{V}_{M,0,n}(x,y)+\sum_{n=0}^{M}
				\bigg( \sum_{i=0}^n  {n\choose i}  {\alpha}^{(n-i)}u^{(0,i)} -g_1^{(n)}  \bigg) G^{V}_{M,1,n}(x,y) +\sum_{(m,n)\in \ind_{M-1}}
				f^{(m,n)} Q^{V}_{M,m,n}(x,y)\\
				& \qquad +\bo(h^{M+2})\\
				&=\sum_{n=0}^{M+1}
				u^{(0,n)} G^{V}_{M,0,n}(x,y) +\sum_{n=0}^{M}
				\sum_{i=0}^n  {n\choose i}  {\alpha}^{(n-i)}u^{(0,i)}  G^{V}_{M,1,n}(x,y) -\sum_{n=0}^{M}
				g_{1}^{(n)}  G^{V}_{M,1,n}(x,y) \\
				& \qquad +\sum_{(m,n)\in \ind_{M-1}}
				f^{(m,n)} Q^{V}_{M,m,n}(x,y) +\bo(h^{M+2})\\
				&=\sum_{n=0}^{M+1}
				u^{(0,n)} G^{V}_{M,0,n}(x,y) +\sum_{i=0}^{M}
				\sum_{n=i}^M  {n\choose i}  {\alpha}^{(n-i)}u^{(0,i)} G^{V}_{M,1,n}(x,y)  -\sum_{n=0}^{M}
				g_{1}^{(n)}  G^{V}_{M,1,n}(x,y)\\
				& \qquad +\sum_{(m,n)\in \ind_{M-1}}
				f^{(m,n)} Q^{V}_{M,m,n}(x,y) +\bo(h^{M+2})\\
				&=u^{(0,M+1)} G^{V}_{M,0,M+1}(x,y)+\sum_{n=0}^{M}
				u^{(0,n)} G^{V}_{M,0,n}(x,y) +\sum_{n=0}^{M}
				\sum_{i=n}^M  {i\choose n}  {\alpha}^{(i-n)}u^{(0,n)} G^{V}_{M,1,i}(x,y)  -\sum_{n=0}^{M}
				g_{1}^{(n)}  G^{V}_{M,1,n}(x,y)\\
				& \qquad +\sum_{(m,n)\in \ind_{M-1}}
				f^{(m,n)} Q^{V}_{M,m,n}(x,y) +\bo(h^{M+2}), \quad \mbox{for } x,y\in (-2h,2h).
			\end{align*}
		}
	}
	So \eqref{stencil:regular:EQ:V：Robin:1} leads to
	\be \label{stencil:regular:EQ2:V：Robin:1}
	\sum_{n=0}^{M+1}
	u^{(0,n)} I^{\mathcal{B}_1}_{n}+
	\sum_{(m,n)\in \ind_{M-1}} f^{(m,n)}
	\left(J^{\B_1}_{m,n}-C^{\B_1}_{f,m,n}\right) +\sum_{n=0}^{M} g_1^{(n)}\left(K^{\B_1}_{n}-C^{\B_1}_{g_{1},n}\right)
	=\bo(h^{M+2}),
	\ee
	as $h \to 0$, where
	\begin{align}
		\nonumber
		&I^{\B_1}_{n}:=\sum_{k=0}^1 \sum_{\ell=-1}^1 C^{\B_1}_{k,\ell} \left( G^{V}_{M,0,n}(kh, \ell h) +
		\sum_{i=n}^M  {i\choose n}  {\alpha}^{(i-n)} G^{V}_{M,1,i}(kh, \ell h) (1-\delta_{n,M+1}) \right),\\
		\label{IB1n}
		& J^{\B_1}_{m,n}:=\sum_{k=0}^1 \sum_{\ell=-1}^1 C^{\B_1}_{k,\ell} Q^{V}_{M,m,n}(kh, \ell h), \quad
		K^{\B_1}_{n}:=-\sum_{k=0}^1 \sum_{\ell=-1}^1 C^{\B_1}_{k,\ell} G^{V}_{M,1,n}(kh, \ell h),
	\end{align}
	$\delta_{a,a}=1$, and $\delta_{a,b}=0$ for $a \neq b$.
\end{proof}
\begin{proof}[Proof of \cref{hybrid:thm:regular:Neu:3}]
The proof is almost identical to the proof of \cref{hybrid:thm:regular:Robin:1}.
\end{proof}	
\begin{proof}[Proof of \cref{hybrid:thm:regular:Robin:4}]
The proof is almost identical to the proof of \cref{hybrid:thm:regular:Robin:1}.
\end{proof}	
\begin{proof}[Proof of \cref{thm:corner:1}]
The proof is similar to the proof of  \cite[Theorem 2.4]{FHM21Helmholtz}. Precisely, replace $\B_1u=\frac{\partial u}{\partial \nv}- \ia \ka u=g_1$ by $\B_1u=\frac{\partial u}{\partial \nv}+\alpha u=g_1$
	in the proof  of  \cite[Theorem 2.4]{FHM21Helmholtz} with $M=M_f=M_{g_1}=M_{g_3}=5$,  and replace  \cite[$G^{V}_{M,m,n}$, $Q^{V}_{M,m,n}$, $G^{H}_{M,m,n}$ and $Q^{H}_{M,m,n}$]{FHM21Helmholtz} by \eqref{GVmn}, \eqref{QVmn}, \eqref{GHmn} and \eqref{QHmn}.
\end{proof}	
\begin{proof}[Proof of \cref{thm:corner:2}]
The proof is similar to the proof of  \cite[Theorem 2.5]{FHM21Helmholtz}. Precisely, replace $\B_1u=\frac{\partial u}{\partial \nv}- \ia \ka u=g_1$ and $\B_4u=\frac{\partial u}{\partial \nv}-\ia \ka u=g_4$ by $\B_1u=\frac{\partial u}{\partial \nv}+\alpha u=g_1$ and $\B_4u=\frac{\partial u}{\partial \nv}+\beta u=g_4$ respectively
	in the proof of \cite[Theorem 2.5]{FHM21Helmholtz} with $M=M_f=M_{g_1}=M_{g_4}=5$  and replace  \cite[$G^{V}_{M,m,n}$, $Q^{V}_{M,m,n}$, $G^{H}_{M,m,n}$ and $Q^{H}_{M,m,n}$]{FHM21Helmholtz} by \eqref{GVmn}, \eqref{QVmn}, \eqref{GHmn} and \eqref{QHmn}. 	
\end{proof}	
\begin{proof}[Proof of \cref{thm:interface}]
	\eqref{T0000} can be obtained by $u_-^{(0,0)}=u_+^{(0,0)}-g_{D}^{(0,0)}$ and
	  \cite[(7.18)]{FHM21b}. The rest of the proof is straightforward and follows from \cite[(7.8), (7.10), (7.16), and (7.18)]{FHM21b}.
\end{proof}
\begin{proof}[Proof of \cref{13point:Inter}]
	Choose $M=4$, replace $\ind_{M+1}^{1}$, $G^{\pm}_{m,n}$, $H^{\pm}_{m,n}$ $d_{i,j}^\pm $ in  \cite[Theorem 3.3]{FHM21b} by $\ind_{M+1}^{V,1}$,  $G^{\pm,V}_{M,m,n}$, $Q^{\pm,V}_{M,m,n}$, $d_{i,j}^\pm \cup e_{i,j}^\pm$  in this paper.
\end{proof}	
\begin{proof}[Proof of \eqref{ckln:irregular}]
	Note that when we use the formulas of \cite{FHM21b} in this proof, we need to replace $\ind_{M+1}^{1}$, $G^{\pm}_{m,n}$, $H^{\pm}_{m,n}$ $d_{i,j}^\pm $ in \cite{FHM21b} by $\ind_{M+1}^{V,1}$,  $G^{\pm,V}_{M,m,n}$, $Q^{\pm,V}_{M,m,n}$, $d_{i,j}^\pm \cup e_{i,j}^\pm$  in this paper.		Consider $I_{0,0}(h)=\bo(h^{M+2})$  in \cite[(3.29)]{FHM21b}.
	According to  \cite[(3.28)]{FHM21b} and \eqref{newJuT} in this paper, $I_{0,0}(h)=\bo(h^{M+2})$ implies
	\be \label{U00:coeff}
	\sum_{(k,\ell)\in d_{i,j}^+\cup e_{i,j}^+}
	C_{k,\ell}(h) G^{+,V}_{M,0,0}(v_0h+kh,w_0h+\ell h)+	\sum_{ \substack{ (m',n')\in \ind_{M+1}^{V,1} \\  m'+n' \ge 0}} I^-_{m',n'}(h) T^{u_+}_{m',n',0,0}=\bo(h^{M+2}).
	\ee
	By \eqref{T0000},
	\eqref{U00:coeff} is equivalent to
	\[\sum_{(k,\ell)\in d_{i,j}^+\cup e_{i,j}^+}
	C_{k,\ell}(h) G^{+,V}_{M,0,0}(v_0h+kh,w_0h+\ell h)+	 I^-_{0,0}(h)=\bo(h^{M+2}),\]
	i.e.,
	\be \label{U00:coeff:2}
	\sum_{(k,\ell)\in d_{i,j}^+\cup e_{i,j}^+}
	C_{k,\ell}(h) G^{+,V}_{M,0,0}(v_0h+kh,w_0h+\ell h)+	 \sum_{(k,\ell)\in d_{i,j}^-\cup e_{i,j}^-}
	C_{k,\ell}(h) G^{-,V}_{M,0,0}(v_0h+kh,w_0h+\ell h)=\bo(h^{M+2}).
	\ee		
	According to the proof of  \cite[Lemma 2.1]{FHM21b} and \eqref{GVmn},
	\be\label{G00}
	G^{\pm,V}_{M,0,0}(x,y):=1.
	\ee
	Consider the coefficients of $h^i$ for $i=0,1,\dots, M+1$ in \eqref{U00:coeff:2}, then \eqref{G00} implies
	\be\label{ckl0}
	\sum_{(k,\ell)\in d_{i,j}^+\cup e_{i,j}^+}
	c_{k,\ell,i} +	 \sum_{(k,\ell)\in d_{i,j}^-\cup e_{i,j}^-}
	c_{k,\ell,i} =0,	\quad \mbox{for} \quad i=0,1,\dots, M+1.
	\ee

This proves \eqref{ckln:irregular}.
\end{proof}

\end{document}